\documentclass{article}

\usepackage{arxiv}

\usepackage{amsmath}
\usepackage{graphicx}
\usepackage{epstopdf}
\usepackage[utf8]{inputenc} % allow utf-8 input
\usepackage[T1]{fontenc}    % use 8-bit T1 fonts
\usepackage{hyperref}       % hyperlinks
\usepackage{url}            % simple URL typesetting
\usepackage{booktabs}       % professional-quality tables
\usepackage{amsfonts}       % blackboard math symbols
\usepackage{nicefrac}       % compact symbols for 1/2, etc.
\usepackage{microtype}      % microtypography
\usepackage{lipsum}
\usepackage{algorithm}
\usepackage{algorithmic}
\usepackage{booktabs}
\usepackage{lscape}
\usepackage{caption}
\usepackage{subcaption}
\usepackage{multirow}
\newtheorem{thm}{Theorem}[section]
\newtheorem{lem}[thm]{Lemma}

\newtheorem{defi}{Definition}[section]

\newcommand{\goto}{\ensuremath{\rightarrow}}

\newcommand{\T}{\ensuremath{\mathcal{T}}}

\newenvironment{Assumptions}
{
\setcounter{enumi}{0}

\begin{enumerate}}
{\end{enumerate} }

\numberwithin{equation}{section} \allowdisplaybreaks

\title{Well-posedness study of a non-linear hyperbolic-parabolic coupled system applied to image speckle reduction}

\author{
 Sudeb Majee \\
  School of Basic Sciences\\ 
  Indian Institute of Technology Mandi\\
  PIN 175005, INDIA\\
  \texttt{sudebmajee@gmail.com} \\
    \And
  Rajendra K. Ray \\
  School of Basic Sciences\\ 
  Indian Institute of Technology Mandi\\
  PIN 175005, INDIA\\
  \texttt{rajendra@iitmandi.ac.in} \\
    \And
   Ananta K. Majee \\
   Department of Mathematics\\
   Indian Institute of Technology Delhi\\
   PIN 110016, INDIA \\
  \texttt{majee@maths.iitd.ac.in} \\ 
}

\begin{document}
\maketitle

\begin{abstract}
In this article, we consider a  non-linear hyperbolic-parabolic coupled system based on telegraph diffusion framework applied to image despeckling. A separate equation is used to calculate the edge variable, which improves the quality of the despeckled images. A well-posedness result of the proposed coupled system is settled via Schauder's fixed point theorem. Numerical experiments are reported to illustrate the effectiveness of the proposed model, with recently developed models, over a set of gray level test images contaminated by speckle noise. 
\end{abstract}


\section{Introduction}
\label{intro}
Beginning with the Perona-Malik model \cite{perona1990scale}, non-linear partial differential equations~(PDEs) are extensively used to develop noise reduction models. Due to the availability of well established numerical schemes and theoretical properties, PDE based image processing is an exciting research area for real-life application purpose as well as for the theoretical study. In the real situation, images are often degraded by different types of noises, e.g., additive, multiplicative, or mixed nature. Hence the noise extraction is a very initial stage for high-level image analysis. In this work, we only consider the multiplicative speckle \cite{majeetdm2019speckle} noise removal process. A Mathematical representation for a degraded image affected by speckle noise\cite{dutt1995statistical}  can be expressed as 
\begin{equation*}
J=I\eta,
\end{equation*}
where $J$ is the noisy image, $I$ is the noise-free image, and $\eta$  signifies the speckle-noise process.

In general, the speckle noise process $\eta$ is Gamma$(L,L)$ distributed, where $L \in {\rm I\!N} $ is the the number of looks corresponding to the noise level in the corrupted
images \cite{argenti2013tutorial,hao2015variational,liu2016modified}. To remove speckle based noise in the images, different types of PDE based models are proposed and resulted in significant momentum both in the development of theoretical as well as numerical aspects of the problems.
Most popular PDE based approaches are
anisotropic diffusion-based methods \cite{jain2019non,jain2018nonlinear,jin2000adaptive,shan2019multiplicative,yu2002speckle,zhou2015doubly,zhou2018nonlinear},
and variational methods \cite{  aubert2008variational,dong2013convex,huang2010multiplicative,jidesh2013complex, jin2010analysis, jin2011variational, liu2013nondivergence,rudin2003multiplicative,shi2008nonlinear}.
Most of the above models take the generalized form
\begin{align}\label{eq:general1}
&I_t  = \text{div}(g(t,x)\nabla I) + \lambda f(J,I)  \hspace{0.3cm} \text{in} \hspace{0.2cm} \Omega_T:=(0,T) \times \Omega,
\end{align}
with the appropriate initial and boundary conditions.
Here $\Omega$ is the domain of the original image $I$ and the observed noisy image $J,$ $T$ is a specified time, $\lambda$ is the weight parameter, $\text{div}$ and $\nabla$ represents the divergence and
gradient operator respectively. The source term $f(J,I)$ is derived from the variational model approach \cite{aubert2008variational,dong2013convex, jin2011variational,rudin2003multiplicative}.
In \eqref{eq:general1}, $g(t,x)$ signifies the degree of denoising which preserves the image characteristics, e.g., textures and edges in the noise removal process. All the above-discussed PDE based models are
parabolic type.
\vspace{.1cm}

Later, V. Ratner, and Y. Zeevi \cite{ratner2007image} introduces the idea of hyperbolic PDE for additive noise removal process. By considering the image as an elastic sheet, the authors in \cite{ratner2007image} suggest the following telegraph diffusion equation (TDE) based model
\begin{align*}
&I_{tt} +\gamma I_{t}- \text{div} \left(   \dfrac{\nabla I}{1+\vert\frac{\nabla I}{K}\vert^2}   \right)=0\,, \hspace{0.2cm} \text{in}\,\,\, \Omega_T\,,
\end{align*}
where $\gamma$ is the damping parameter and $K$ is a threshold constant. Even though the TDE model can effectively preserve the sharp edges but failed to produce satisfactory smoothing in the presence of a large level of noise. To overcome this issue, several non-linear telegraph diffusion-based method have been proposed \cite{cao2010class,jain2016edge,sun2016class,yang2014kernel,zhang2015spatial}.
However, in spite of their impressive applications in the field of additive noise removal process, hyperbolic PDE based approaches have not successfully used for speckle noise removal process. Recently Sudeb
et al. suggest a couple of hyperbolic PDE based models \cite{fuzzy2019ttvmodel,majeetdm2019speckle} for speckle noise removal process. 
The authors in \cite{fuzzy2019ttvmodel} developed a model in a telegraph total variation framework as
\begin{align*}
&I_{tt}+\gamma I_t = \text{div}\left(\theta(I)\frac{\nabla I}{|\nabla I|}\right)-\lambda \left( 1-\frac{J}{I}\right),
\hspace{0.6cm} \text{in} \hspace{0.2cm} \Omega_{T},  \\
&\partial_n I =0,               \hspace{5.8cm} \text{in}
\hspace{0.2cm} \partial \Omega_T:=(0,T)\times \partial \Omega,\\
&I(0,x)=J(x), \hspace{0.1cm}   I_t(0,x)=0,	\hspace{2.8cm} \text{in} \hspace{0.2cm}\Omega,
\end{align*}
where $\theta$ is the fuzzy edge detector function \cite{chaira2008new}. In \cite{majeetdm2019speckle}, the authors developed a model in a telegraph diffusion based framework of the form
\begin{align*}
I_{tt} +\gamma I_{t}- \text{div} \left( g\left(I_\xi,\vert \nabla I_\xi \vert \right)    \nabla I\right)=0\,, \hspace{1.1cm} \text{in}\,\,\, \Omega_T\,,
\end{align*}
with the similar initial and boundary conditions as in \cite{fuzzy2019ttvmodel}, where the diffusion control function $g$ is given by
\begin{align*}
g\left(I_\xi,\vert \nabla I_\xi \vert \right)=\dfrac{ 2\vert I_\xi \vert^\nu}{\big(M^{I}_{\xi}\big)^\nu+\vert I_\xi \vert^\nu}\cdot \dfrac{1}{1+ \left(\frac{|\nabla I_{\xi}|}{K} \right)^2 }\,.
\end{align*}
In the above, $I_\xi=G_\xi\ast I$, $M^{I}_{\xi}= \underset{ x \in \Omega}{\text{max}}\vert I_\xi(t,x) \vert,$ $\nu \geq 1,$ $\xi>0,$ and $``\ast"$ represents convolution in $x$ only and $G_{\xi}$ is a two dimensional Gaussian kernel.
\vspace{.1cm}

To the best of our knowledge, most of the PDE based models for speckle noise removal are single and parabolic types. Inspired by the ideas of \cite{jain2018nonlinear,majeetdm2019speckle}, we propose the following improved nonlinear and coupled hyperbolic-parabolic model
\begin{align}
&I_{tt} + \gamma I_{t}- {\rm div}\Big( \frac{s^\alpha}{ 1 + s^\alpha} \frac{1}{ 1 + \iota|u_{\xi}|^\beta} \nabla I\Big)=0, \hspace{2.5cm} {\rm in}~~~ \Omega_T\,, \label{maina} \\
& u_t = h(|\nabla I_\xi|) - u + \frac{\nu^2}{2}\Delta u,  \hspace{4.9cm} \text{in}~~~\Omega_T\,, \label{mainb} \\
& \begin{cases} \label{mainc}
 I(0,x)= I_0(x)\,,~ I_t(0,x)=0\,,~ u(0,x)= G_\xi \ast |\nabla I_0|^2 \hspace{0.9cm} {\rm in}~~~\Omega\,, \\
\partial_n I =0=\partial_n u \hspace{6.4cm} {\rm on}~~~\partial \Omega_T\,,
\end{cases}
\end{align}
with $s:=\frac{|I_\xi|}{M_\xi^I} \in [0,1]$. 
In the above, $I_0$ is the observed noisy image, $ \alpha \geq 1, \beta \geq 1, \gamma>0,\nu>0,$ and $\iota>0$ are constants.  $\Delta$ is Laplace operator, $h:\mathbb{R}^+ \rightarrow \mathbb{R}^+$ is a bounded,
Lipschitz continuous function. Moreover, $u$ represents the edge strength at each scale. Here we utilize an extra equation to calculate the edge variable, which improves the present model over our previous
model \cite{majeetdm2019speckle}. For simulation purpose, we opt an explicit numerical method to solve the present model and then apply it on different types of gray level test images. A comparison study regarding
the quality of the despeckled image is carried out with recently developed models \cite{majeetdm2019speckle,shan2019multiplicative}. Moreover, we compare the quantitative and qualitative results at
different noise levels.
\vspace{.1cm}

The rest of the paper is organized as follows. In section \ref{sec:analysis}, we study the wellposedness of the proposed model. Section \ref{sec:numerical} describes the numerical implementation
and despeckling performance of the proposed model. We conclude the paper in Section \ref{sec:Conclusion} with a scope on future work.

\section{Existence and Uniqueness of weak solution}
\label{sec:analysis}
This section is devoted to the wellposedness result of the proposed system \eqref{maina}-\eqref{mainc}. Due to the nonlinearity in the
system \eqref{maina}-\eqref{mainc}, we first consider the associated linearized problem and then use Schauder's fixed-point theorem \cite{LCEvans1998} to complete the proof.
Without loss of generality, we assume that $\gamma=1,\iota=1,$ and $\nu=1$ in the equations \eqref{maina} and \eqref{mainb}.
%%%%%%%%%%%%%%%%%%%%%%%%%%%%%%%%%%%%%%%%%%%%%%%%%%%%%%%%%%%%%%%%%%%%%%%%%%%%%%%%%%%%%%%%%%%%%%%%%%%%%%%%%%%%%%%%%%%%%%%%%%%%%%%%%%%%%%%%%%%%%%%%%%%%%%%
\subsection{Technical framework and statement of the main result}
Throughout this article, we consider $C>0$ as a generic constant. By $(L^p, \|\cdot\|_{L^p})$ with $1\le p\le \infty$, we denote the standard spaces of $p$-th order integrable
functions on $\Omega$. Moreover, for $r\in \mathbb{N}$ 
we write $(H^r, \|\cdot\|_{H^r})$ as the usual Sobolev spaces on $\Omega$, and $(H^{1})^\prime$ as the dual space of $H^1$. 
We consider the solution space $W(0,T)$ for the underlying problem \eqref{maina}-\eqref{mainc} as $W(0,T)= W_1(0,T)\times W_2(0,T)$, where 
\begin{align*}
 W_1(0,T)&=\Big\{w\in L^\infty(0,T; H^1)\,, w_t \in L^\infty(0,T; L^2); \,w_{tt} \in L^2(0,T; (H^1)') \Big\}\,, \\
 W_2(0,T)&= \Big\{w: w\in L^\infty(0,T; H^1);\,  w_t \in L^\infty(0,T; L^2)\Big\}\,.
\end{align*}

\begin{defi}[Weak solution]\label{defi:weak}
A pair $(I,u)$ is said to be a weak solution of \eqref{maina}-\eqref{mainc}, if
\begin{itemize}
 \item[a)] $I \in W_1(0,T),~u\in W_2(0,T)$ and \eqref{mainc} holds.
 \item[b)] For all $\phi \in H^1$ and a.e $t\in (0,T)$, there hold
 \begin{align*}
  & \big\langle  I_{tt}, \phi \big\rangle + \int_{\Omega} I_t \phi \,dx
  + \int_{\Omega} \frac{s^\alpha}{ 1 + s^\alpha} \frac{1}{ 1 + |u_{\xi}|^\beta} \nabla I\cdot \nabla \phi\,dx = 0\,, \\
   &\int_{\Omega} u_t \phi \,dx + \frac{1}{2}\int_{\Omega} \nabla u\cdot \nabla \phi\,dx 
   + \int_{\Omega}u\phi\,dx = \int_{\Omega} h(|\nabla I_\xi|) \phi\,dx\,.
 \end{align*}
\end{itemize}
\end{defi}
\begin{thm}\label{thm:existence-uniqueness}
The system \eqref{maina}-\eqref{mainc} admits a unique weak solution $(I,u) \in W$ in the sense of Definition \ref{defi:weak}, provided the following two conditions hold:
\begin{Assumptions}
 \item \label{A1} $I_0 \in H^2$ satisfying $0< \rho:=\underset{x\in \Omega}\inf I_0(x)$. 
\item \label{A2}  $h:\mathbb{R}^+ \rightarrow \mathbb{R}^+$ is a bounded, Lipschitz continuous function with Lipschitz constant $c_h$ such that
\begin{align*}
 0 \leq h(\tilde{u})\leq 1\, \quad \forall\, \tilde{u}\in \mathbb{R}^+\,.
 \end{align*}
\end{Assumptions}
\end{thm}
%===================================================================================================================================
\subsection{Linearized problem $ \& $ its Well-posedness} For any positive constants $M_1, M_2>0$, define the convex set
\begin{align*}
 \mathcal{B}_{M_1,M_2}=
 \begin{cases}
 \bar{I}\in W_1(0,T):~~~\|\bar{I}\|_{L^\infty(0,T;H^1)} + \|\bar{I}_t\|_{L^\infty(0,T; L^2)} \le M_1\|I_0\|_{H^1}\,, \\
\hspace{0.1cm} 0<\rho\le \bar{I}(t,x)~~{\rm  for~a.e.}~~(t,x)\in \Omega_T\,, ~~ \\
  \bar{u}\in W_2(0,T):~~~\|\bar{u}\|_{L^\infty(0,T; L^2)} + \|\bar{u}_t\|_{L^\infty(0,T; L^2)}\le M_2\|I_0\|_{L^2}\,.
\end{cases}
\end{align*}
For any fixed $(\bar{I}, \bar{u})\in \mathcal{B}_{M_1, M_2}$, consider the linearized problem:
\begin{align}
&I_{tt} + I_{t}- {\rm div}\big( \bar{g}(t,x) \nabla I\big)=0 \hspace{4.5cm} {\rm in}~~~ \Omega_T\,, \label{linmaina} \\
& u_t = h(|\nabla \bar{I}_\xi|) - u + \frac{1}{2}\Delta u  \hspace{5.3cm} \text{in}~~~\Omega_T\,, \label{linmainb} 
\end{align}
with the condition \eqref{mainc}, where the function $\bar{g}$ is given by
\begin{align*}
 \bar{g}(t,x) \equiv g_{\bar{I}, \bar{u}}(t,x):= \frac{|\bar{I}_\xi|^\alpha}{ \big(M_\xi^{\bar{I}}\big)^\alpha + |\bar{I}_\xi|^\alpha}\cdot \frac{1}{ 1 + |\bar{u}_{\xi}|^\beta}\,.
\end{align*}
Since $(\bar{I}, \bar{u})\in \mathcal{B}_{M_1, M_2}$, a similar argument as in the proof of \cite[Claim $2.1$]{majeetdm2019speckle} revels that 
\begin{equation}\label{bound:g_w}
 \begin{aligned}
  &{\rm i)}~ 0< \kappa \le \bar{g}\le 1\,, \\
  & {\rm ii)}~ |\bar{g}_t| \le C\,,
 \end{aligned}
\end{equation}
where $\kappa, C >0$ are constants depending only on $G_\xi, I_0, M_1, M_2,\beta$, $\alpha$ and $\rho$. Hence, thanks to the classical Galerkin method \cite{LCEvans1998}, one can show that there exists a 
unique weak solution $(I,u) \in W(0,T)$ of the linearized problem \eqref{linmaina}-\eqref{linmainb} with the condition \eqref{mainc}. 

\begin{lem}\label{lem:a-priori}
The unique solution $(I,u) \in W(0,T)$ of  the linearized problem \eqref{linmaina}-\eqref{linmainb} with the condition \eqref{mainc} satisfies the following:
\begin{itemize}
\item[a)] $ \|I\|_{L^\infty(0,T; H^1)} + \|I_t\|_{L^\infty(0,T; L^2)} \le C \|I_0\|_{H^1}$, \\
\item[b)] $\int_0^T \|I_{tt}\|_{(H^1)^\prime}^2\,dt \le C T \|I_0\|_{H^1}^2$, \\
\item[c)] $ \|u\|_{L^\infty(0,T; H^1)} + \|u_t\|_{L^\infty(0,T; L^2)} \le C \|I_0\|_{H^1}$,
\end{itemize}
where $C>0$ is a constant, depends only on $G_\xi, I_0, h, M_1, M_2, \alpha,\beta$ and $\rho$. 
\end{lem}
\textbf{Proof:}
%\begin{proof}
 Since $\|\bar{u}_t\|_{L^\infty(0,T;L^2)}\le C \|I_0\|_{H^1}$, by following computations as in Sudeb at el. \cite[Lemma 3.2]{majeetdm2019speckle}, one can  show the validation of the estimates ${\rm a)}$ and ${\rm b)}$ of 
 Lemma \ref{lem:a-priori}. To prove ${\rm c)}$, we proceed as follows: multiply \eqref{linmainb} by $u_t $, integrate by parts over $\Omega$, use Cauchy-Schwarz and Young's inequalities, and then
integrate w.r.t time between $0$ to $t$. We have, for a.e. $t\in (0,T)$
\begin{align*}
\|u\|^2_{H^1}+\int_{0}^{t}\|u_t\|^2_{L^2}\, ds \le C\big( 1+t\,|\Omega |\big)\,.
\end{align*}
Moreover, since  $u_0 \in H^2$ and $ \bar{h}_t  \in L^\infty(0,\T;L^2) $, by regularity theory \cite{LCEvans1998}, $ u_t \in L^{\infty}(0,\T;L^2)$ with
\begin{align}
\| u\|^2_{H^1} + \|u_t\|^2_{L^2} \le C\,e^{t}\big( 1+t\,|\Omega |\big)\,. \label{bound-uh1-utLinf}
\end{align}
Hence ${\rm c)}$ of Lemma \ref{lem:a-priori} follows from \eqref{bound-uh1-utLinf}.
%\end{proof}
%========================================================================================================================================================
\subsection{Proof of Theorem \ref{thm:existence-uniqueness}}
As mentioned earlier, we show the well-posedness of the system \eqref{maina}-\eqref{mainc} via Schauder's fixed-point theorem. To do so, we introduce a non-empty, convex and weakly compact subset $W_0$
of $W(0, T)$ defined by
\begin{align*}
W_0=\Bigg\{ &(w,v) \in W(0,T):\,  \|w\|_{L^\infty(0,T; H^1)} + \|w_t\|_{L^\infty(0,T; L^2)} \leq C\|I_0\|_{H^1}^2\,, \\
& \|v\|_{L^\infty(0,T; H^1)} + \|v_t\|_{L^\infty(0,T; L^2)} \le C \|I_0\|_{H^1};  \\
& \hspace{2cm} ~~ 0<\rho \le w(t,x)~{\rm for ~a.e.}~(t,x)\in \Omega_T\,, 
~~\text{and}~~(w,v)~{\rm satisfies}~\eqref{mainc}\Bigg\}\,.
\end{align*}
 Consider a mapping
\begin{align*}
\mathcal{P}:~ & W_0 \goto W_0 \\
& (w,v)\mapsto (I_w, u_v)\,.
\end{align*}
If we show that the mapping $\mathcal{P}:(w,v) \rightarrow (I_w, u_v)$ is weakly continuous from $W_0$ into $W_0$, then by Schauder's fixed-point theorem, there exists $(w,v)\in W_0$ such that 
$(w,v)=\mathcal{P}(w,v)$ . In other words, the coupled system \eqref{maina}-\eqref{mainb} has a weak solution. In order to prove weak continuity of $\mathcal{P}$, let
$(w_k, v_k)$ be a sequence that converges weakly to some $(w,v)$ in $W_0$ and let $(I_k, u_k) = \big(I_{w_k}, u_{v_k}\big)$. We have to show that $\mathcal{P}(w_k, v_k):= (I_k, u_k)$ converges weakly
to $\mathcal{P}(w,v): = (I_w, u_v)$.
\vspace{.1cm}

Thanks to Lemma \ref{lem:a-priori}, one can use classical results of compact inclusion in Sobolev spaces \cite{raadams1975}
to extract subsequences $\{w_{k_n}\}$ of $\{w_k\}$, $\{v_{k_n}\}$ of $\{v_k\}$, $\{I_{k_n}\}$ of $\{I_k\}$  and $\{u_{k_n}\}$ of $\{u_k\}$, still denoted by same sequences $\{w_k\}$, $\{v_k\}$, $\{I_k\}$
and $\{u_k\}$, such that
for some $(I,u)\in W_0$, the followings hold as $k\rightarrow \infty$:
\begin{align*}
 & w_{k}\rightarrow  w\,,~~ v_{k}\rightarrow  v\ \hspace{0.2cm} {\rm in}~~L^2(0,T;L^2) \hspace{0.2cm} {\rm and~a.e.~on}~~\Omega_T\,, \\
& G_{\xi}\ast w_k  \rightarrow G_{\xi}\ast w  \hspace{0.2cm} {\rm in}~~L^2(0,T;L^2) \hspace{0.2cm} {\rm and~ a.e.~ on}~~\Omega_T\,,\\
& | G_{\xi}\ast w_k |^\alpha \rightarrow | G_{\xi}\ast w |^\alpha \hspace{0.2cm} {\rm in}~~L^2(0,T;L^2) \hspace{0.2cm} {\rm and~a.e.~on}~~\Omega_T\,,\\
& \dfrac{| G_{\xi}\ast w_k |^\alpha}{ \big(M_\xi^{w_k}\big)^\alpha + |G_{\xi}\ast w_k |^\alpha} \rightarrow \dfrac{|G_{\xi}\ast w |^\alpha}{ \big(M_\xi^{w}\big)^\alpha + |G_{\xi}\ast w |^\alpha}
\hspace{0.2cm} {\rm in}~~L^2(0,T;L^2) \hspace{0.2cm} {\rm and~a.e.~on}~~\Omega_T\,,\\
 & \partial_{x_i} G_{\xi}\ast w_k  \rightarrow \partial_{x_i} G_{\xi}\ast w ~(i=1,2) \hspace{0.2cm} {\rm in}~~L^2(0,T;L^2) \hspace{0.2cm} {\rm and~a.e.~on}~~\Omega_T\,,\\
 & h(|\nabla G_\xi \ast w_k|)\rightarrow h(|\nabla G_\xi \ast w|)\hspace{0.2cm} {\rm in}~~L^2(0,T;L^2) \hspace{0.2cm} {\rm and~ a.e.~ on}~~\Omega_T\,,\\
&  | G_\xi \ast v_k| \rightarrow | G_\xi \ast v|  \hspace{0.2cm} {\rm in}~~L^2(0,T;L^2) \hspace{0.2cm} {\rm and~ a.e.~ on}~~\Omega_T\,,\\
& \dfrac{1}{1 +| G_{\xi}\ast v_k|^\beta} \rightarrow \dfrac{1}{1 + |G_{\xi}\ast v|^\beta}  \hspace{0.2cm} {\rm in}~~L^2(0,T;L^2) \hspace{0.2cm} {\rm and~ a.e.~ on}~~\Omega_T\,,\\
& \displaystyle I_{k} \rightarrow  I\,, ~~ u_k \rightarrow u \hspace{0.2cm} \text{weakly} *~ \text{in}~~L^{\infty}(0,T;H^1)\,,\\
& \displaystyle I_{k} \rightarrow  I\,, ~~ u_k \rightarrow u \hspace{0.2cm} \text{in}~~L^{2}(0,T; L^2)\,,\\
 & \partial_t I_k \rightarrow  \partial_t I\,,~ \partial_t u_k \rightarrow  \partial_t u \hspace{0.2cm} \text{weakly} * ~\text{in}~~L^{\infty}(0,T;L^2)\,,\\
 & \partial_{tt} I_k \rightarrow \partial_{tt}I \hspace{0.2cm} \text{weakly} *~ \text{in}~~L^{2}(0,T;(H^1)^\prime)\,.
\end{align*}
In view of the above convergences, one can pass to the limit in \eqref{linmaina}-\eqref{linmainb} and obtain $(I,u)=\mathcal{P}(w,v)$. Moreover, since the
solution of \eqref{linmaina}-\eqref{linmainb} is unique, the whole
sequence $(I_k,u_k)=\mathcal{P}(w_k, v_k)$ converges weakly in $W_0$ to $(I,u)=\mathcal{P}(w,v)$. Hence $\mathcal{P}$ is weakly continuous. Therefore, the problem \eqref{maina}-\eqref{mainc} admits a weak solution.
%%%%%%%%%%%%%%%%%%%%%%%%%%%%%%%%%%%%%%%%%%%%%%%%%%%%%%%%%%%%%%%%%%%%%%%%%%%%%%%%%%%%%%%%%%%%%%%%%%%%%%%%%%%%%%%%%%%%%%%%%%%%%%%%%%%%%%%%%%%%%%%%%%%%%%%%%%%%%%%%
\vspace{.1cm}

\noindent{\bf Uniqueness of weak solution:}
To prove the uniqueness of weak solutions of the underlying problem \eqref{maina}-\eqref{mainc}, we use here a standard methodology \cite{LCEvans1998}. Let $(I_{1},u_1)$ and $(I_{2},u_2)$ be two weak 
solutions of \eqref{maina}-\eqref{mainc}.
Then, we have
\begin{align}
&I_{tt}+ I_t-\text{div} \big(g_{I_1, u_1} \nabla I\big)  = {\rm div}\big( \big(g_{I_1,u_1}-g_{I_2,u_2}\big) \nabla I_2 \big)\hspace{1cm}\text{in}~~\Omega_T\,, \label{eq:maina}\\
& u_t -\Delta u + u = h(|\nabla G_\xi \ast I_1|)- h(|\nabla G_\xi \ast I_2|) \hspace{2cm} {\rm in}~~ \Omega_T\,, \label{eq:mainb} \\
& \begin{cases} \label{eq:mainc}
 I(0,x)= 0\,,~ I_t(0,x)=0\,,~ u(0,x)= 0 \hspace{2.9cm} {\rm in}~~~\Omega\,, \\
\partial_n I =0=\partial_n u \hspace{6.3cm} {\rm on}~~~\partial \Omega_T\,,
\end{cases}
\end{align}
where $I=I_1-I_2$ and $u=u_1-u_2$. 
It suffices to show that $ (I,u) \equiv(0, 0)$ . To verify this, fix $ 0 < s < T$, and set for $i=1,2$, 
\begin{align}\label{relationvi}
v_{i}(t,\cdot)= \begin{cases}
\displaystyle \int_{t}^{s} I_{i}( \tau,\cdot)d\tau, \hspace{0.5cm} 0<t\leq s\,, \\ 
 0 \hspace{2.5cm} s \leq t < T\,.
\end{cases}
\end{align}
Note that, for $t\in (0,T)$,
\begin{align}\label{eq:fact-1}
 \begin{cases}
 \partial_t v_i(t,x)=-I_i(t,x) \quad i=1,2\,, \\
  v_{i}(t,\cdot) \in H^1\,,~~~\partial_n v_{i}=0~~~\text{on}\,\, \partial \Omega\,\,\text{in the sence of distribution}.
 \end{cases}
\end{align}
Set $v=v_1-v_2$. Then $v(s,\cdot)=0$. Now one can follow the same argumentation as in \cite[Section $3.3$]{majeetdm2019speckle} to arrive at 
\begin{align}\label{unique9}
&\frac{1}{2}\|I(s)\|_{L^2}^2+\int_{0}^{s}\|I(t)\|_{L^2}^2\,dt  + \frac{1}{2}\int_{\Omega} g_{I_1,u_1}(0,x) |\nabla v(0,x)|^2\, dx \nonumber  \\
&\leq \frac{1}{2} \Big|\int_{0}^{s}\int_{\Omega} |\nabla v|^2 \partial_t g_{I_1,u_1}\, dx\,dt\Big| + \int_0^s \|(g_{I_1, u_1}-g_{I_2, u_2})(t)\|_{L^\infty} \|\nabla I_2(t)\|_{L^2}\|\nabla v(t)\|_{L^2}\,dt\,.
\end{align}
Like in \eqref{bound:g_w}, there exist positive constants $\kappa_1, C_1>0$ such that
\begin{align*}
\kappa_1 \leq g_{I_1,u_1}\leq 1\,,\quad |\partial_t g_{I_1,u_1}|\le C\,.
\end{align*}
Moreover, by using property of convolution and the positive lower bound $\rho$ of the solutions $I_i$, we get
\begin{align*}
\|(g_{I_1,u_1}-g_{I_2,u_2})(t)\|_{L^{\infty}} \leq C(\xi, \alpha, I_0, \rho) \big(||I(t)||_{L^{2}}^\alpha + \|u\|_{L^2}\big)\,.
\end{align*}
Thus we have, from \eqref{unique9}
\begin{align*}
\frac{1}{2}\|I(s)\|_{L^2}^2+\int_{0}^{s}\| I(t)\|_{L^2}^2\,dt  + C \|\nabla v(0)\|_{L^2}^2 
&\le   C\Big( \int_0^s \big( \|\nabla v(t)\|_{L^2}^2 + \|I(t)\|_{L^2}^{2\alpha} + \|u(t)\|_{L^2}^2\big)\,dt\Big) \notag \\
& \le  C\Big( \int_0^s \big( \|v(t)\|_{H^1}^2 + \|I(t)\|_{L^2}^{2} + \|u(t)\|_{L^2}^2\big)\,dt\Big)\,,
\end{align*}
where in the last inequality, we have used the fact that $\alpha \ge 1$. Set
\begin{align*}
w_{i}(t,\cdot)=& \int_{0}^{t} I_{i}(\tau,\cdot)d\tau\, ; \quad w(t,\cdot)=(w_1-w_2)(t,\cdot)\,, \hspace{0.5cm} 0<t\leq T.
\end{align*}
Then, by using a similar argument as in \cite[Section $3.3$]{majeetdm2019speckle}, we obtain
\begin{align*}
&\frac{1}{2}\|I(s)\|_{L^2}^2+\int_{0}^{s}\| I(t)\|_{L^2}^2\,dt  + C \|w(s)\|_{H^1}^2 \notag \\
& \le \tilde{C} s\,\|w(s)\|_{H^1}^2 + C\,\int_0^s \Big( \|w(t)\|_{H^1}^2 + \|I(t)\|_{L^2}^{2} + \|u(t)\|_{L^2}^2\Big)\,dt\,.
\end{align*}
Choose $T_1$ sufficiently small such that  $C-\tilde{C} T_1 >0$. 
Then, for $0<s\leq T_1,$ we have
\begin{align}
\| I(s)\|_{L^2}^2 + \|w(s)\|_{H^1}^2 \le C \int_0^s\Big( \|w(t)\|_{H^1}^2 + \|I(t)\|_{L^2}^{2} + \|u(t)\|_{L^2}^2\Big)\,dt\,. \label{unique17}
\end{align}
Now, by multiplying \eqref{eq:mainb} by $u$ and integrating over $\Omega$, we have 
\begin{align*}
\frac{d}{dt} \|u\|_{L^2}^2 + 2 \|\nabla u\|_{L^2}^2 \le C  \Big(\| h(|\nabla G_\xi \ast I_1|)-h(|\nabla G_\xi \ast I_2|)\|_{L^2}^2 
+ \|u\|_{L^2}^2\Big)\,.
\end{align*}
Since $h$ is Lipschitz continuous, by using Young's inequality for convolution, we see that
\begin{align*}
\| h(|\nabla G_\xi \ast I_1|)-h(|\nabla G_\xi \ast I_2|)\|_{L^2}^2 \le C(c_h, \xi) \|I\|_{L^2}^2\,.
\end{align*}
Thus, we have,  for $0<s\le T_1$,
\begin{align}
u(s) \le C \int_0^s \|I(t)\|_{L^2}^2\,dt\,. \label{unique18}
\end{align}
Adding \eqref{unique17} and \eqref{unique18}, we finally get , for $0<s\le T_1$,
\begin{align*}
\| I(s)\|_{L^2}^2  + \|u(s)\|_{L^2}^2+ \|w(s)\|_{H^1}^2 \le C \int_0^s\Big( \|w(t)\|_{H^1}^2 + \|I(t)\|_{L^2}^{2} + \|u(t)\|_{L^2}^2\Big)\,dt\,. 
\end{align*}
Hence by Gronwall's lemma, we see that $ (I,u) \equiv (0,0)$ on $[0,T_1]$.
We repeatedly use the above argument on the intervals $(T_1, 2T_1]$, $(2T_1,3T_1],\ldots$ step by step, and arrive at the conclusion that 
$I_{1} = I_{2}$ and $u_1=u_2$ on $(0,T)$. This completes the proof of Theorem \ref{thm:existence-uniqueness}.
%%%%%%%%%%%%%%%%%%%%%%%%%%%%%%%%%%%%%%%%%%%%%%%%%%%%%%%%%%%%%%%%%
\vspace{.1cm}

For any weak solution $(I,u)$ of \eqref{maina}-\eqref{mainc}, we next show the boundedness of $I$ under the assumption that initial image $I_0$ has a
finite upper bound, whose proof follows from the proof of \cite[Lemma $3.3$]{majeetdm2019speckle}. 
\begin{lem}
 Let $(I,u)$ be a weak solution of the system \eqref{maina}-\eqref{mainc}, and $\varrho:= \underset{x\in \Omega}\sup I_0(x) < \infty$.  Then 
 \begin{align}
  0<\rho \le I(t,x)\le \varrho \quad {\rm for ~a.e.}~(t,x)\in \Omega_T\,. \label{boundedness-weak-solution}
 \end{align}
\end{lem}

%===================================================================================================================================
\section{Numerical method and Experimental Results}
\label{sec:numerical}
In this section, we show the image despeckling performance of the suggested model over two existing approaches \cite{majeetdm2019speckle,shan2019multiplicative}. To solve the model  \eqref{maina}-\eqref{mainc}  numerically, we opt an explicit finite difference scheme. We replace the derivative terms in the model \eqref{maina}-\eqref{mainc} using the following finite difference formulas:
\begin{align*}
& \dfrac{\partial I_{i,j}^n }{\partial t} \approx  \displaystyle\frac{I_{i,j}^{n+1}-I_{i,j}^n}{\tau}\,, \quad 
\dfrac{\partial^2 I_{i,j}^n }{\partial t^2} \approx \displaystyle\frac{I_{i,j}^{n+1}-2I_{i,j}^n+I_{i,j}^{n-1}}{\tau ^2}\,,\\
& \nabla_x I_{i,j}^n \approx \displaystyle\frac{I_{i+1,j}^n-I_{i-1,j}^n}{2\tilde{h}}\,, \quad 
\nabla_y I_{i,j}^n \approx \displaystyle\frac{I_{i,j+1}^n-I_{i,j-1}^n}{2\tilde{h}}\,, \\
& \Delta_x I_{i,j}^n \approx  \frac{{I_{i+1,j}^n-2I_{i,j}^n+I_{i-1,j}^n}}{\tilde{h}^2}\,, \quad 
\Delta_y I_{i,j}^n \approx \frac{{I_{i,j+1}^n-2I_{i,j}^n+I_{i,j-1}^n}}{\tilde{h}^2}\,, \\
& |\nabla I_{i,j}^n| \approx \sqrt{(\nabla_x I_{i,j}^n)^2 + (\nabla_y I_{i,j}^n)^2}\,.
\end{align*}
In the above, $\tau$ resp. $\tilde{h}$ denotes the time step size resp. the spatial step size. $I^n_{i,j}=I(t_n,x_i,y_j)$ where $x_i=i\tilde{h}~(i=0,1,2...,N)$, 
$y_j=j\tilde{h}~(j=0,1,2...,M)$, $t_n=n\tau~(n=0,1,2,\ldots)$ where $n$ is the number of iterations and $M \times N$ is the image dimension. Then, the discrete form of the equation \eqref{maina}
could be written as 
\begin{align}\label{disc:I}
&(1+\gamma \tau)I_{i,j}^{n+1}=(2+\gamma \tau)I_{i,j}^n -I_{i,j}^{n-1} + {\tau ^2} \big\{ \nabla_x \left( g_{i,j}^n \nabla_x I_{i,j}^n  \right) 
 + \nabla_y \left( g_{i,j}^n \nabla_y I_{i,j}^n \right) \big\} ,
\end{align}
where
\begin{align*}
g_{i,j}^n=b(s^{n}_{i,j})\cdot  \dfrac{1}{1+ \iota|G_{\xi} \ast u^n_{i,j} |^\beta }~~~\text{with}~~b(s)=\frac{s^\alpha}{1 + s^\alpha}\,.
\end{align*}
Moreover, $u(t_n,x_i,y_j)=u^n_{i,j}$ is calculated from the discretized equation of \eqref{mainb} as follows
\begin{align}\label{disc:u}
&u_{i,j}^{n}=u_{i,j}^{n-1} + \tau \big\{h_{i,j}^{n} -u_{i,j}^{n-1} +   \frac{\nu^2}{2} \big(\Delta_x u_{i,j}^{n-1} +  \Delta_y u_{i,j}^{n-1}\big) \big\}\,, 
\end{align}
where $ h_{i,j}^n= h\big(|\nabla(G_{\xi} \ast I^n_{i,j})| \big)$. We choose the function $h$ as $h(\theta)= \epsilon +\text{min}\{ \theta^2, K \}$ for numerical experiments, where $K$ is square of the maximum gray level value of the image $I$ and $\epsilon >0$ is a very small number. The boundary and initial conditions are given as follows:
\begin{align*}
& I_{-1,j}^n=I_{0,j}^n,~~I_{N+1,j}^n=I_{N,j}^n\,, \quad  I_{i,-1}^n=I_{i,0}^n\,,~~I_{i,M+1}^n=I_{i,M}^n\,, \\
& I_{i,j}^0 =I_0(x_i,y_j)\,, \quad  I_{i,j}^1=I_{i,j}^0,\,\,\,\,	0 \leq i \leq N\,, ~~ 0 \leq j \leq M\,.
\end{align*}
We choose similar boundary conditions for edge variable $u$, and $u_{i,j}^0 =G_{\xi} \ast|\nabla  I^0_{i,j} |^2.$ We start the simulation with the initial value $I_0$ and utilize
the equations \eqref{disc:I} and \eqref{disc:u} repeatedly to find a sequence of values of ${I(t,x)};t>0$, which represents the filtered versions of $I_0$. We made a stopping criterion for
the noise elimination process when the best PSNR \cite{gonzalez2002digital} value for the restored image $I(t,x)$ is reached. 
\vspace{.1cm}

\begin{figure}[]
       \centering
       \begin{subfigure}[b]{0.21\textwidth}           
           \includegraphics[scale=0.31]{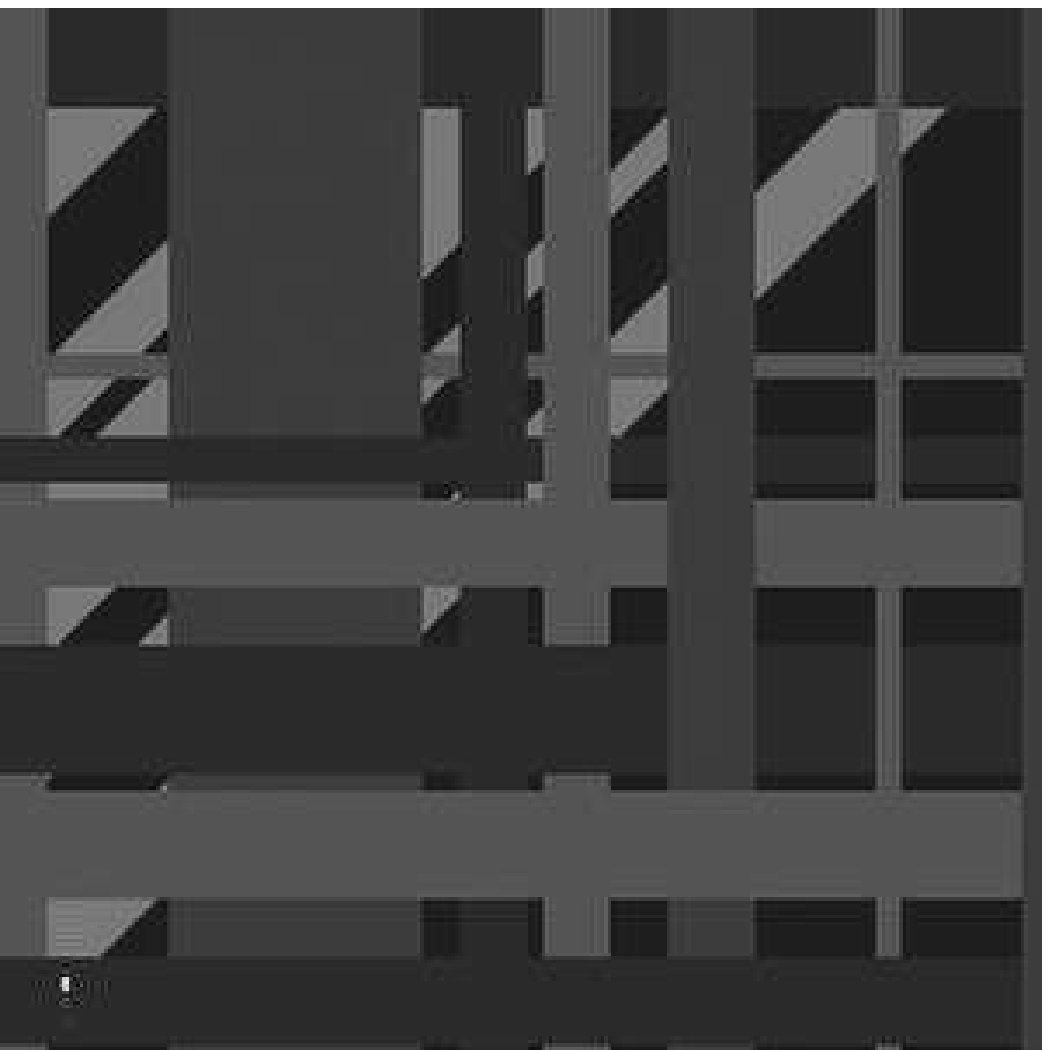}           
                \caption{Texture}
                \label{fig:boat_clear}
       \end{subfigure}%
              \begin{subfigure}[b]{0.21\textwidth}           
                \includegraphics[scale=0.36]{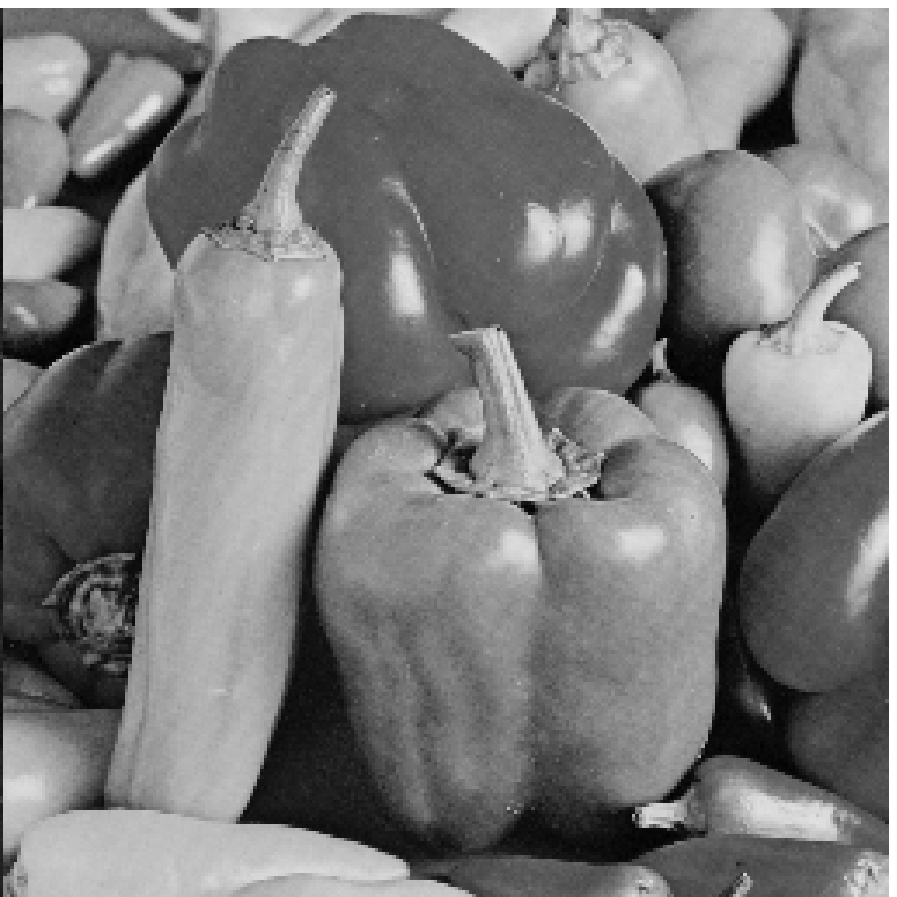}               
                \caption{Peppers}
                \label{fig:brick_clear}
       \end{subfigure}% 
      \begin{subfigure}[b]{0.21\textwidth}           
                \includegraphics[scale=0.309]{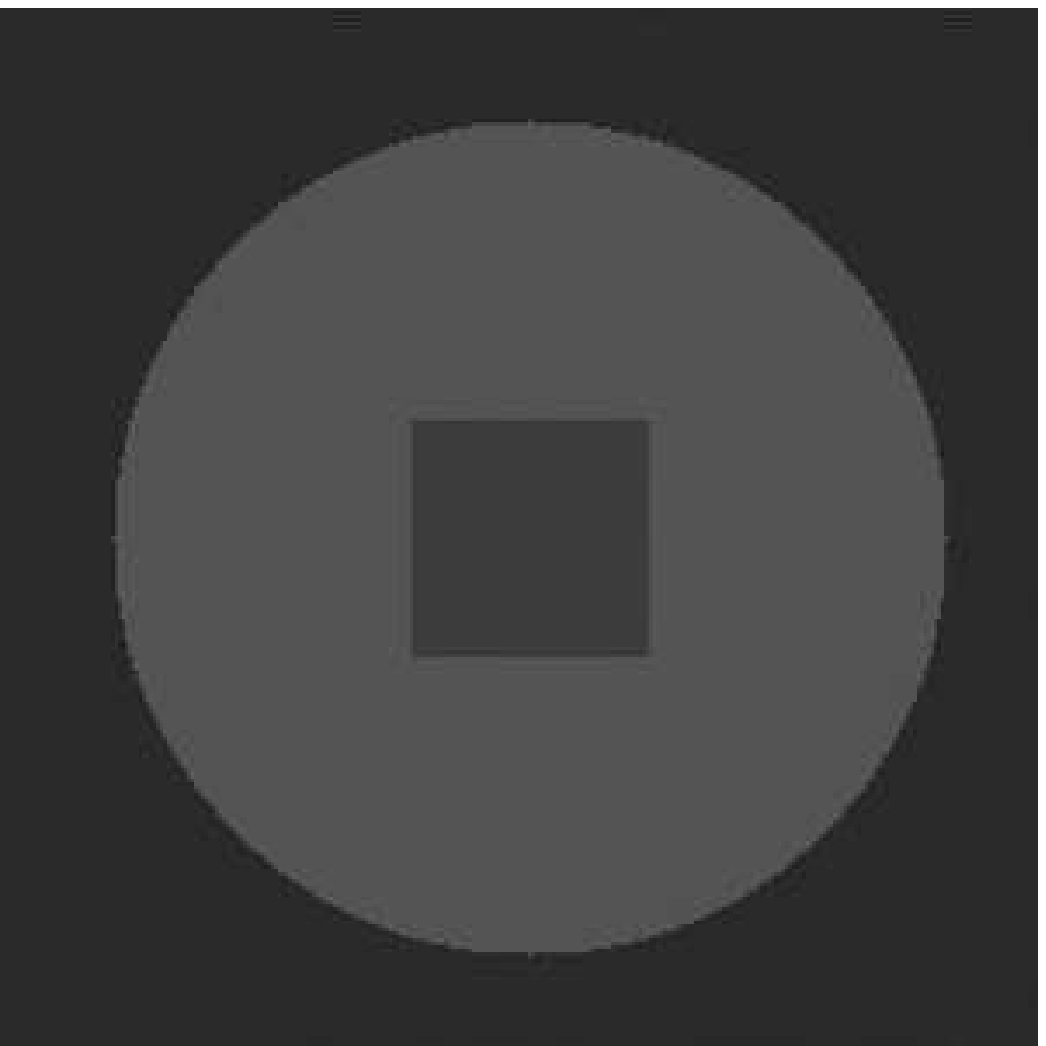}               
                \caption{Circle}
                \label{fig:circle_clear}
       \end{subfigure}%
       
 \caption{Test Images: (a) Texture Image, (b) Natural Image, (c) Synthetic Image.}\label{fig:all_clear_images}
\end{figure}

We perform all the experiments on three standard test images \ref{fig:all_clear_images}, which are initially degraded with different level of speckle noise.  All the numerical experiments are
performed under windows $7$ and MATLAB version $R2018b$ running on a desktop with an Intel Core $i5$ dual-core CPU at $2.53$ GHz with $4$ GB of memory. In this process, we use same numerical scheme 
as done it for the proposed model to discretize the considered existing models. We choose an uniform time step size $\tau = 0.2$, spatial step size $\tilde{h}=1$, and $\xi=1$ for each models.

\begin{figure}
       \centering

      \begin{subfigure}[b]{0.24\textwidth}           
                \includegraphics[scale=0.32]{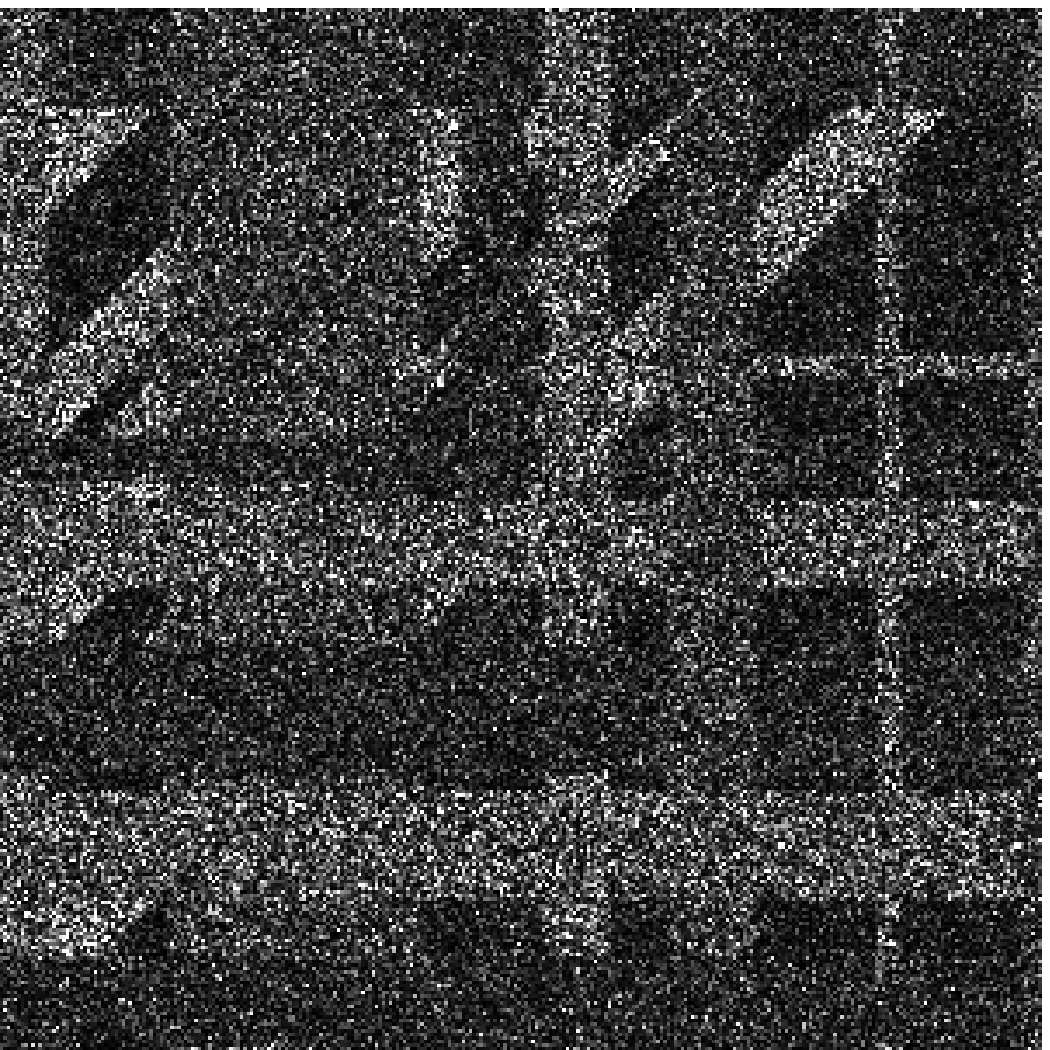}               
                \caption{Noisy}
                \label{fig1:texture_1}
       \end{subfigure}% 
       \begin{subfigure}[b]{0.24\textwidth}           
                \includegraphics[scale=0.32]{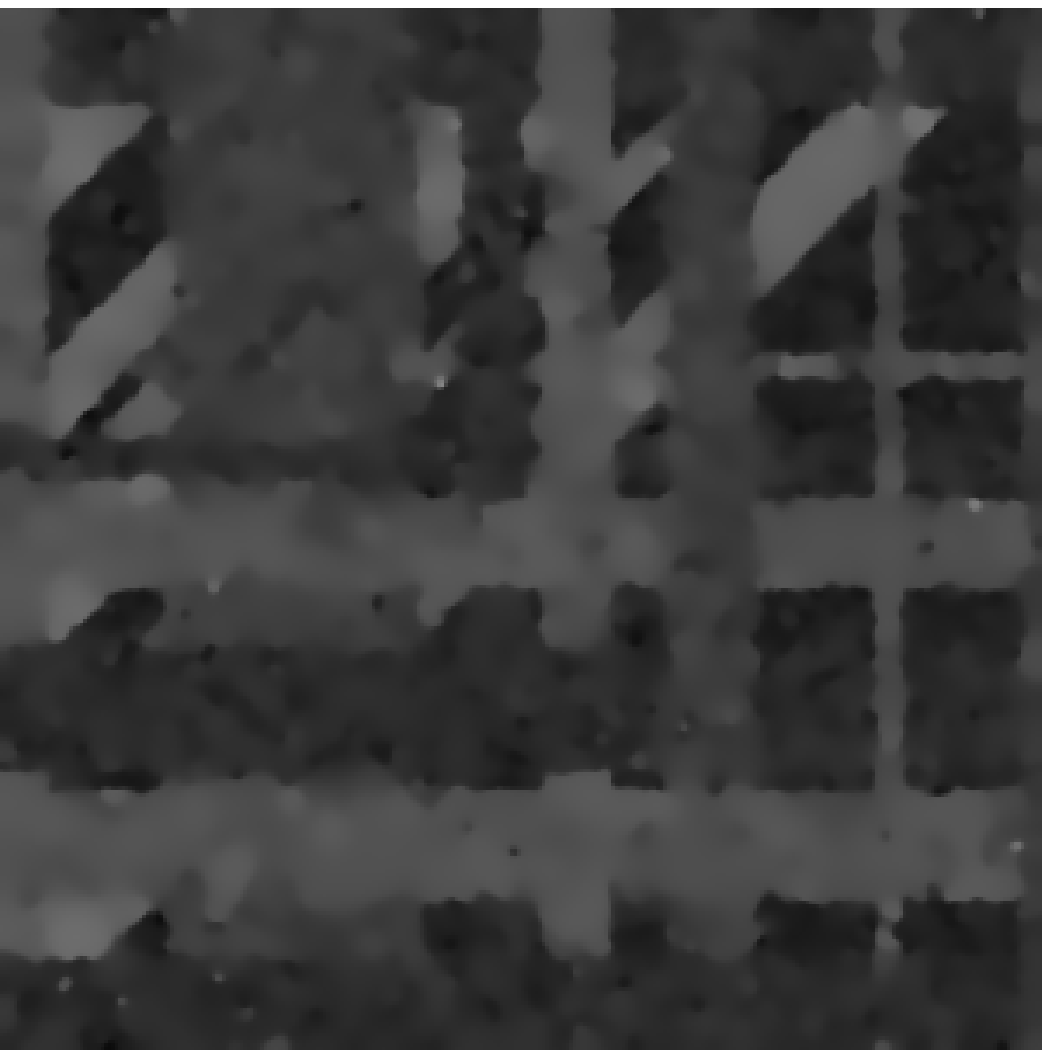}               
                \caption{Shan}
                \label{fig1:texture_1_shan}
       \end{subfigure}%
         \begin{subfigure}[b]{0.24\textwidth}           
                \includegraphics[scale=0.32]{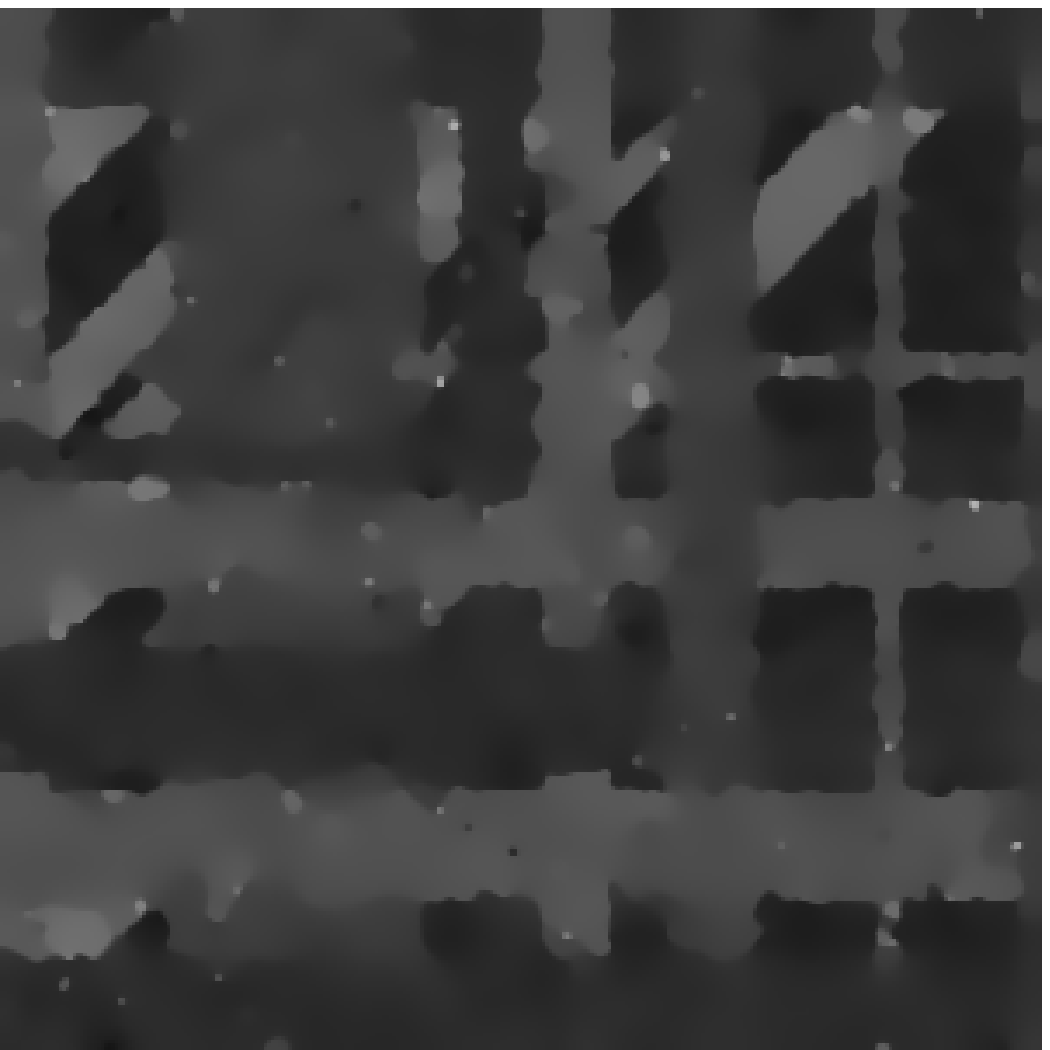}               
                \caption{TDM}
                \label{fig1:texture_1_tdm}
       \end{subfigure}% 
        \begin{subfigure}[b]{0.24\textwidth}           
                \includegraphics[scale=0.32]{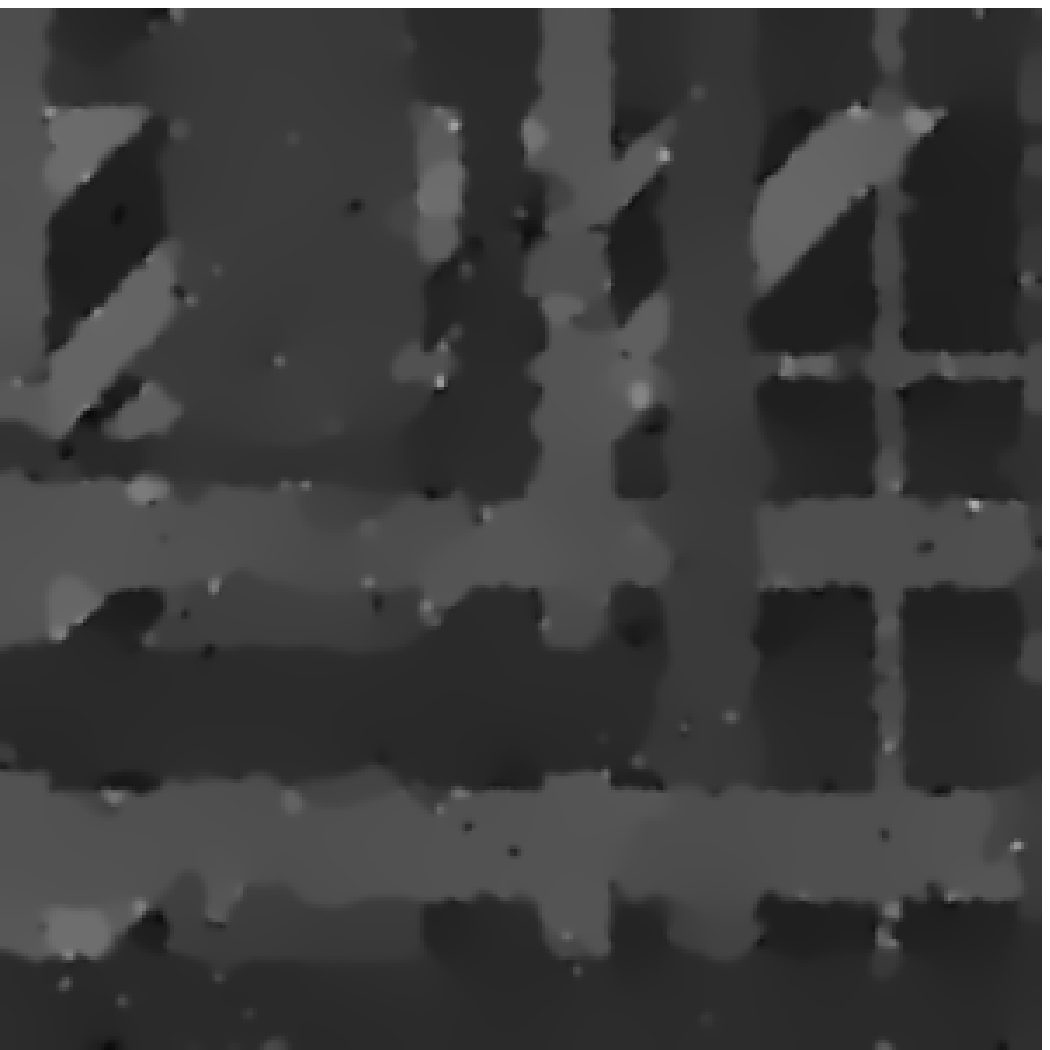}               
                \caption{Proposed}
                \label{fig1:texture_1_tdm_hp}
       \end{subfigure}
     
\caption{Image corrupted with speckle look L=1 and restored by different models.}\label{texture_1}
\end{figure}

\begin{figure}
       \centering

      \begin{subfigure}[b]{0.24\textwidth}           
                \includegraphics[scale=0.32]{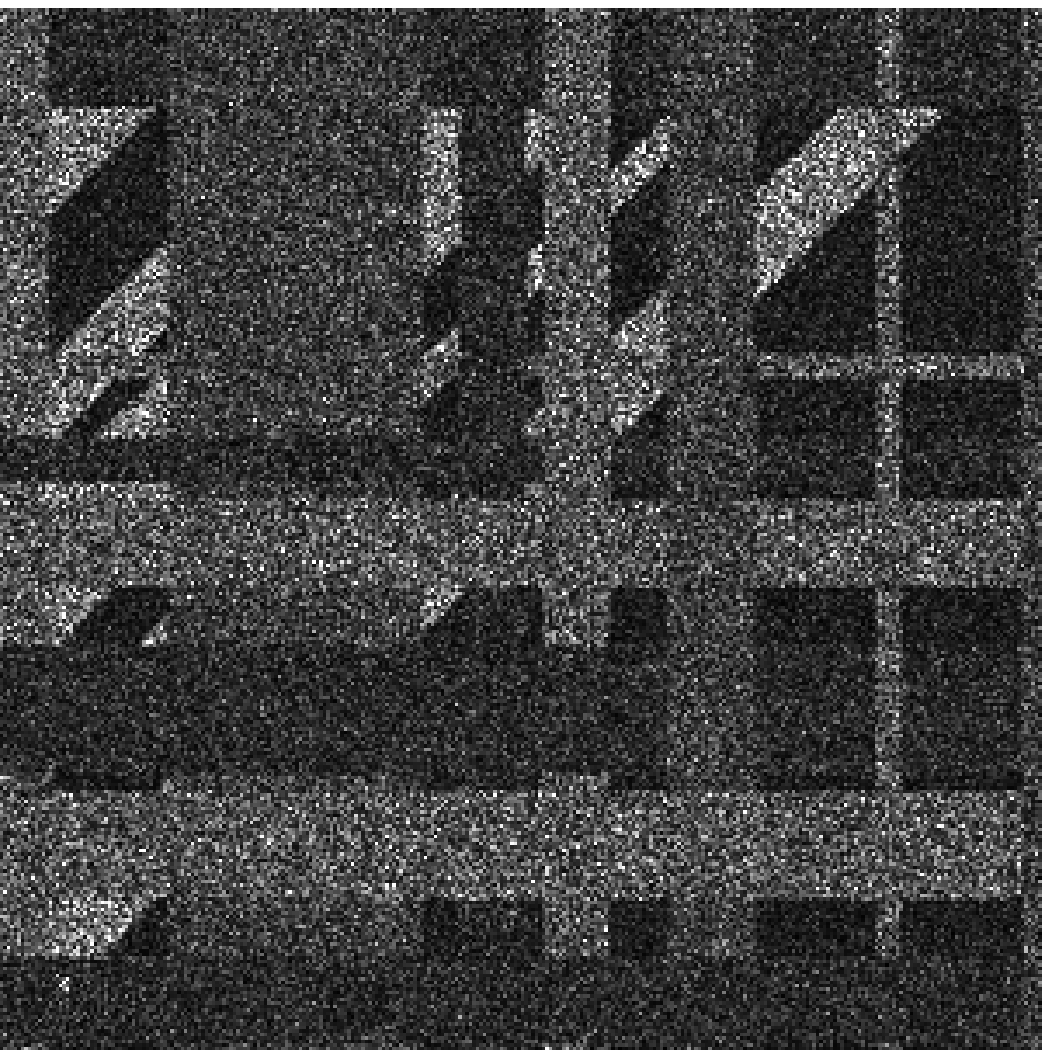}               
                \caption{Noisy}
                \label{fig2:texture_3}
       \end{subfigure}% 
       \begin{subfigure}[b]{0.24\textwidth}           
                \includegraphics[scale=0.32]{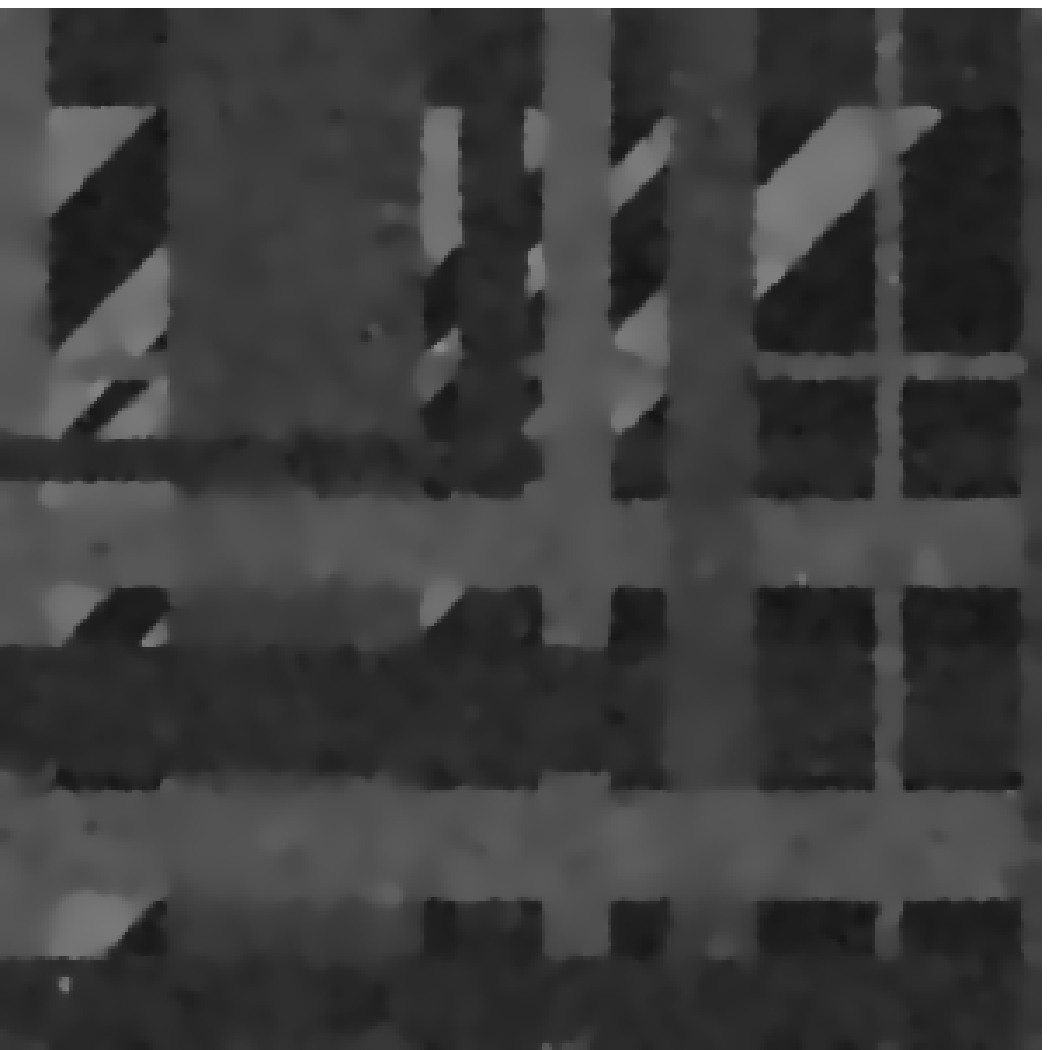}               
                \caption{Shan}
                \label{fig2:texture_3_shan}
       \end{subfigure}%
         \begin{subfigure}[b]{0.24\textwidth}           
                \includegraphics[scale=0.32]{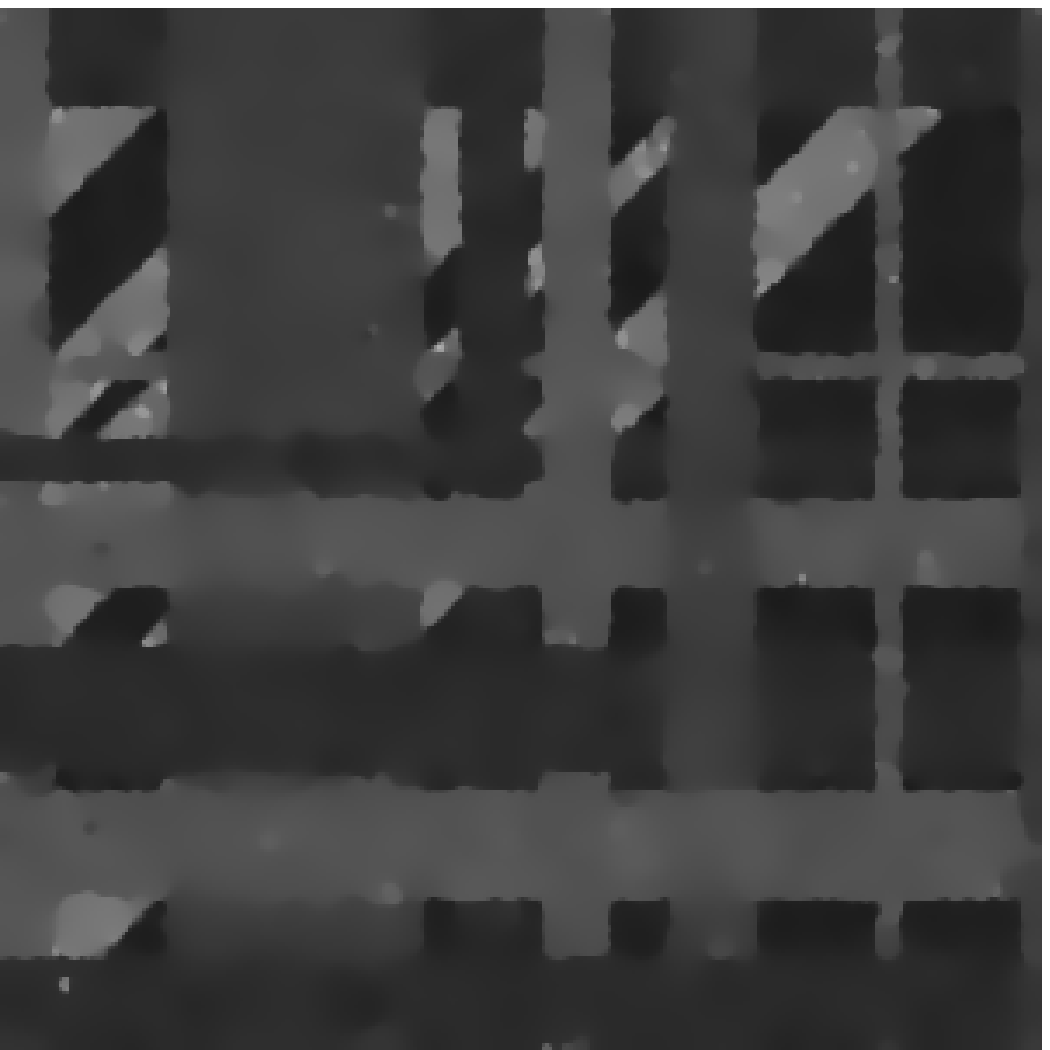}               
                \caption{TDM}
                \label{fig2:texture_3_tdm}
       \end{subfigure}% 
        \begin{subfigure}[b]{0.24\textwidth}           
                \includegraphics[scale=0.32]{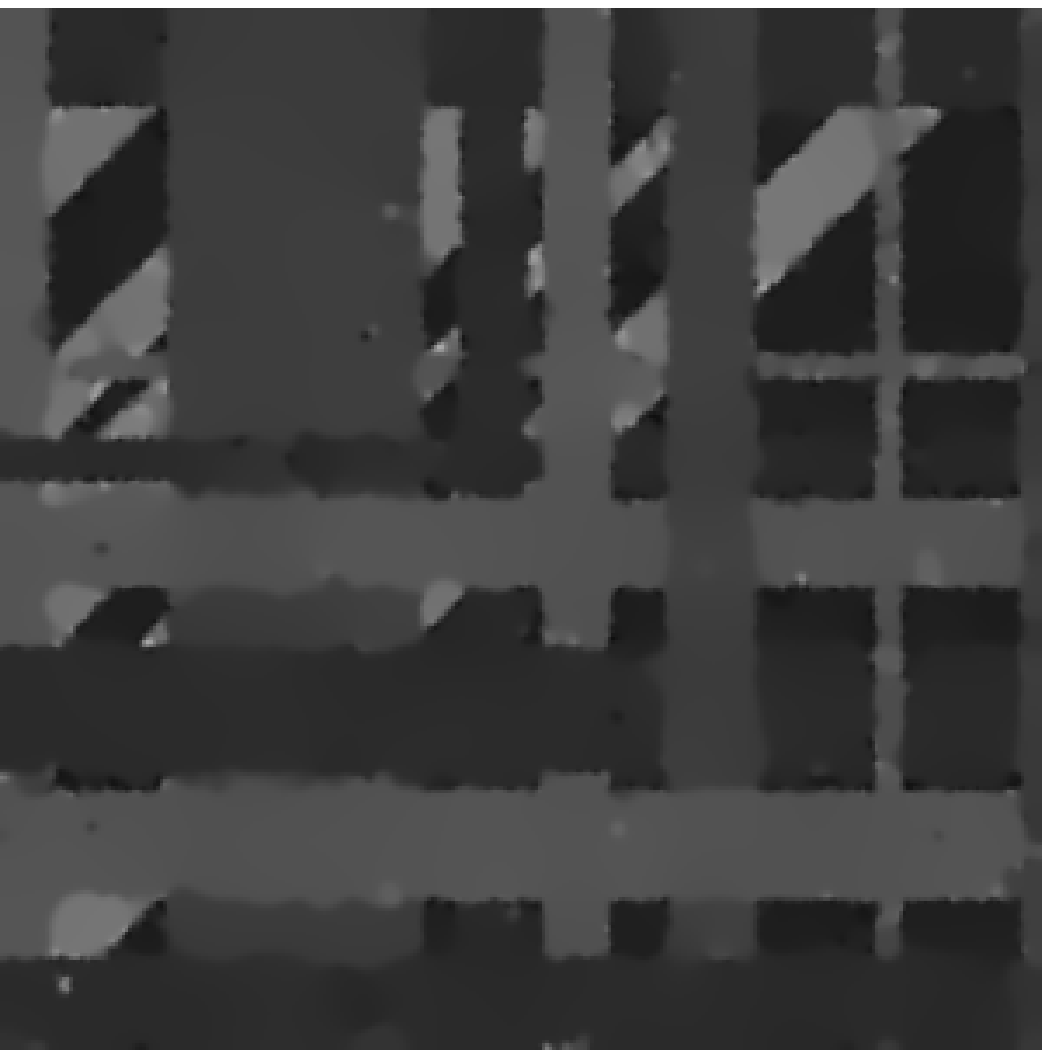}               
                \caption{Proposed}
                \label{fig2:texture_3_tdm_hp_new}
       \end{subfigure}
     
\caption{Image corrupted with speckle look L=3 and restored by different models.}\label{texture_3}
\end{figure}

\begin{figure}
       \centering

      \begin{subfigure}[b]{0.24\textwidth}           
                \includegraphics[scale=0.32]{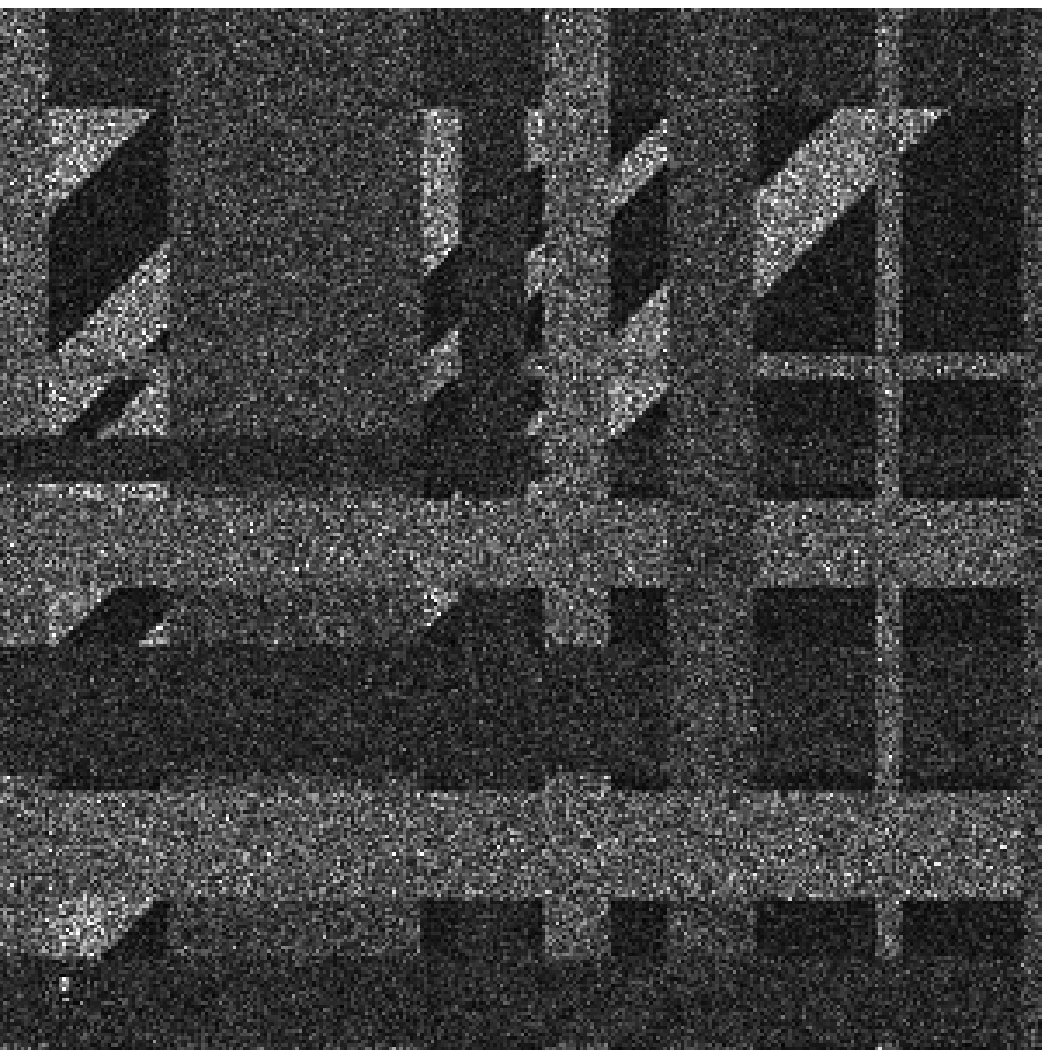}               
                \caption{Noisy}
                \label{fig3:texture_5}
       \end{subfigure}% 
       \begin{subfigure}[b]{0.24\textwidth}           
                \includegraphics[scale=0.32]{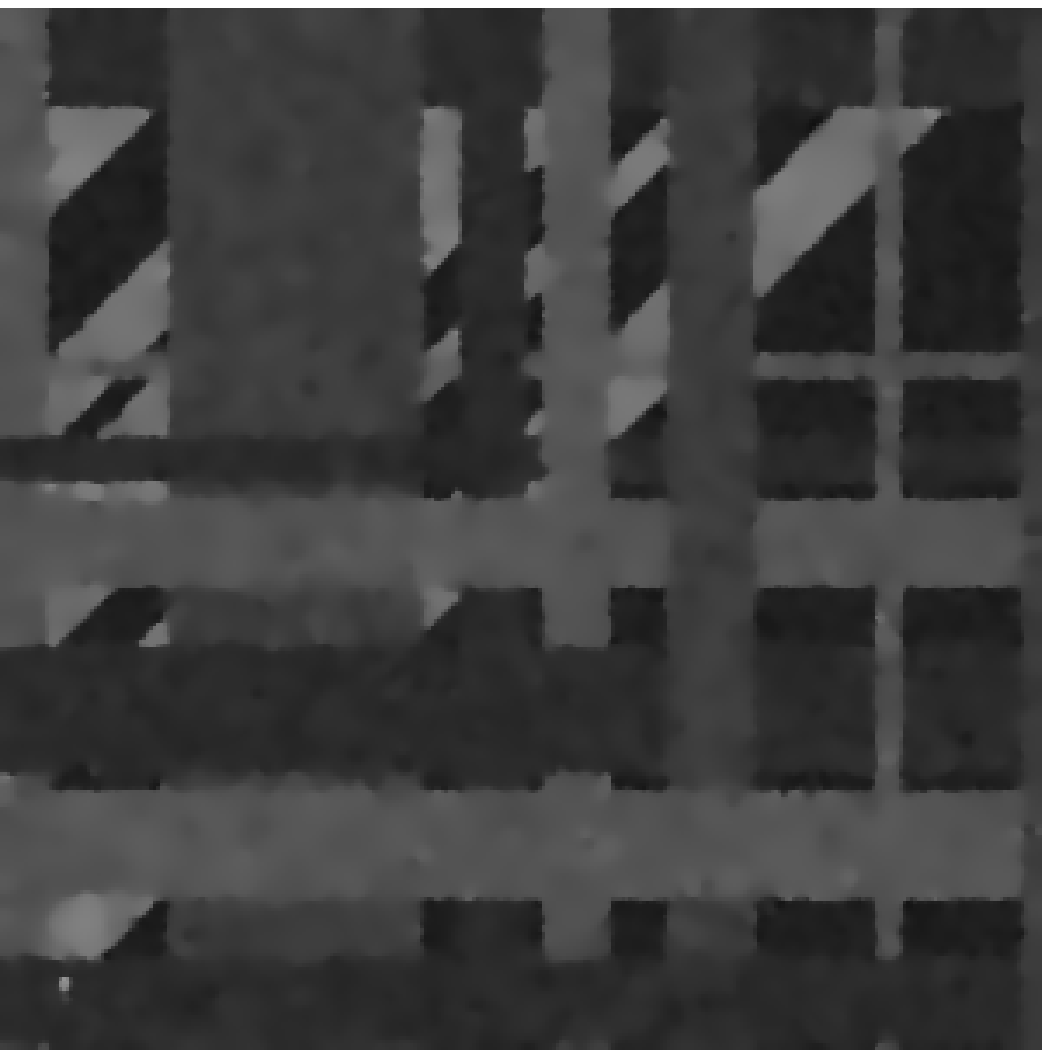}               
                \caption{Shan}
                \label{fig3:texture_5_shan}
       \end{subfigure}%
         \begin{subfigure}[b]{0.24\textwidth}           
                \includegraphics[scale=0.32]{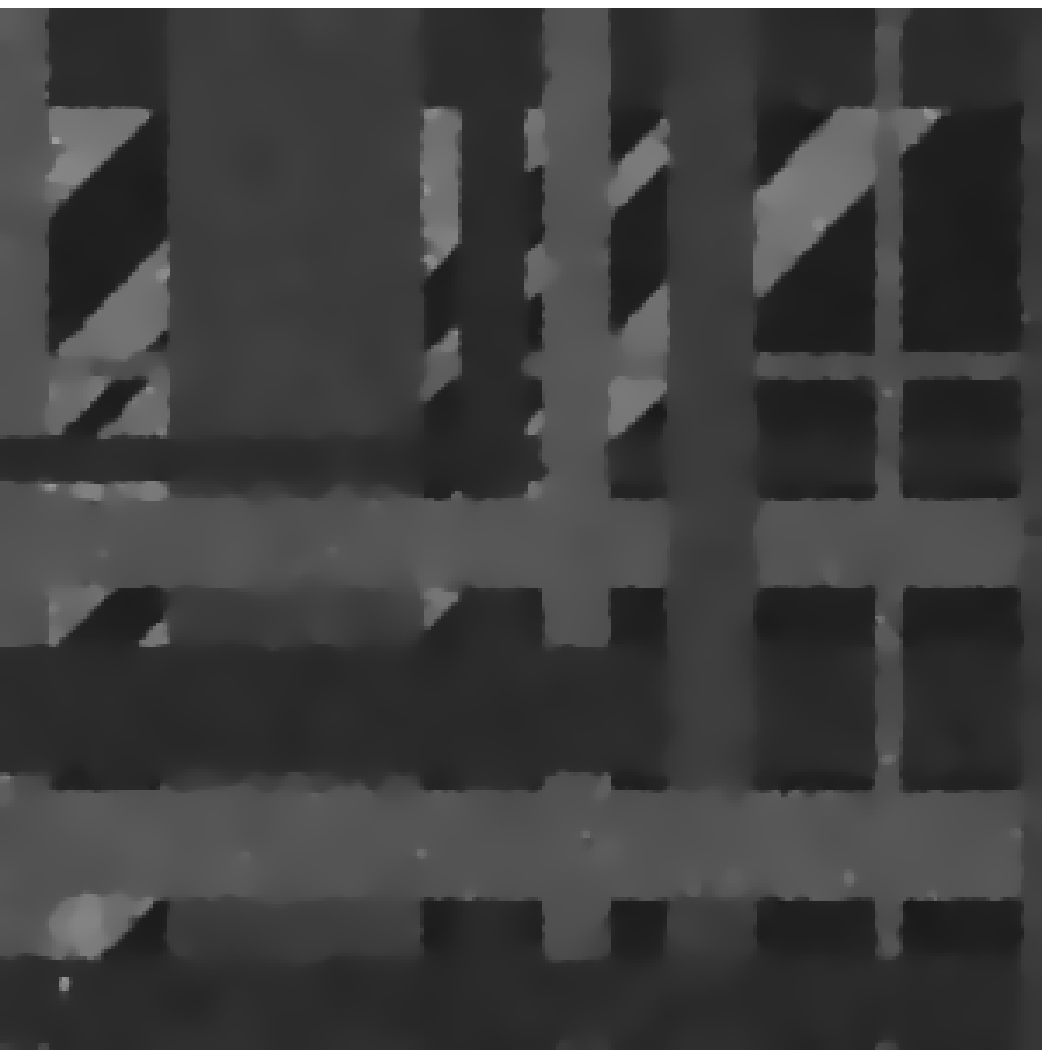}               
                \caption{TDM}
                \label{fig3:texture_5_tdm}
       \end{subfigure}% 
        \begin{subfigure}[b]{0.24\textwidth}           
                \includegraphics[scale=0.32]{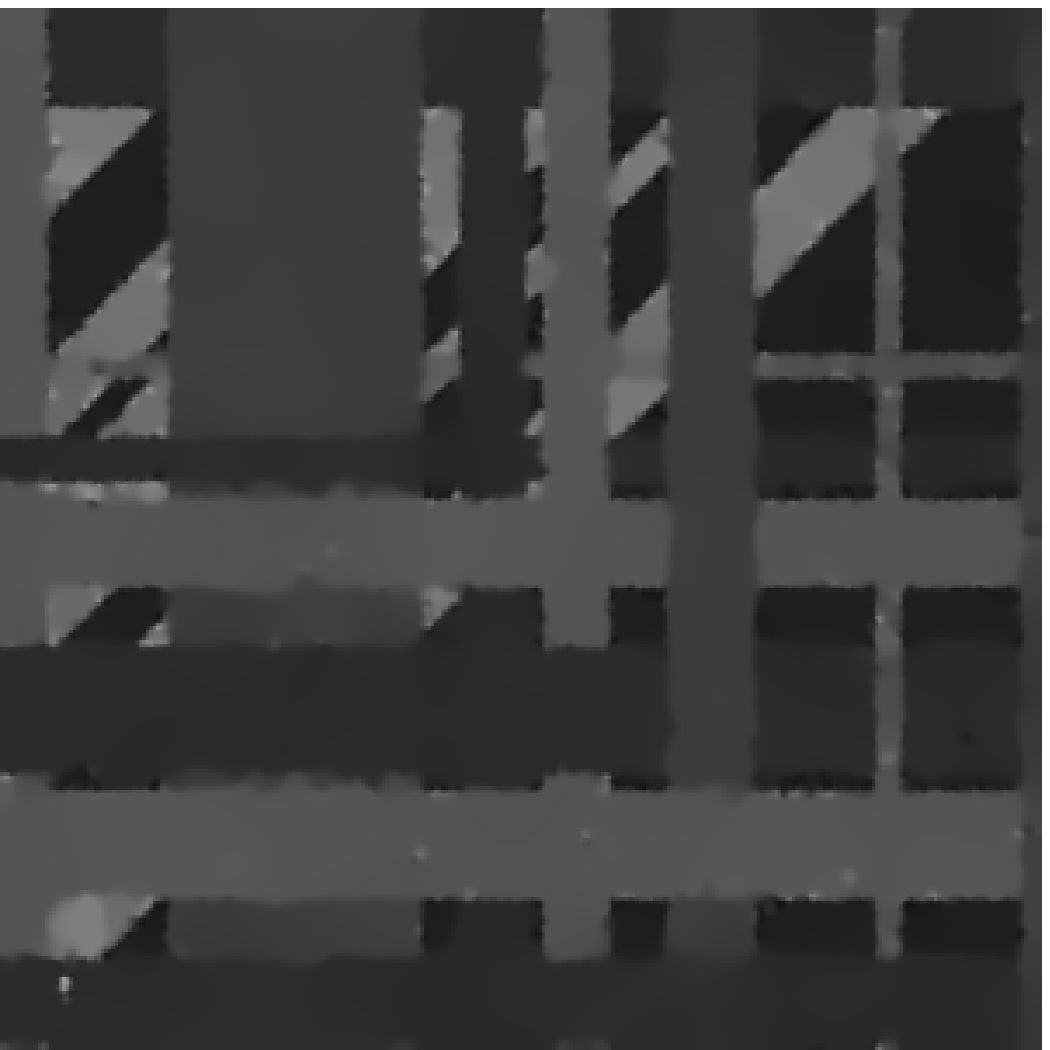}               
                \caption{Proposed}
                \label{fig3:texture_5_tdm_hp}
       \end{subfigure}
     
\caption{Image corrupted with speckle look L=5 and restored by different models.}\label{texture_5}
\end{figure}

To compare the quantitative results, we compute the values of the two standard parameters PSNR\cite{gonzalez2002digital} and MSSIM \cite{wang2004image}. A higher numerical value of MSSIM and PSNR 
suggests that the despeckled image is closer to the noise-free image. Apart from the despeckled image, for qualitative observations,  we compute the 2D contour plot, 3D surface plot for the better
visualization of the computational results.
%=========================================================================================================================
%\iffalse
%=========================================================================================================================
%============================================================================================================================
In figures \ref{texture_1}-\ref{texture_5}, we represent the restored results of a Texture image which is contaminated by multiplicative speckle noise with  $L=\lbrace1,3,5\rbrace$. From figure \ref{fig1:texture_1_shan}, we can see that the Shan model failed to preserve the fine edges for very high noise level. TDM model works better than the Shan model, but the present model preserves the fine edges better than Shan and TDM models. 

\begin{figure}
       \centering
      
      \begin{subfigure}[b]{0.24\textwidth}           
                \includegraphics[scale=0.4]{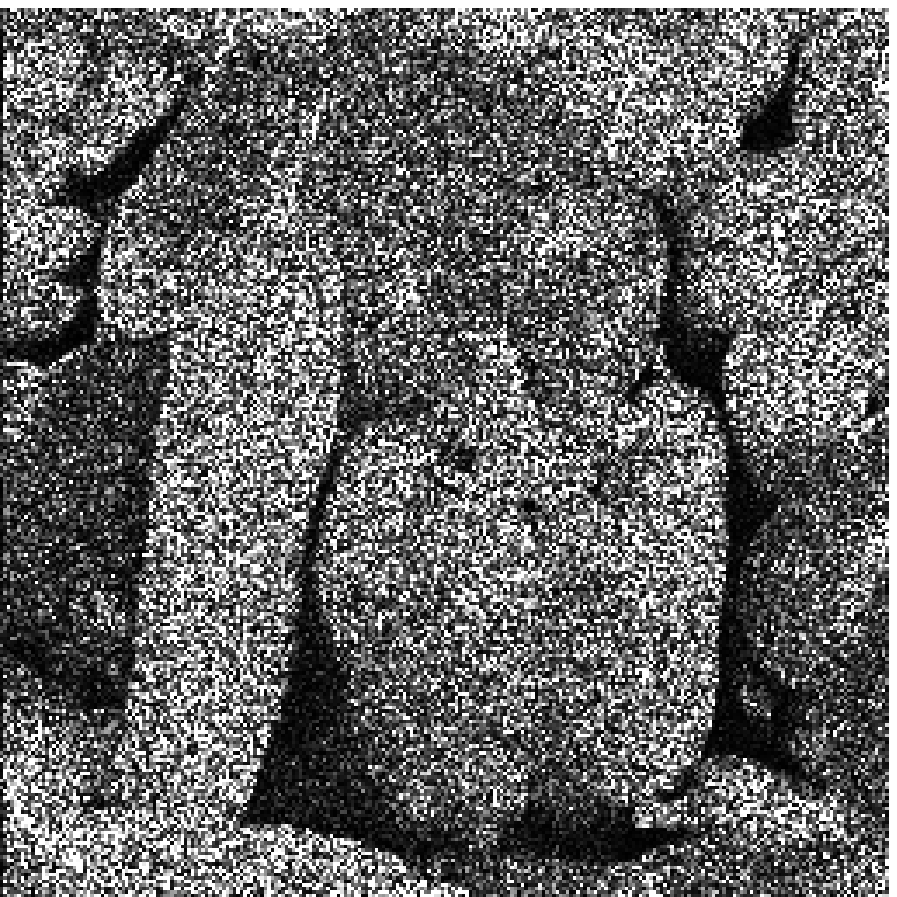}               
                \caption{Noisy}
                \label{fig1:peppers_1}
       \end{subfigure}% 
       \begin{subfigure}[b]{0.24\textwidth}           
                \includegraphics[scale=0.4]{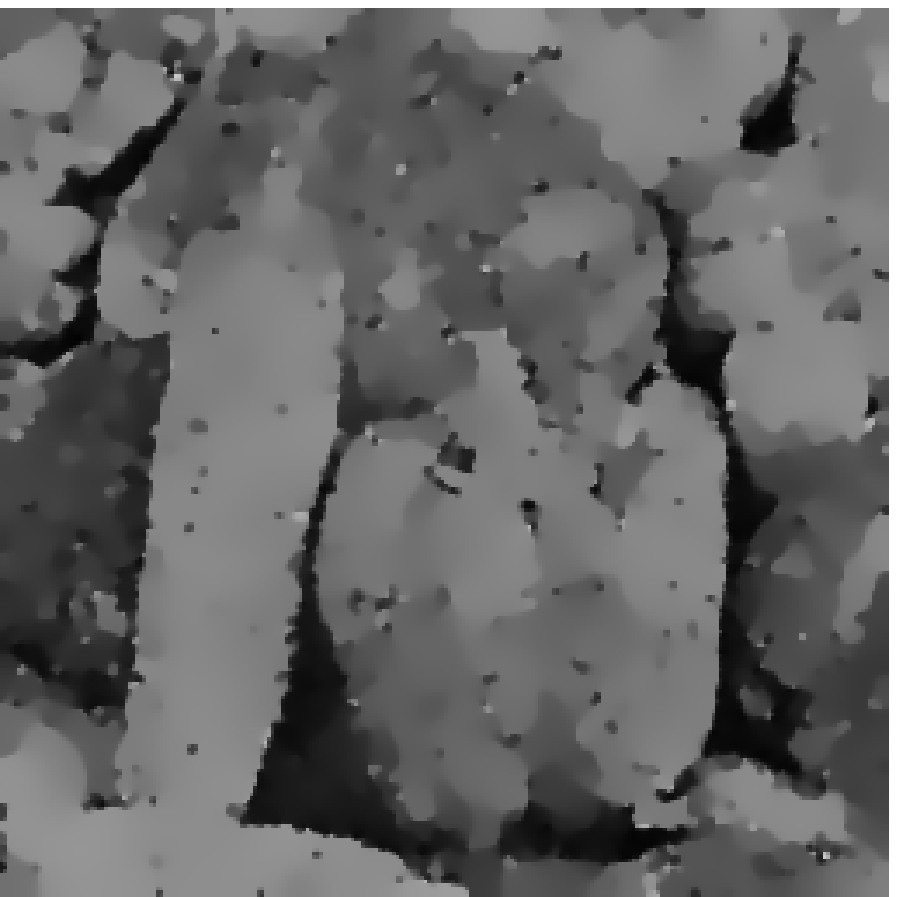}               
                \caption{Shan}
                \label{fig1:peppers_1_shan}
       \end{subfigure}%
         \begin{subfigure}[b]{0.24\textwidth}           
                \includegraphics[scale=0.4]{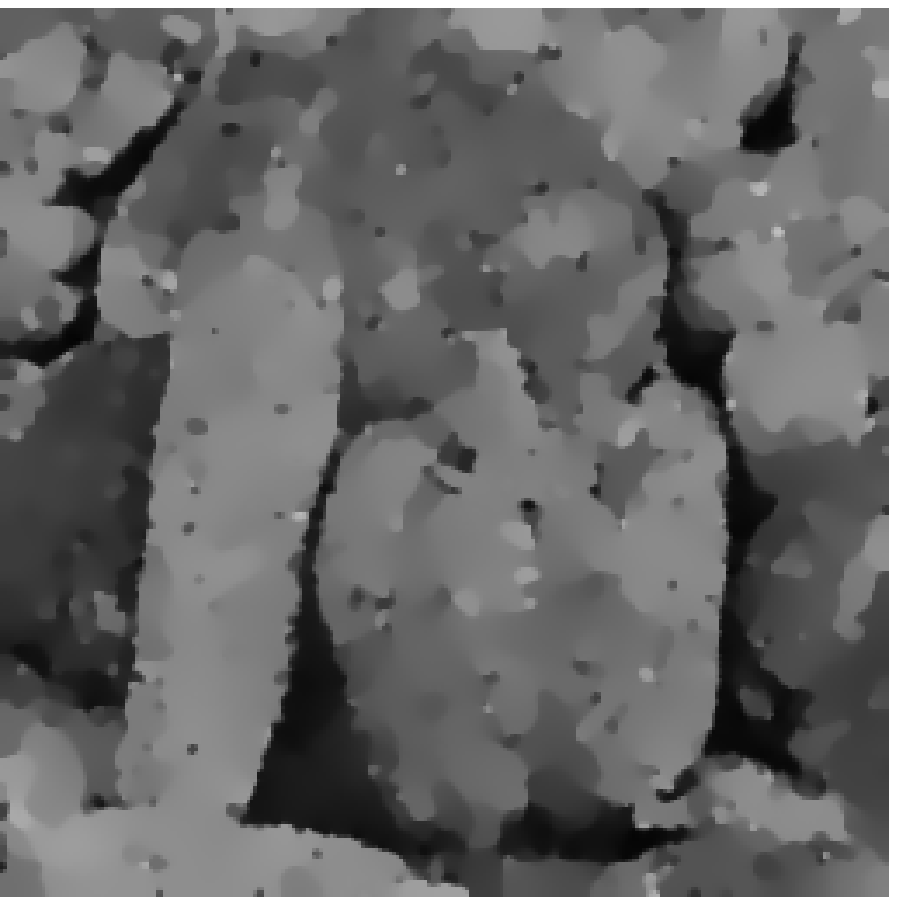}               
                \caption{TDM}
                \label{fig1:peppers_1_tdm}
       \end{subfigure}% 
        \begin{subfigure}[b]{0.24\textwidth}           
                \includegraphics[scale=0.4]{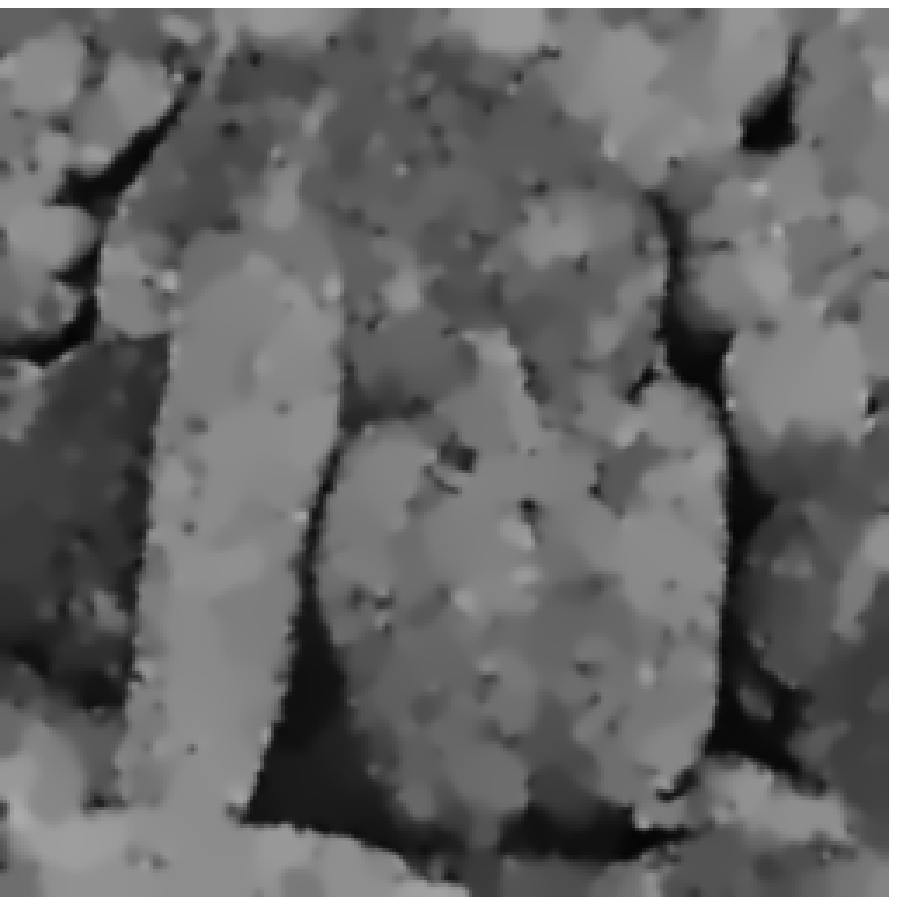}               
                \caption{Proposed}
                \label{fig1:peppers_1_tdm_hp}
       \end{subfigure}
     
\caption{Image corrupted with speckle look L=1 and restored by different models.}\label{peppers_1}
\end{figure}

\begin{figure}
       \centering

       \begin{subfigure}[b]{0.24\textwidth}           
                \includegraphics[scale=0.4]{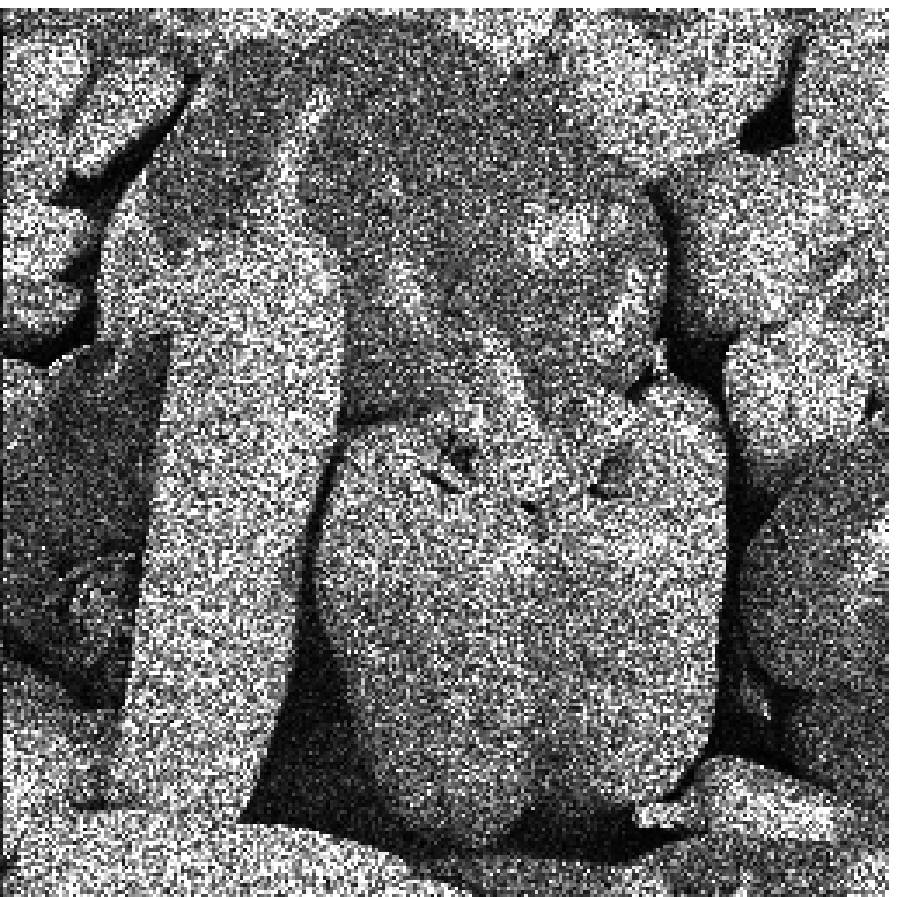}               
                \caption{Noisy}
                \label{fig2:peppers_3}
       \end{subfigure}% 
   \begin{subfigure}[b]{0.24\textwidth}           
                \includegraphics[scale=0.4]{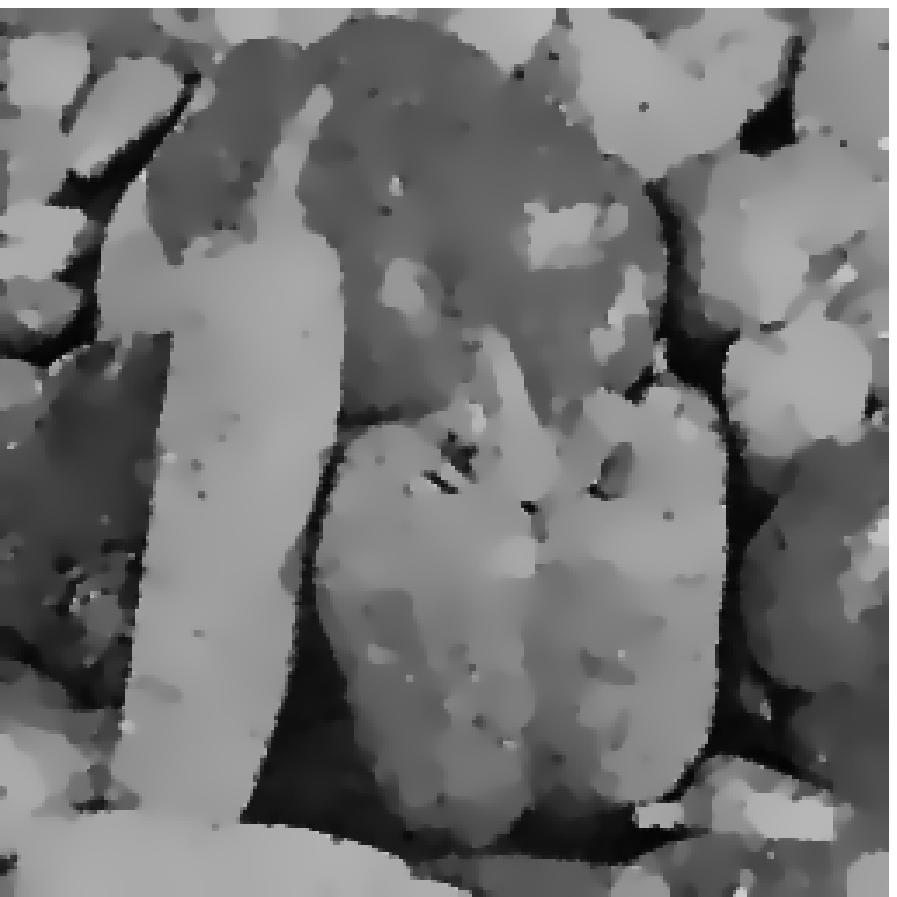}               
                \caption{Shan}
                \label{fig2:peppers_3_shan}
       \end{subfigure}%
         \begin{subfigure}[b]{0.24\textwidth}           
                \includegraphics[scale=0.4]{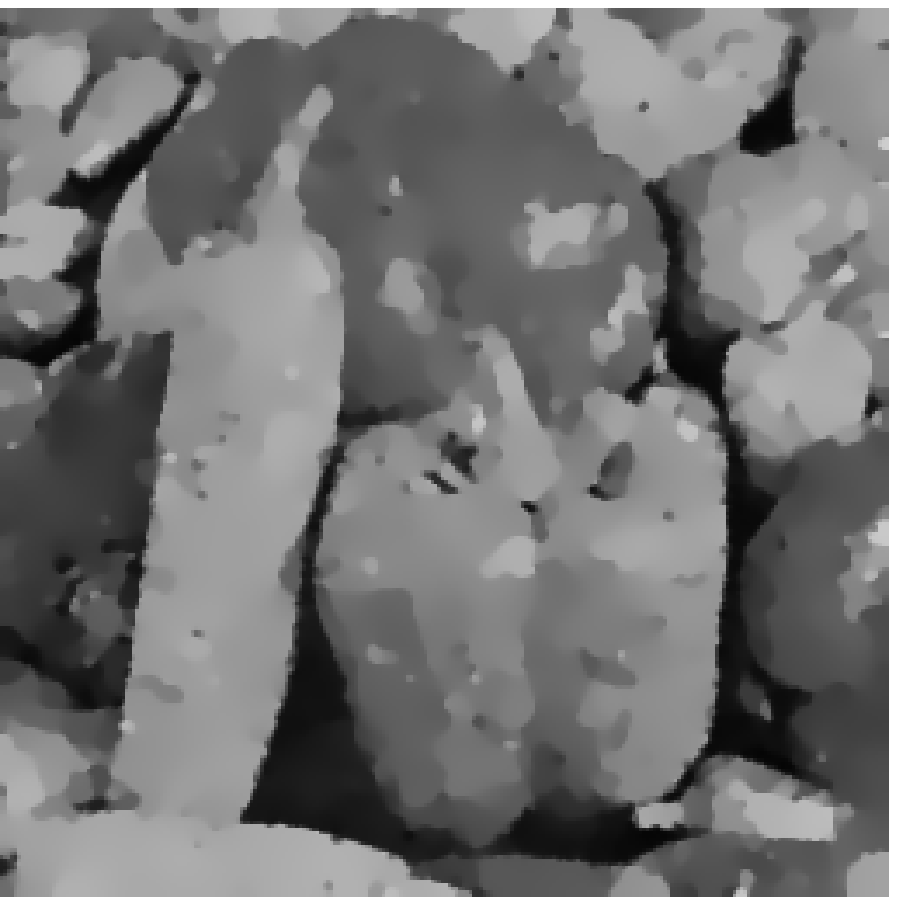}               
                \caption{TDM}
                \label{fig2:peppers_3_tdm}
       \end{subfigure}% 
        \begin{subfigure}[b]{0.24\textwidth}           
                \includegraphics[scale=0.4]{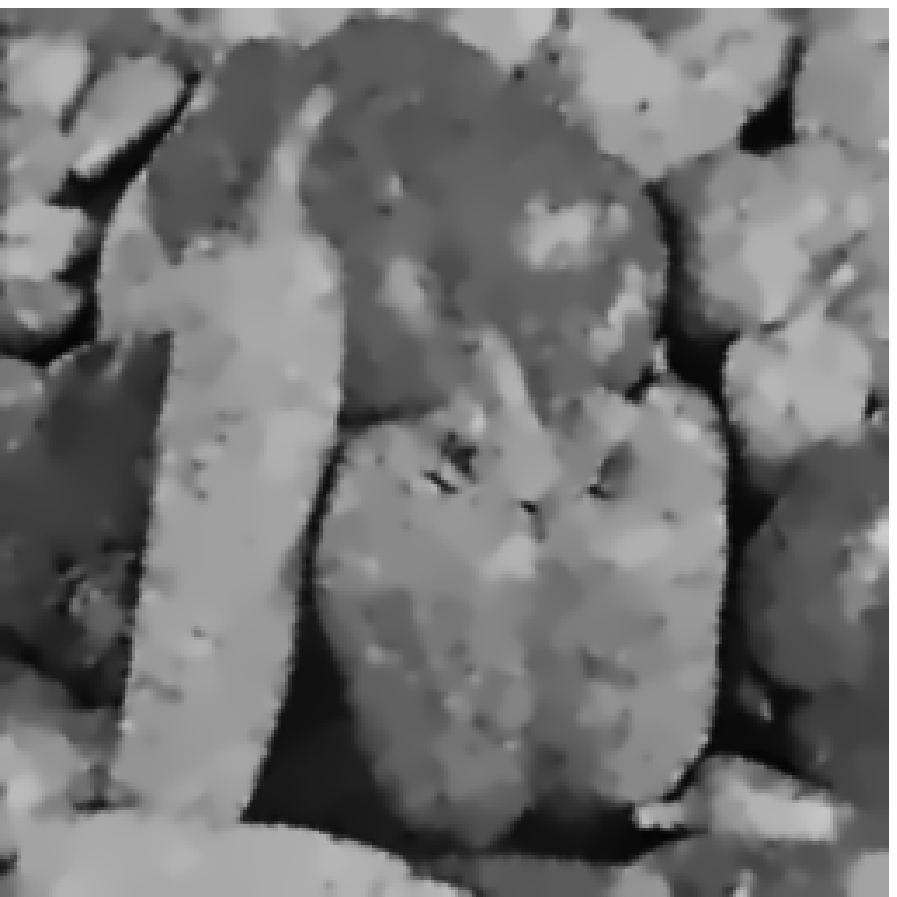}               
                \caption{Proposed}
                \label{fig2:peppers_3_tdm_hp}
       \end{subfigure}
     
\caption{Image corrupted with speckle look L=3 and restored by different models.}\label{peppers_3}
\end{figure}

\begin{figure}
       \centering

        \begin{subfigure}[b]{0.24\textwidth}           
                \includegraphics[scale=0.4]{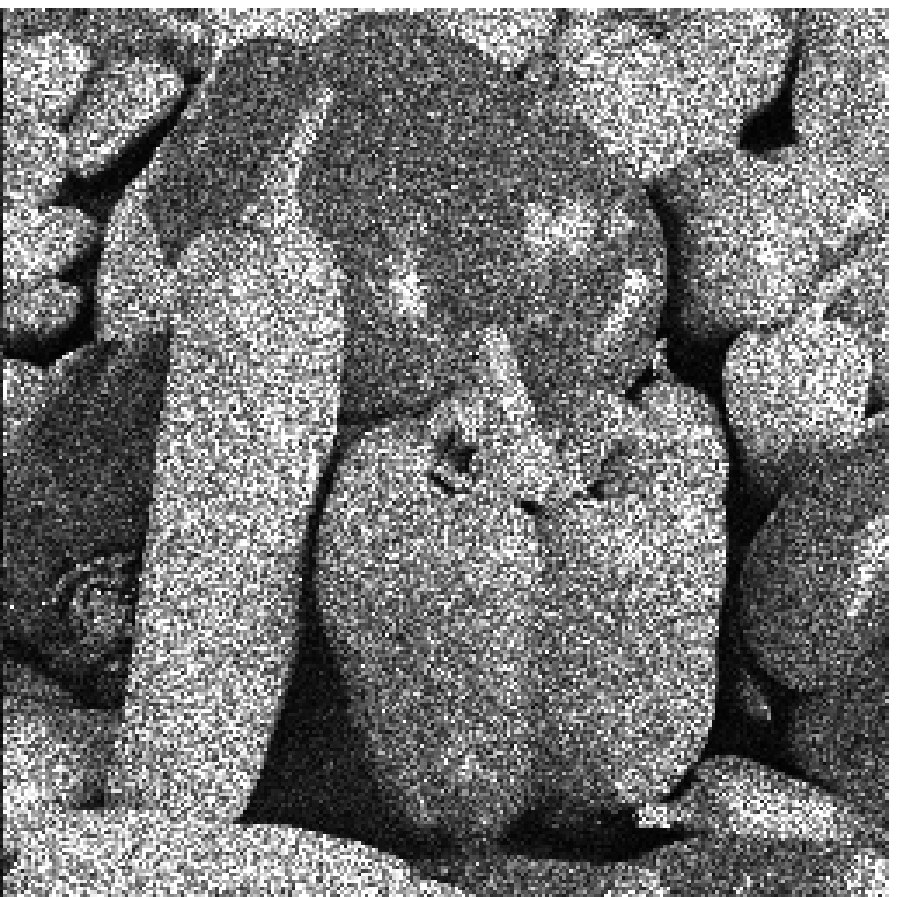}               
                \caption{Noisy}
                \label{fig3:peppers_5}
       \end{subfigure}% 
   \begin{subfigure}[b]{0.24\textwidth}           
                \includegraphics[scale=0.4]{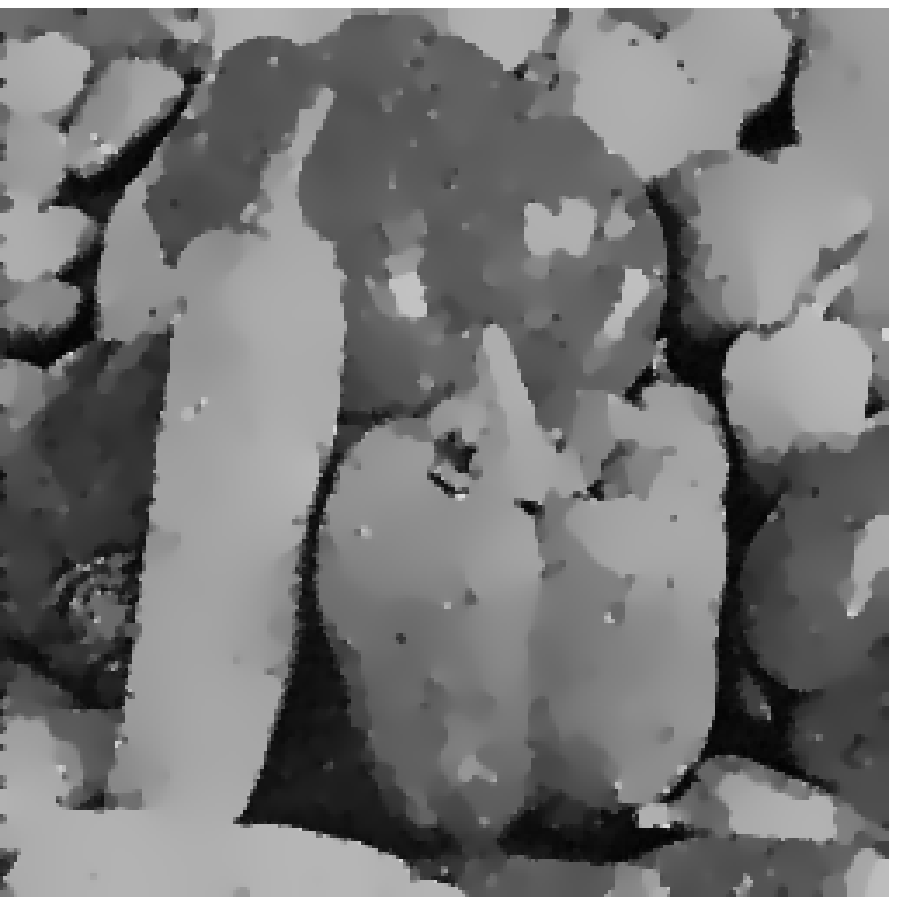}               
                \caption{Shan}
                \label{fig3:peppers_5_shan}
       \end{subfigure}%
         \begin{subfigure}[b]{0.24\textwidth}           
                \includegraphics[scale=0.4]{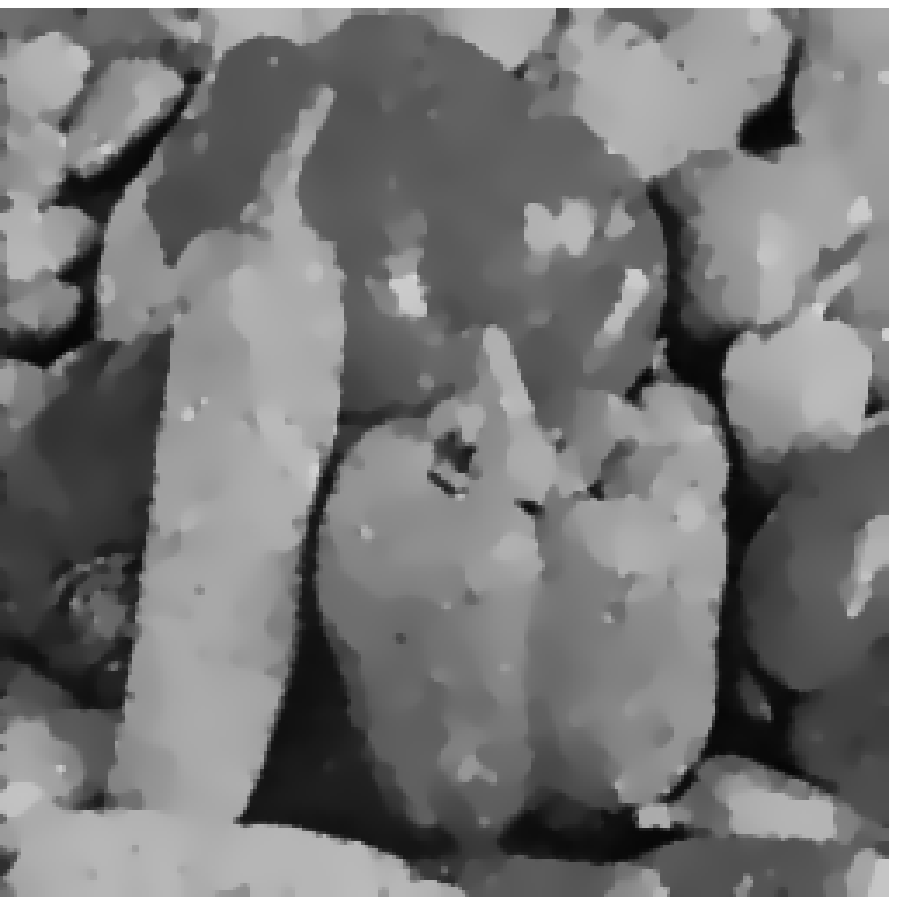}               
                \caption{TDM}
                \label{fig3:peppers_5_tdm}
       \end{subfigure}% 
        \begin{subfigure}[b]{0.24\textwidth}           
                \includegraphics[scale=0.4]{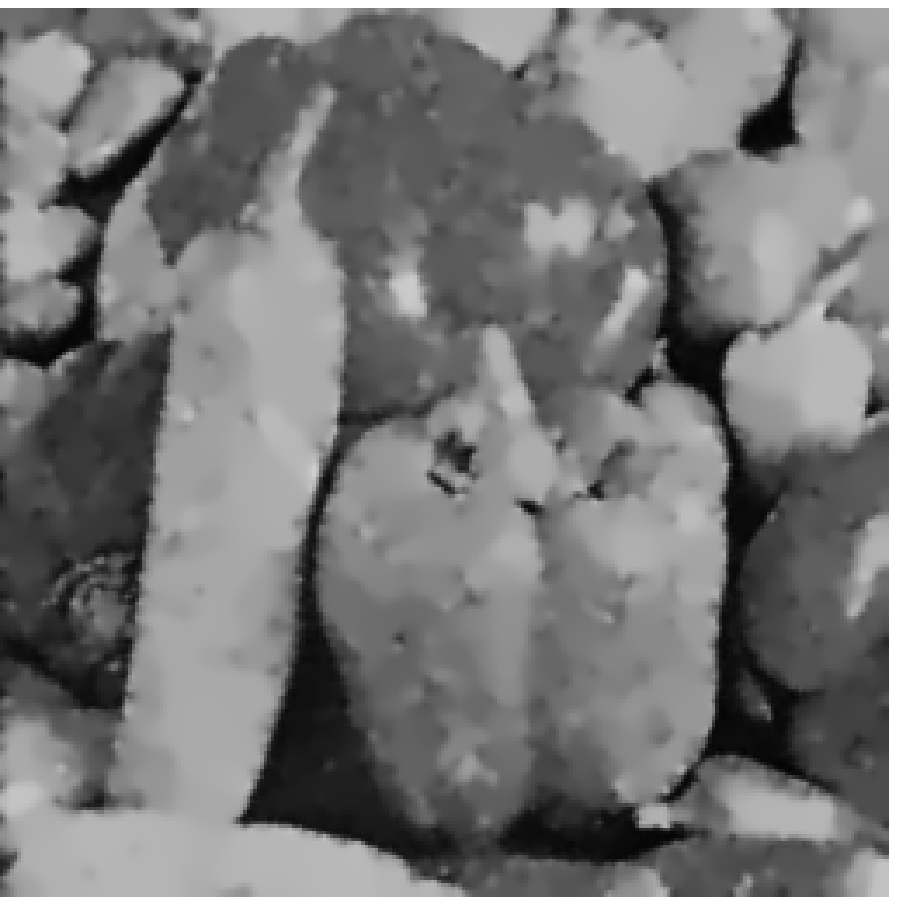}               
                \caption{Proposed}
                \label{fig3:peppers_5_tdm_hp}
       \end{subfigure}
     
\caption{Image corrupted with speckle look L=5 and restored by different models.}\label{peppers_5}
\end{figure}

In figures \ref{peppers_1}-\ref{peppers_5}, we represent the reconstructed results of a Peppers image (Natural Image) which is corrupted by speckle noise with  $L=\lbrace1,3,5\rbrace$. From figure \ref{peppers_1}, we see that the present model leave less speckle than the other two models.

\begin{figure}
       \centering
     
       \begin{subfigure}[b]{0.23\textwidth}           
                \includegraphics[scale=0.34]{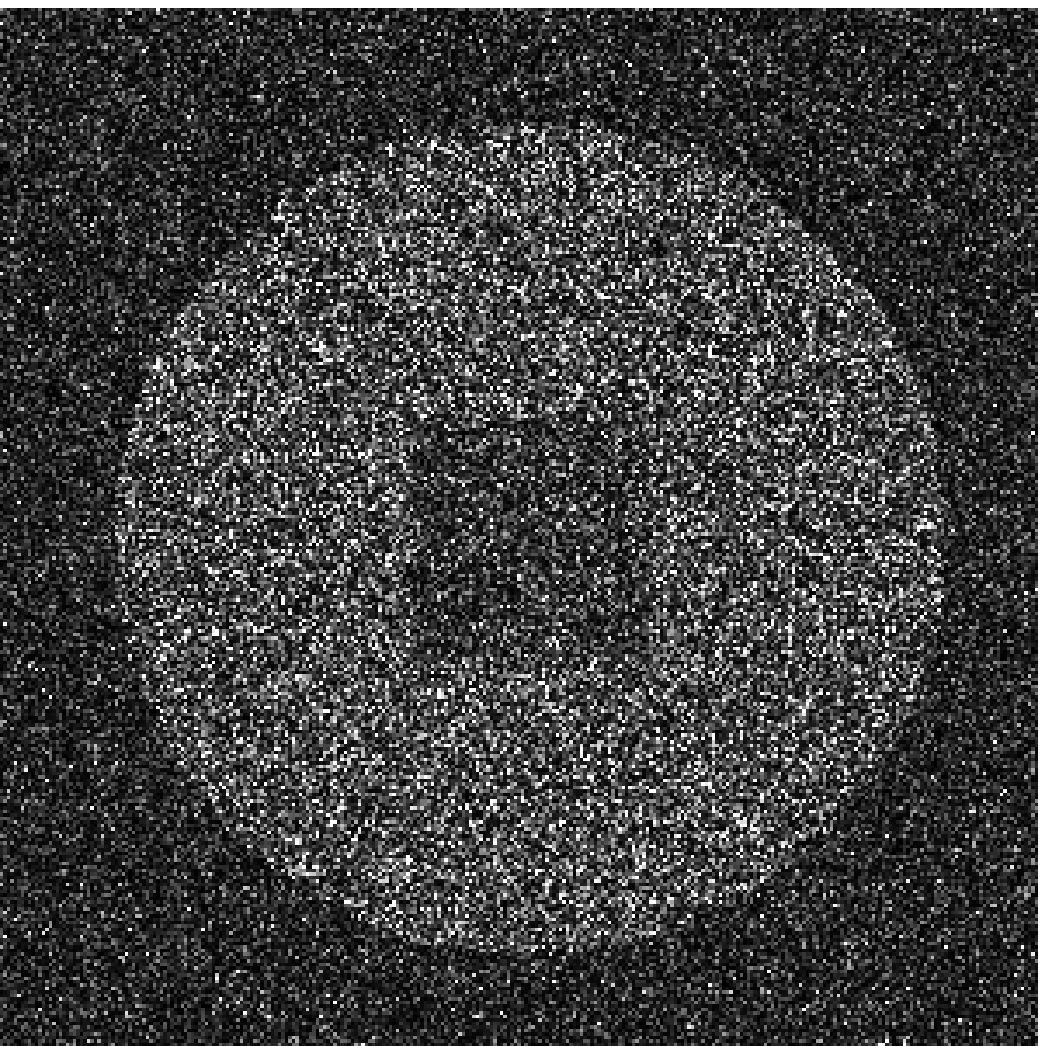}               
                \caption{Noisy}
                \label{fig1:circle_1}
       \end{subfigure}%
         \begin{subfigure}[b]{0.23\textwidth}           
                \includegraphics[scale=0.34]{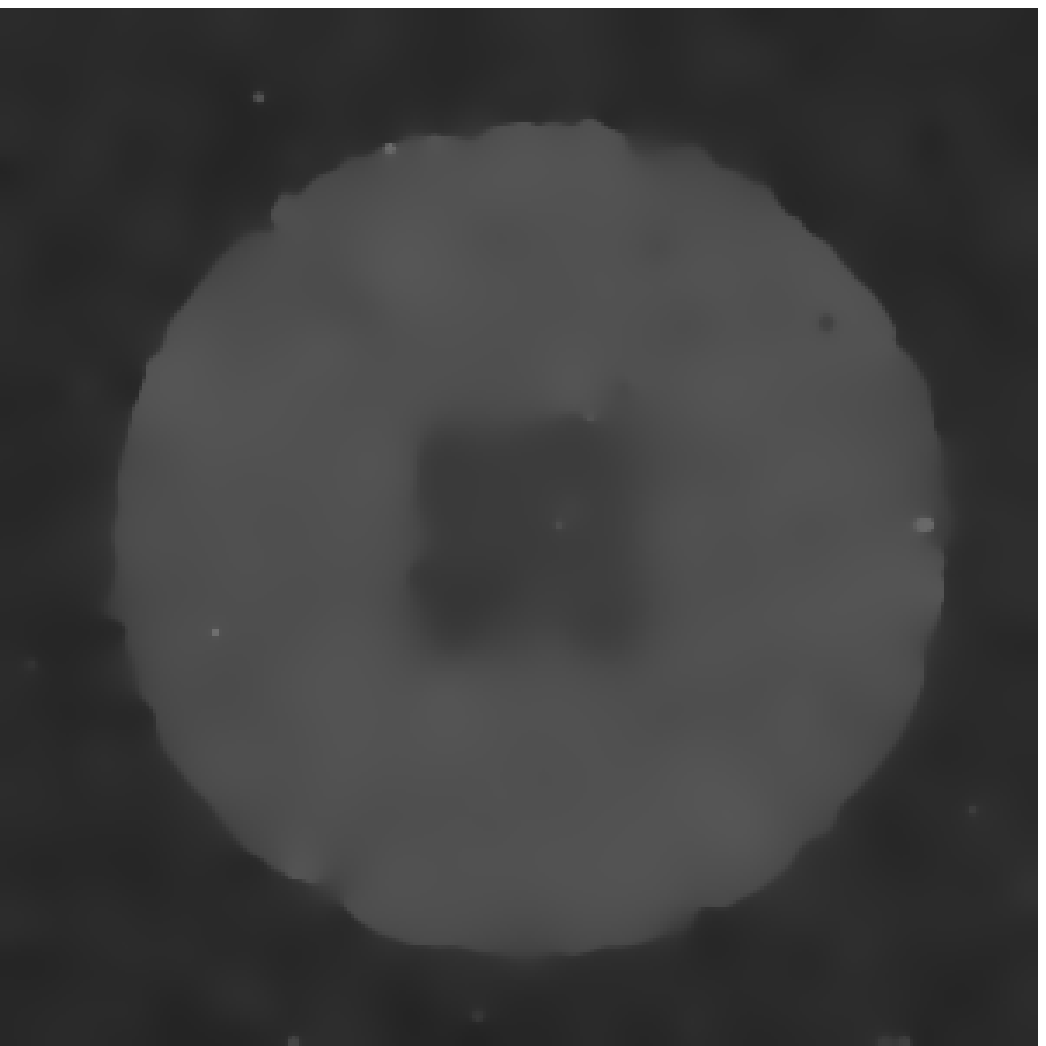}               
                \caption{Shan}
                \label{fig1:circle_1_shan}
       \end{subfigure}% 
        \begin{subfigure}[b]{0.23\textwidth}           
                \includegraphics[scale=0.34]{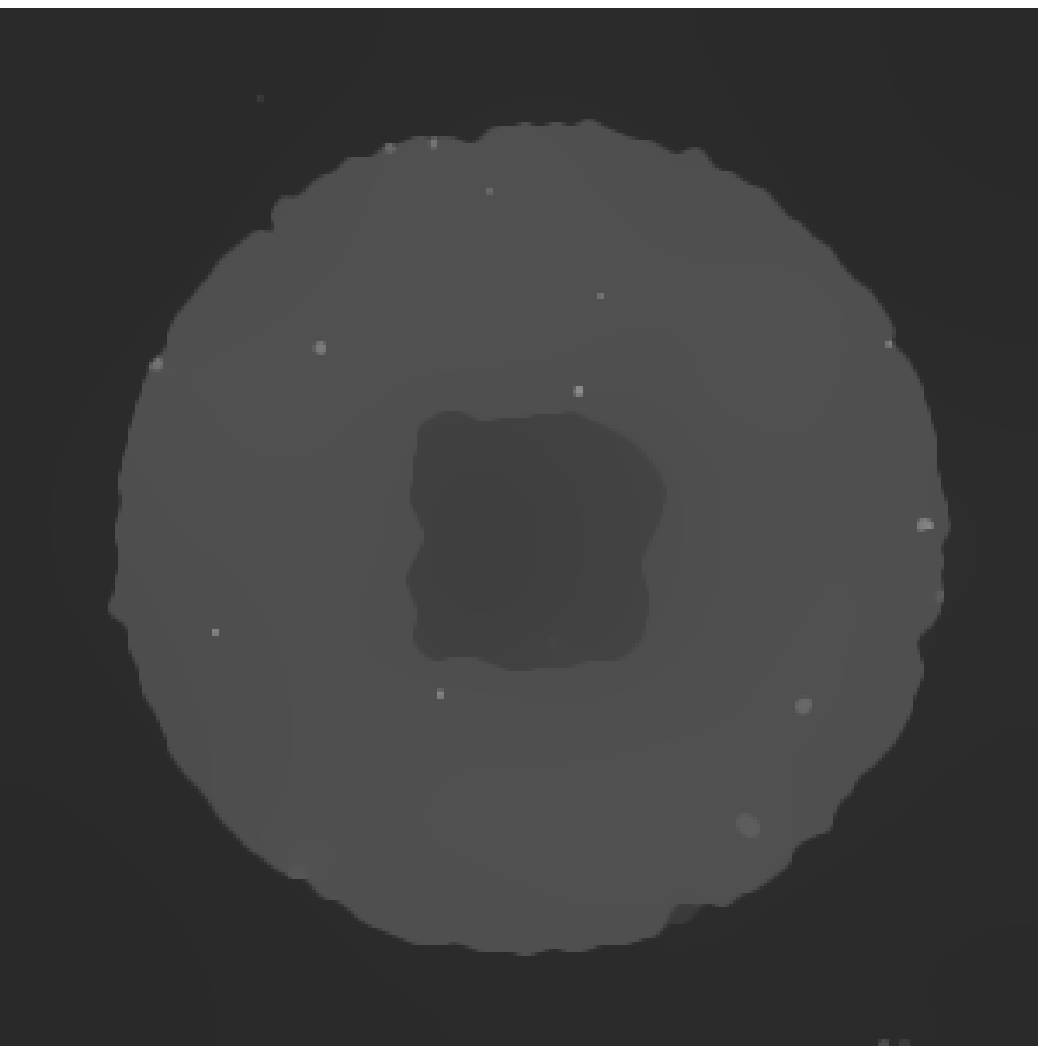}               
                \caption{TDM}
                \label{fig1:circle_1_tdm}
       \end{subfigure}
             \begin{subfigure}[b]{0.23\textwidth}           
                \includegraphics[scale=0.34]{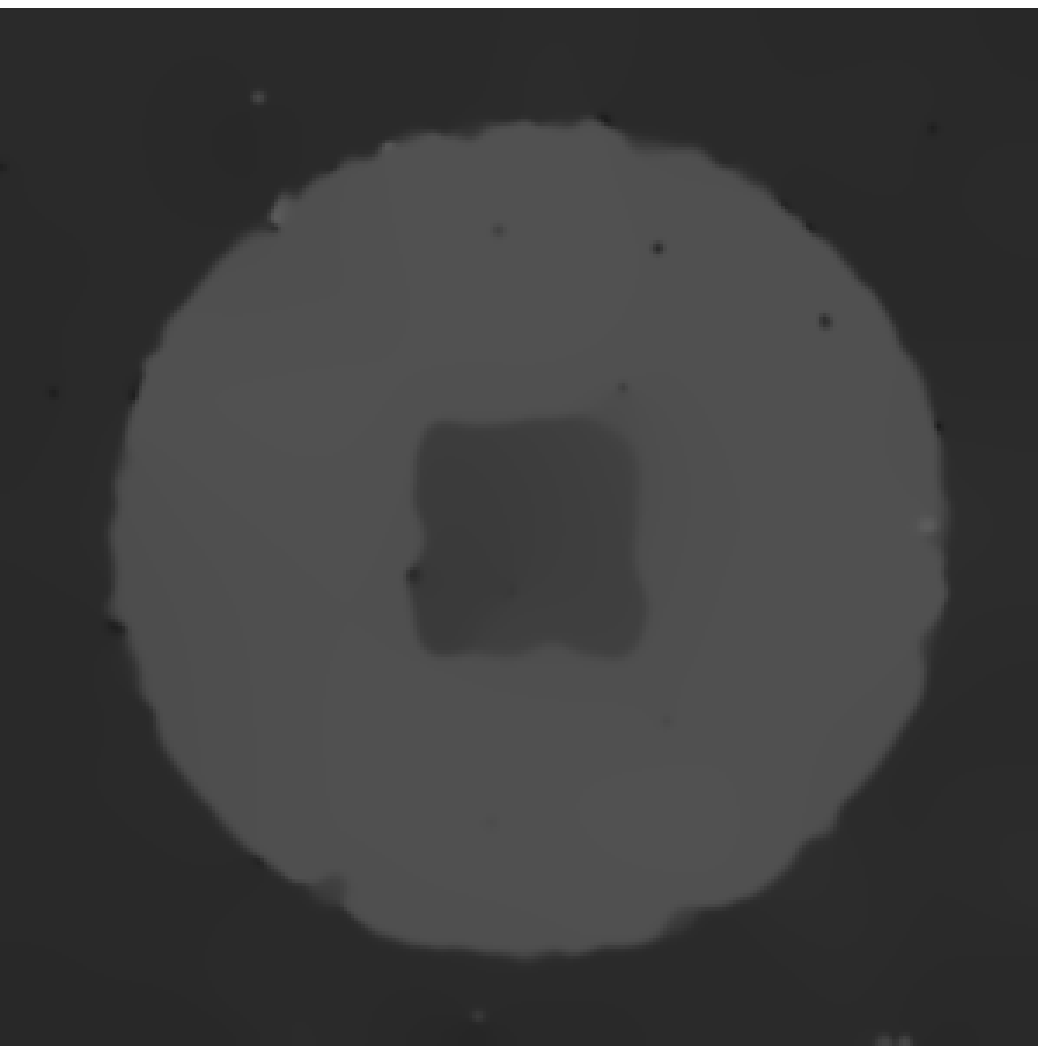}               
                \caption{Proposed}
                \label{fig1:circle_1_tdm_hp}
       \end{subfigure}
     
 \caption{Image corrupted with speckle look L=1 and restored by different models. }\label{circle_1}
\end{figure}

\begin{figure}
       \centering

      \begin{subfigure}[b]{0.23\textwidth}           
                \includegraphics[scale=0.34]{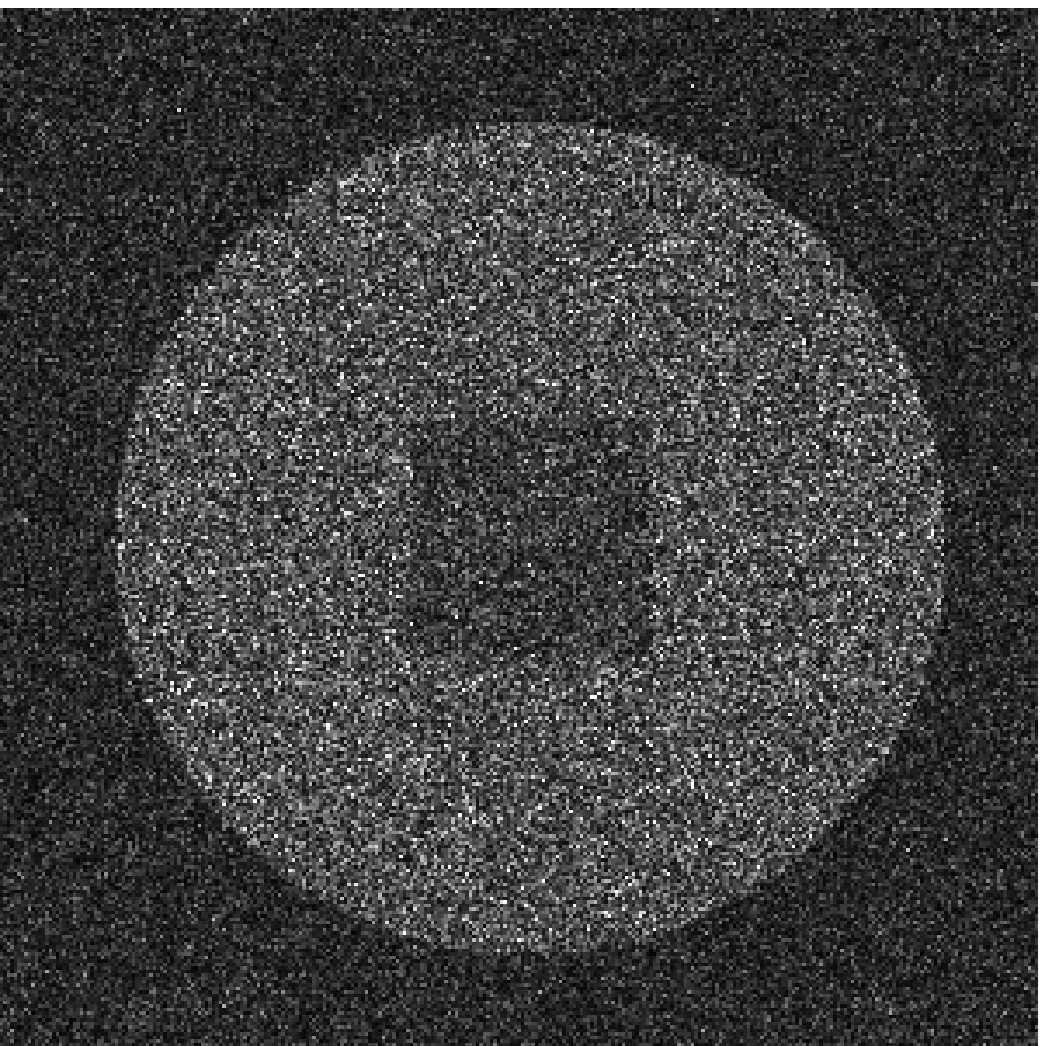}               
                \caption{Noisy}
                \label{fig2:circle_3}
       \end{subfigure}%
         \begin{subfigure}[b]{0.23\textwidth}           
                \includegraphics[scale=0.34]{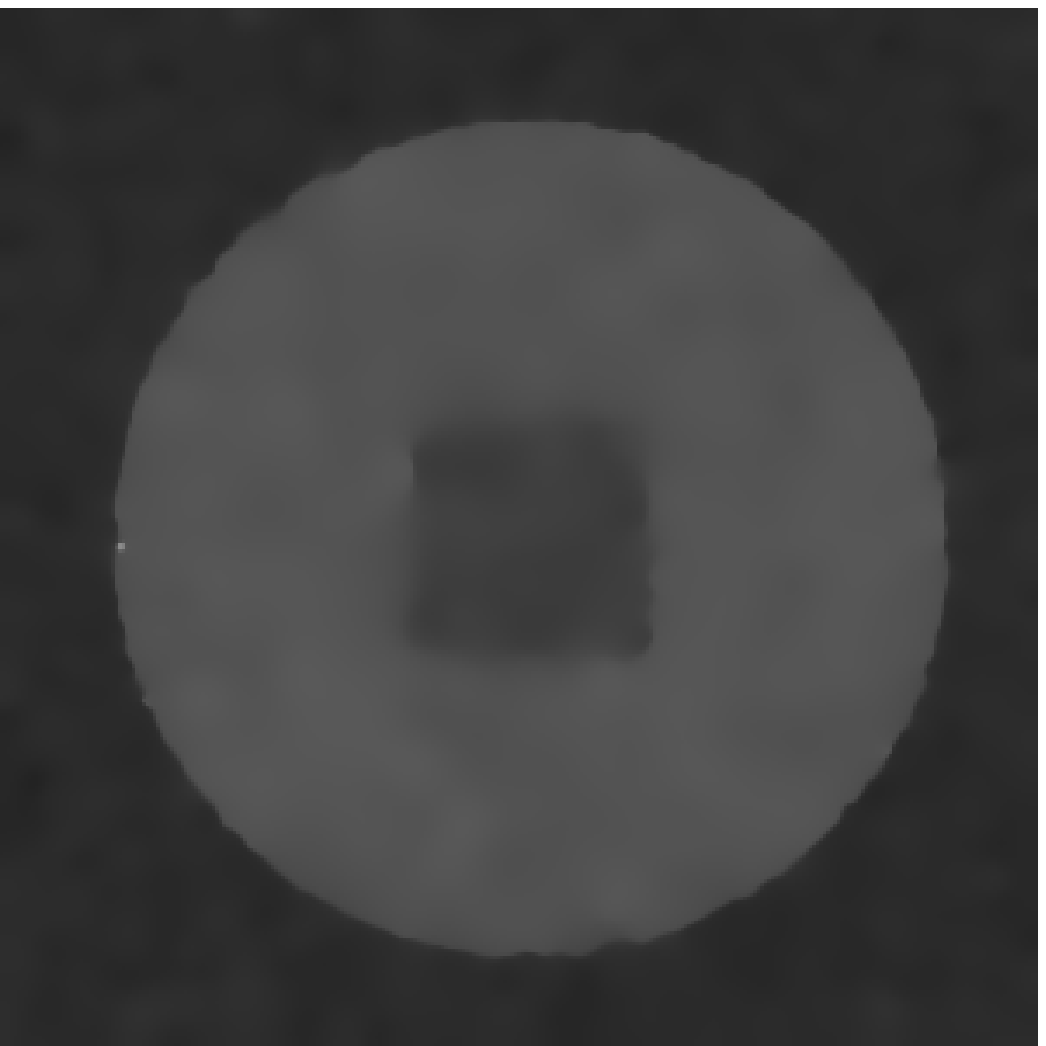}               
                \caption{Shan}
                \label{fig2:circle_3_shan}
       \end{subfigure}% 
        \begin{subfigure}[b]{0.23\textwidth}           
                \includegraphics[scale=0.34]{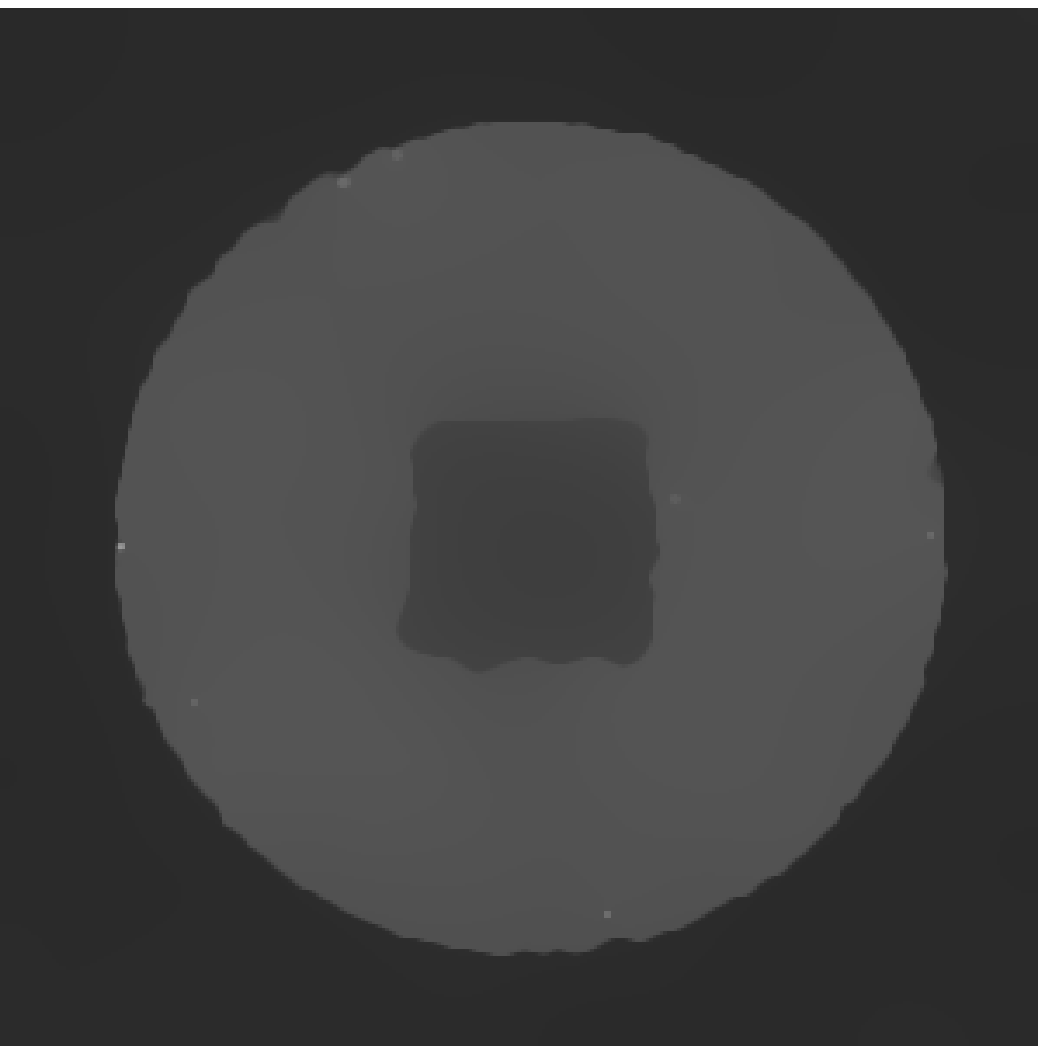}               
                \caption{TDM}
                \label{fig2:circle_3_tdm}
       \end{subfigure}
             \begin{subfigure}[b]{0.23\textwidth}           
                \includegraphics[scale=0.34]{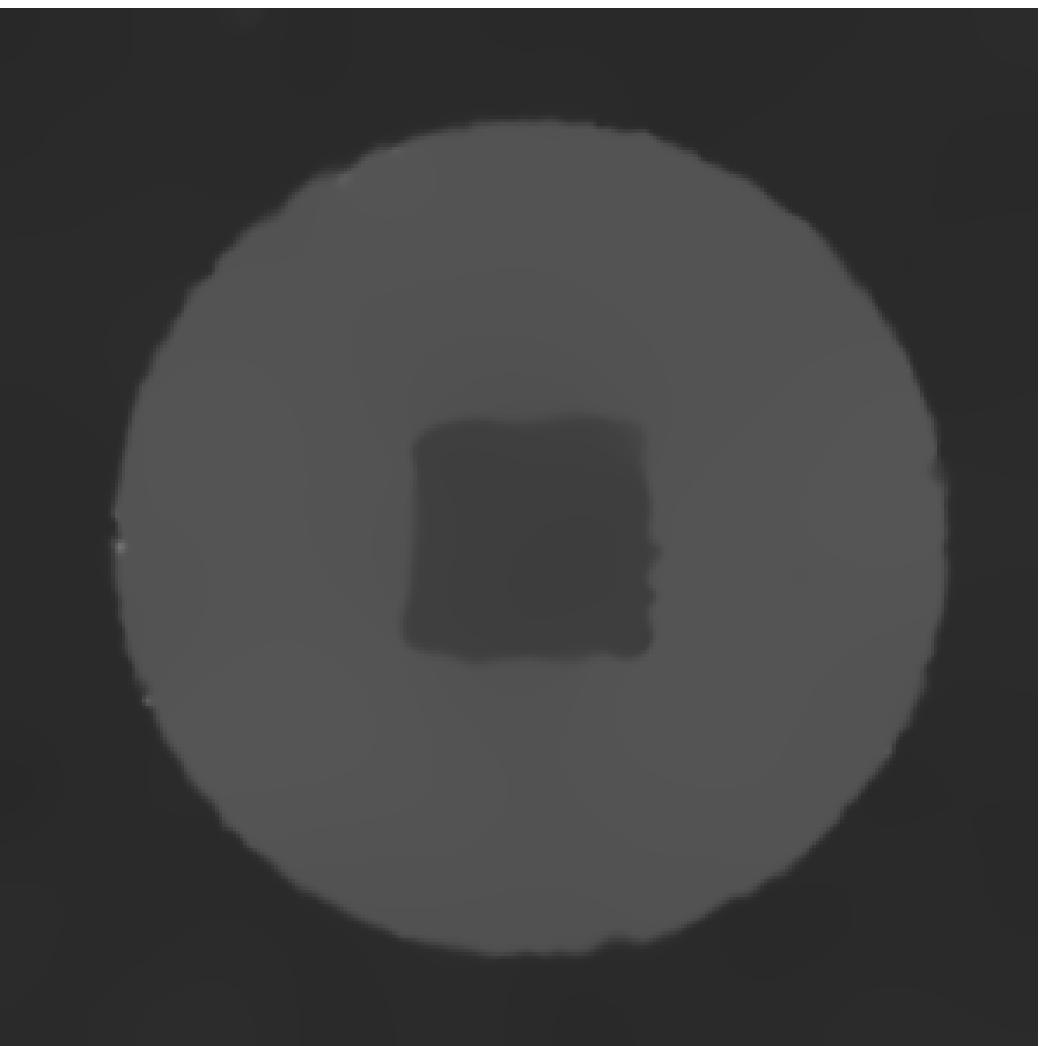}               
                \caption{Proposed}
                \label{fig2:circle_3_tdm_hp}
       \end{subfigure}
     
 \caption{Image corrupted with speckle look L=3 and restored by different models. }\label{circle_3}
\end{figure}

\begin{figure}
       \centering

         \begin{subfigure}[b]{0.23\textwidth}           
                \includegraphics[scale=0.34]{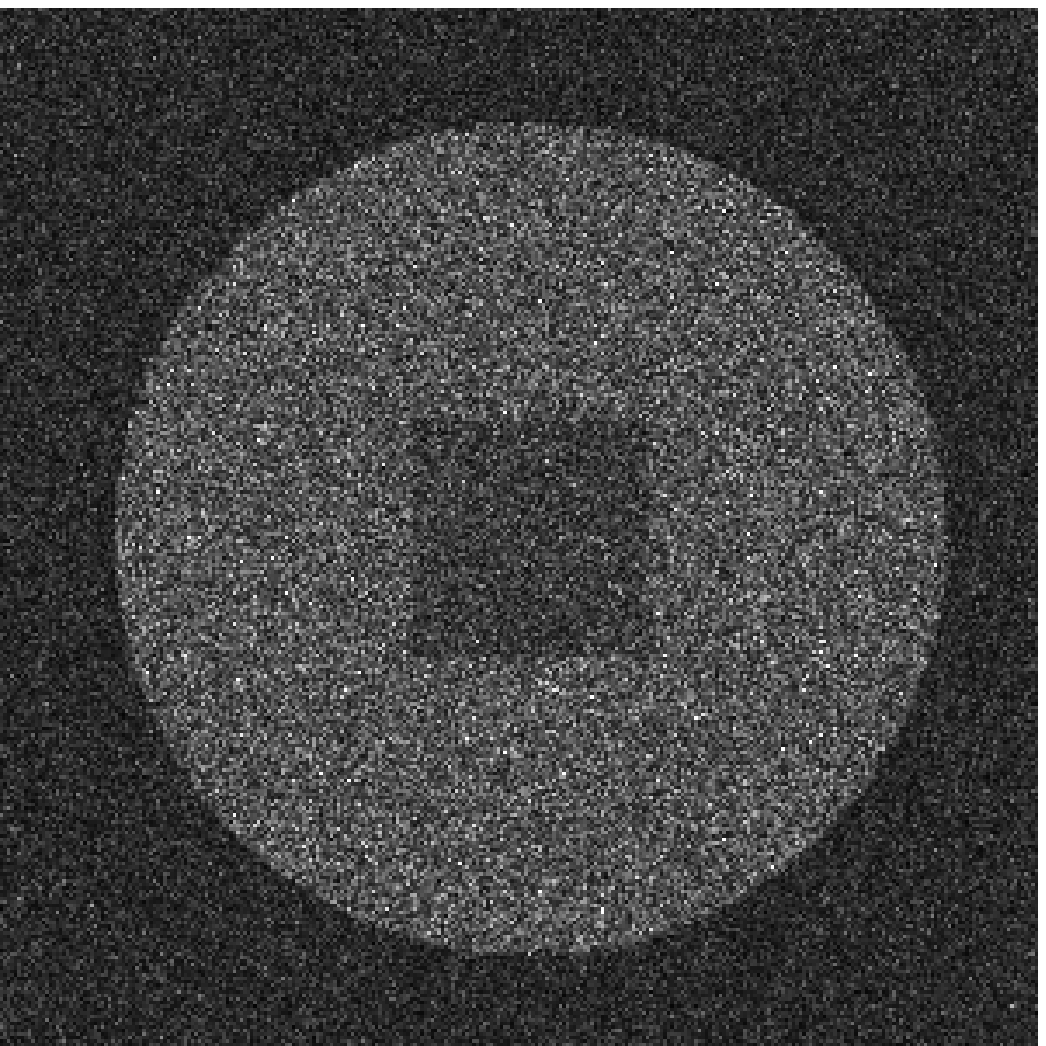}               
                \caption{Noisy}
                \label{fig3:circle_5}
       \end{subfigure}%
         \begin{subfigure}[b]{0.23\textwidth}           
                \includegraphics[scale=0.34]{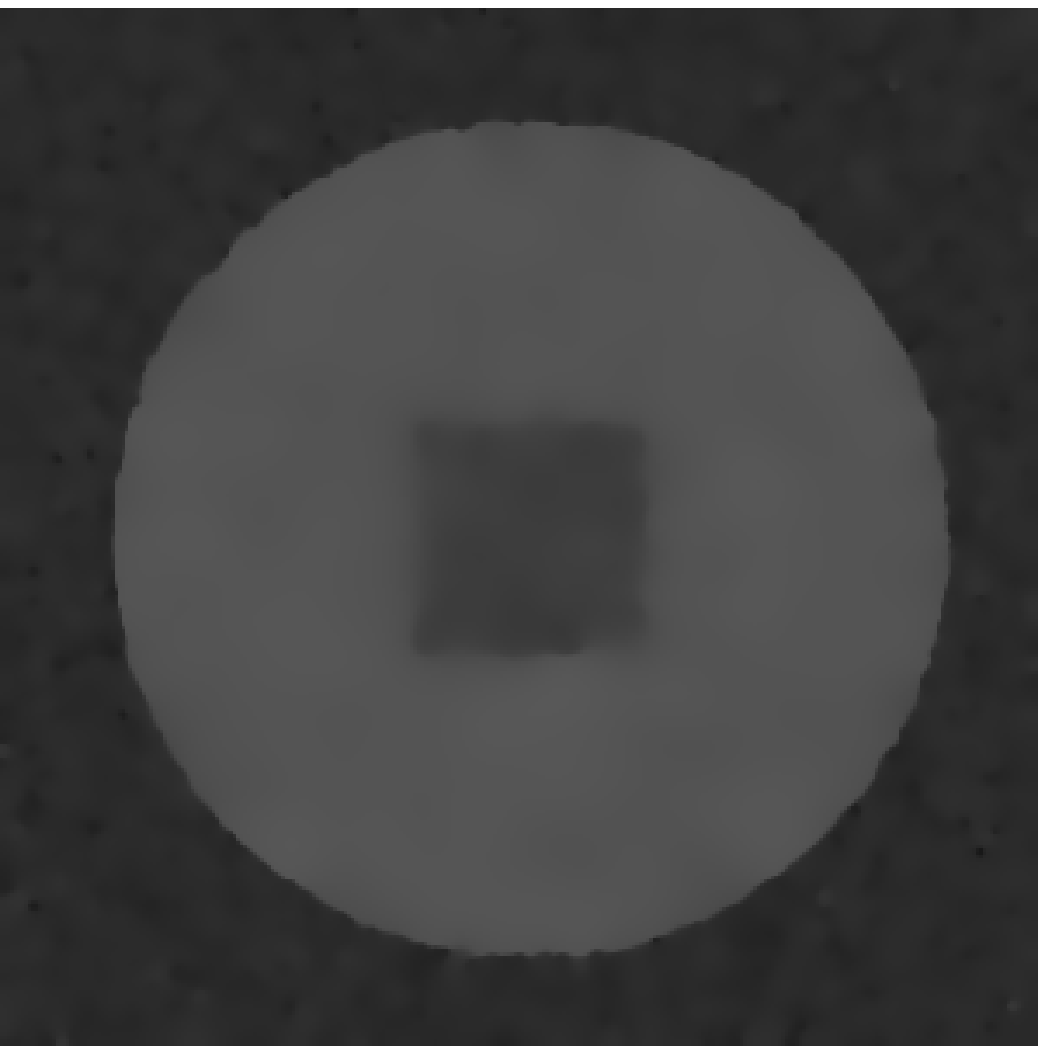}               
                \caption{Shan}
                \label{fig3:circle_5_shan}
       \end{subfigure}% 
        \begin{subfigure}[b]{0.23\textwidth}           
                \includegraphics[scale=0.34]{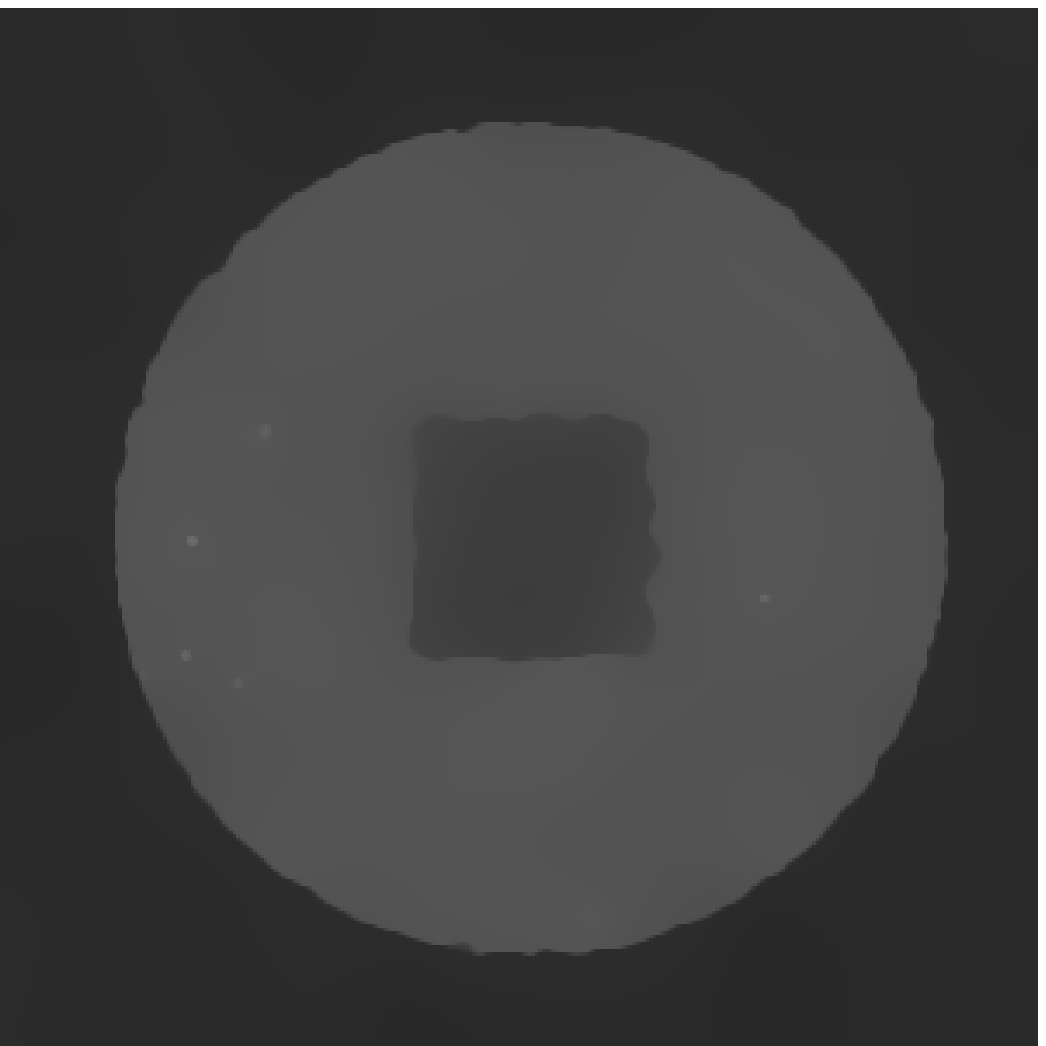}               
                \caption{TDM}
                \label{fig3:circle_5_tdm}
       \end{subfigure}
             \begin{subfigure}[b]{0.23\textwidth}           
                \includegraphics[scale=0.34]{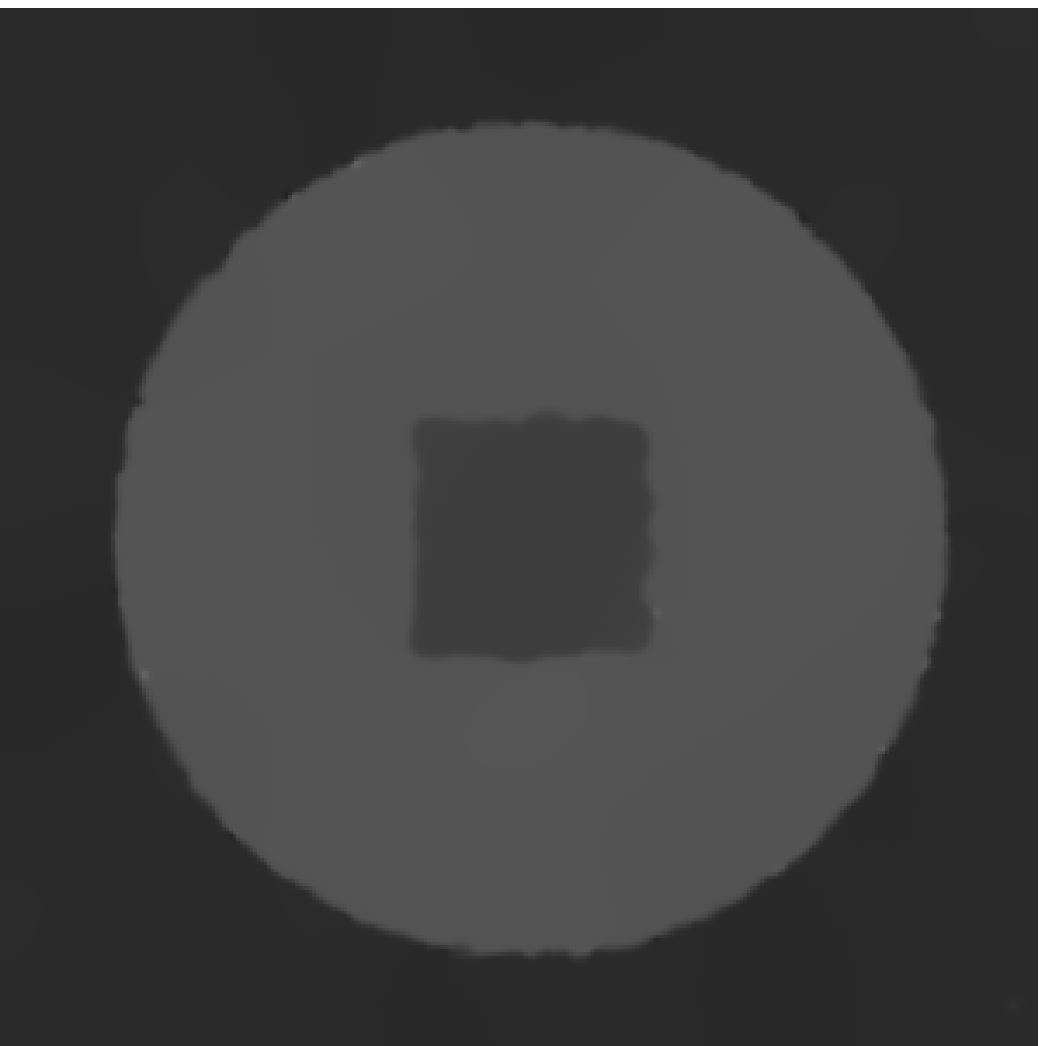}               
                \caption{Proposed}
                \label{fig3:circle_5_tdm_hp}
       \end{subfigure}
     
 \caption{Image corrupted with speckle look L=5 and restored by different models. }\label{circle_5}
\end{figure}

\begin{figure}
       \centering
       
         \begin{subfigure}[b]{0.35\textwidth}           
           \includegraphics[scale=0.4]{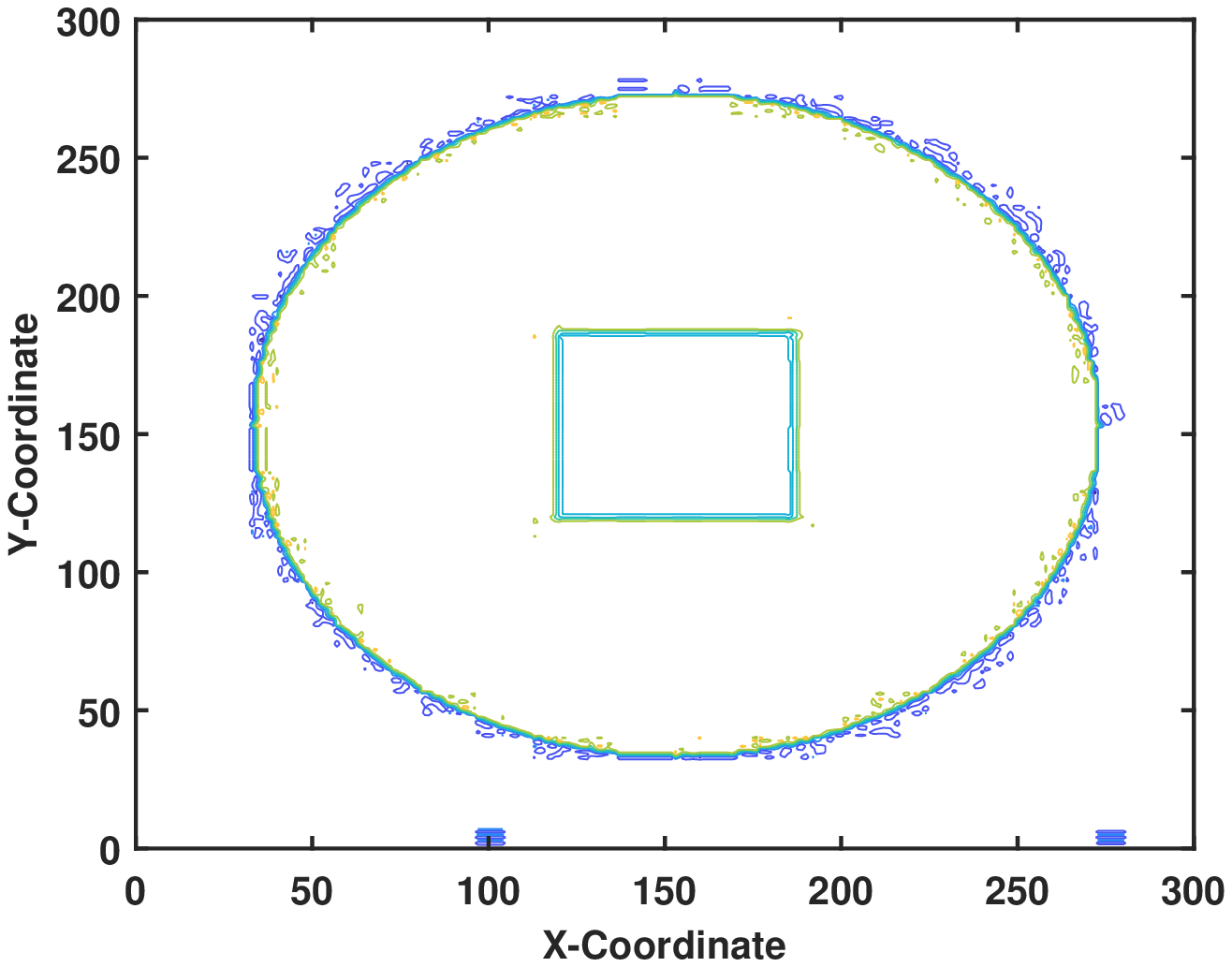}           
                \caption{Original}
                \label{fig:3a_cont}
        \end{subfigure}
       \begin{subfigure}[b]{0.35\textwidth}           
                \includegraphics[scale=0.4]{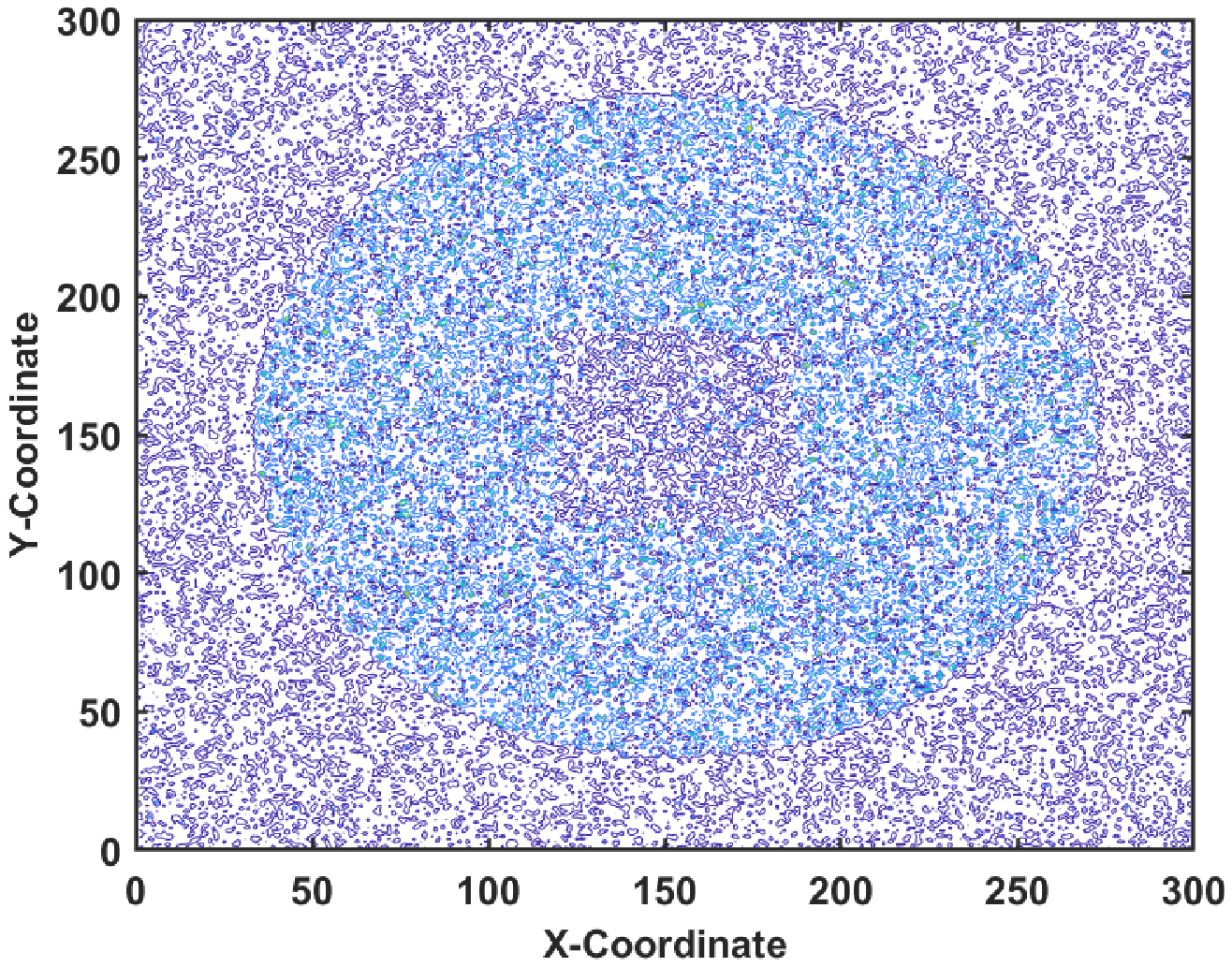}               
                \caption{Noisy}
                \label{fig:3b_cont}
       \end{subfigure}% 
    
       \begin{subfigure}[b]{0.35\textwidth}           
                \includegraphics[scale=0.4]{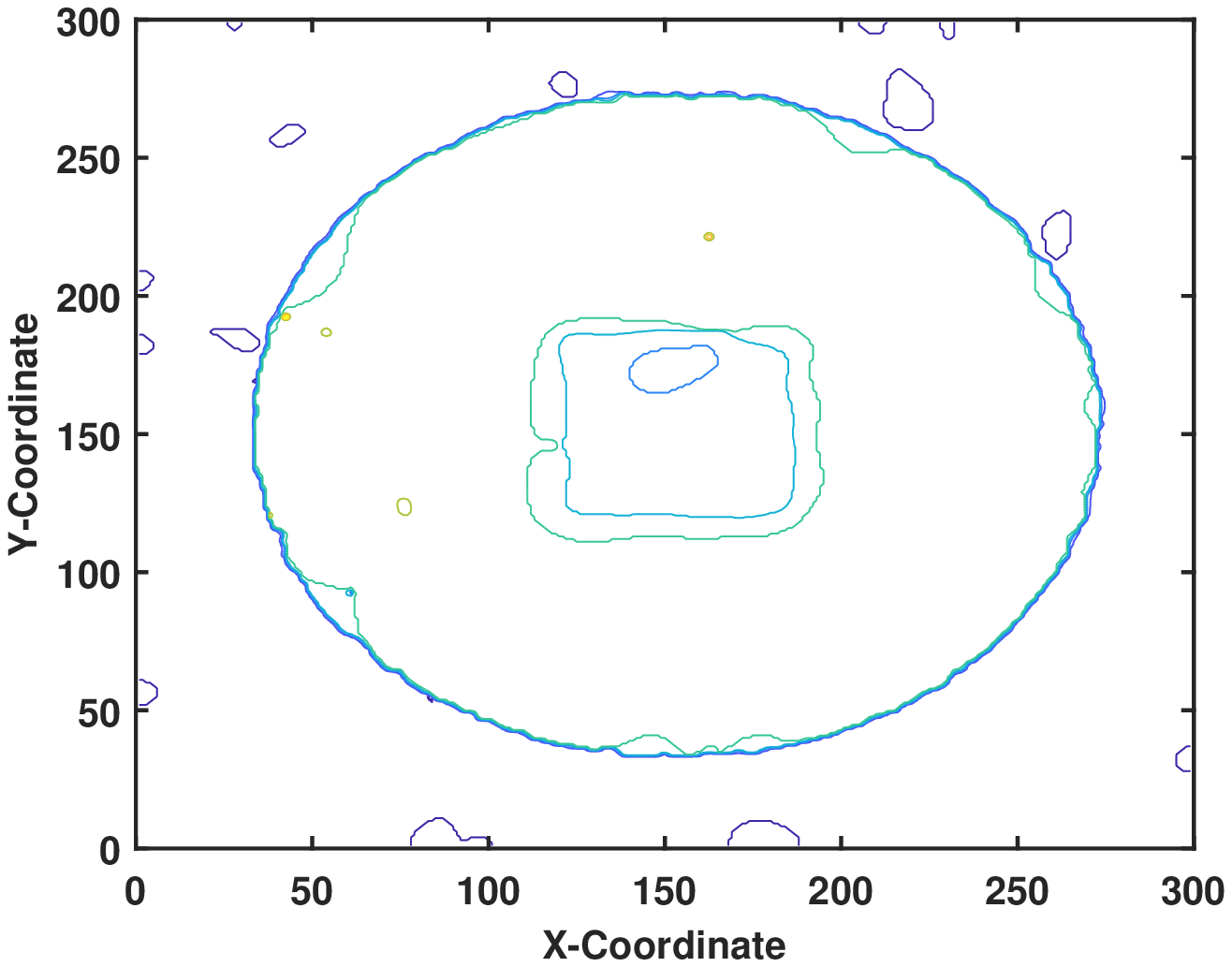}               
                \caption{Shan}
                \label{fig:3c_cont}
        \end{subfigure}% 
        \begin{subfigure}[b]{0.35\textwidth}           
                \includegraphics[scale=0.4]{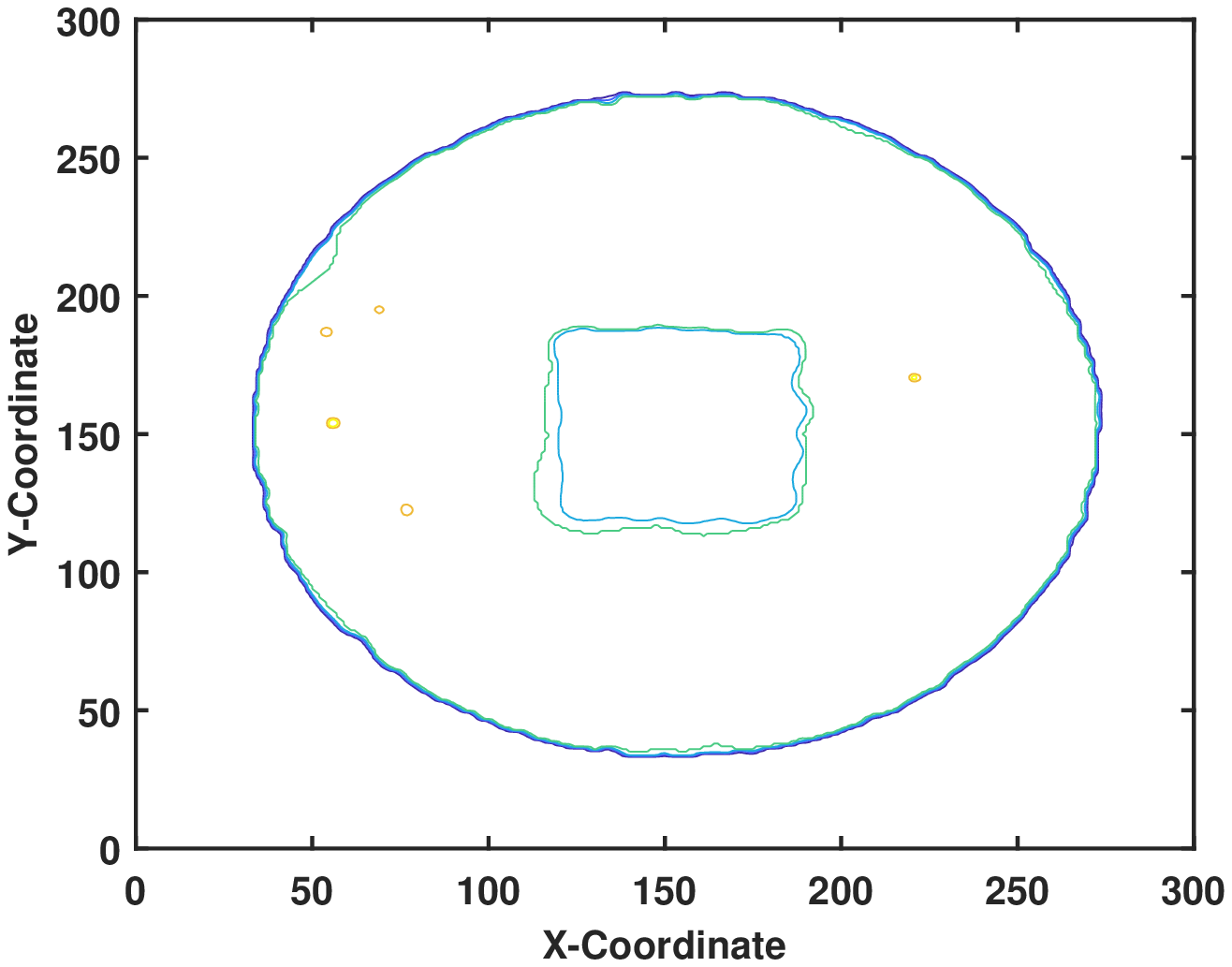}               
                \caption{TDM}
                \label{fig:3d_cont}
         \end{subfigure}% 
        \begin{subfigure}[b]{0.35\textwidth}           
                \includegraphics[scale=0.4]{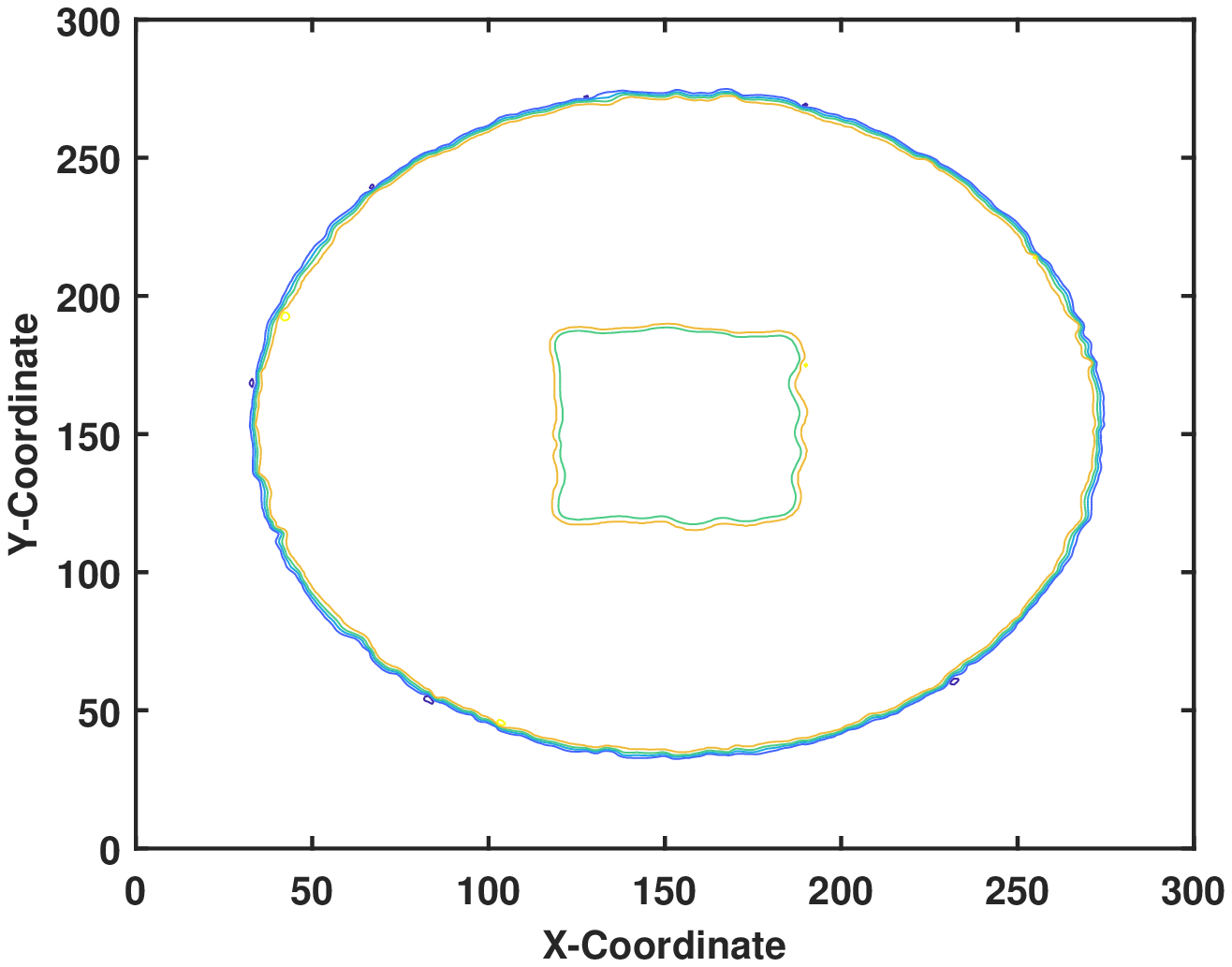}               
                \caption{Proposed}
                \label{fig:3e_cont}
        \end{subfigure}% 
       
\caption{Contour plots of the restored images in figure \ref{circle_5}.}\label{circle_5_cont}
\end{figure}

\begin{figure}
       \centering
       
           \begin{subfigure}[b]{0.35\textwidth}           
           \includegraphics[scale=0.4]{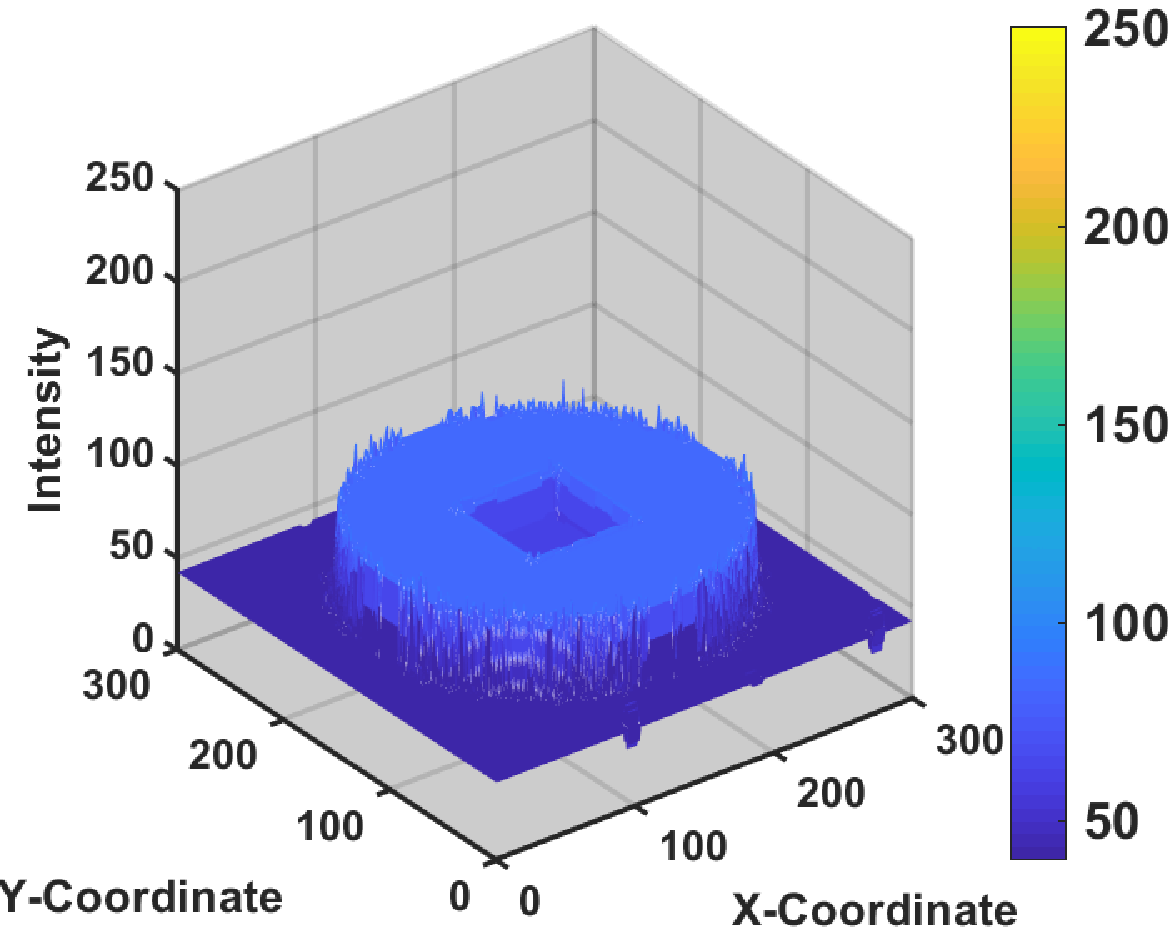}           
                \caption{Original}
                \label{fig:2a_3d}
        \end{subfigure}
        \begin{subfigure}[b]{0.35\textwidth}           
                \includegraphics[scale=0.4]{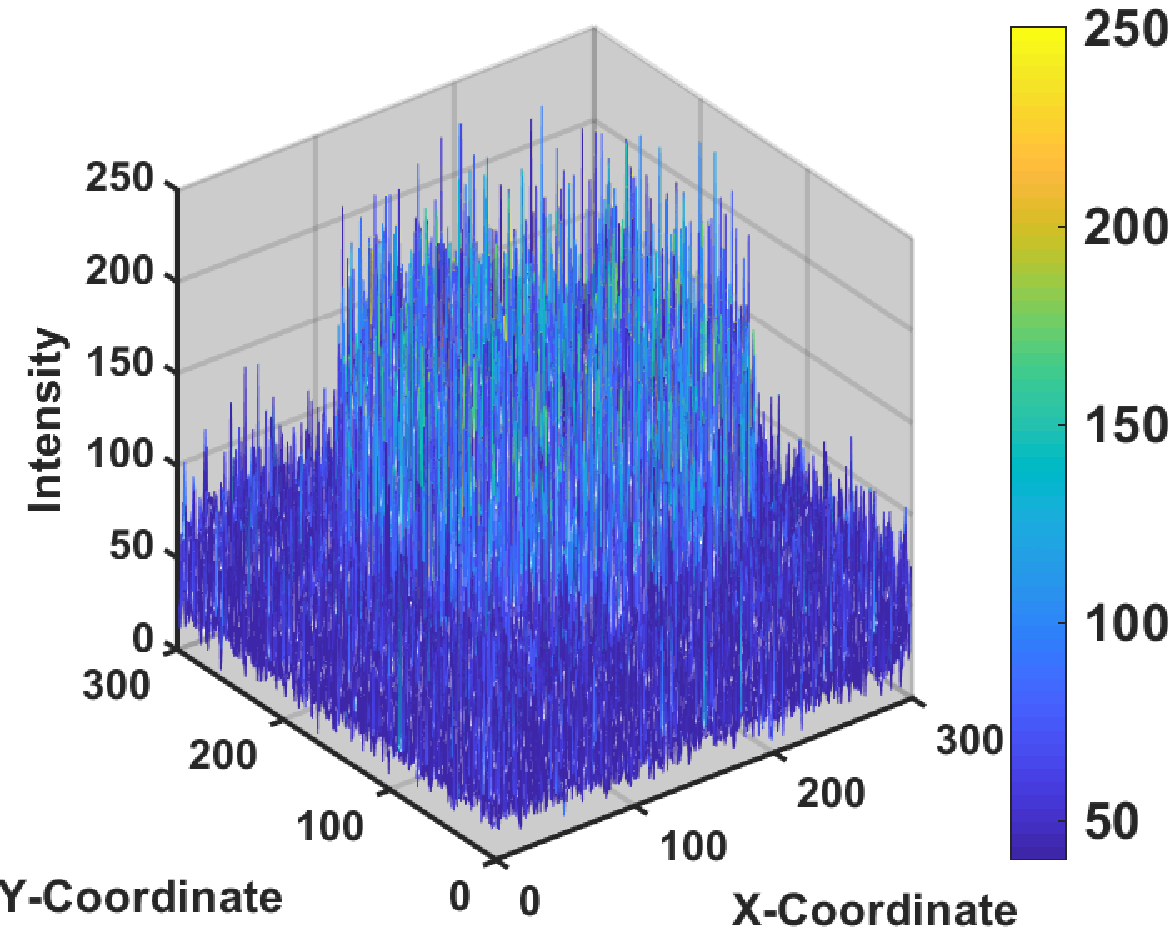}               
                \caption{Noisy}
                \label{fig:2b_3d}
       \end{subfigure}% 
       
       \begin{subfigure}[b]{0.35\textwidth}           
                \includegraphics[scale=0.4]{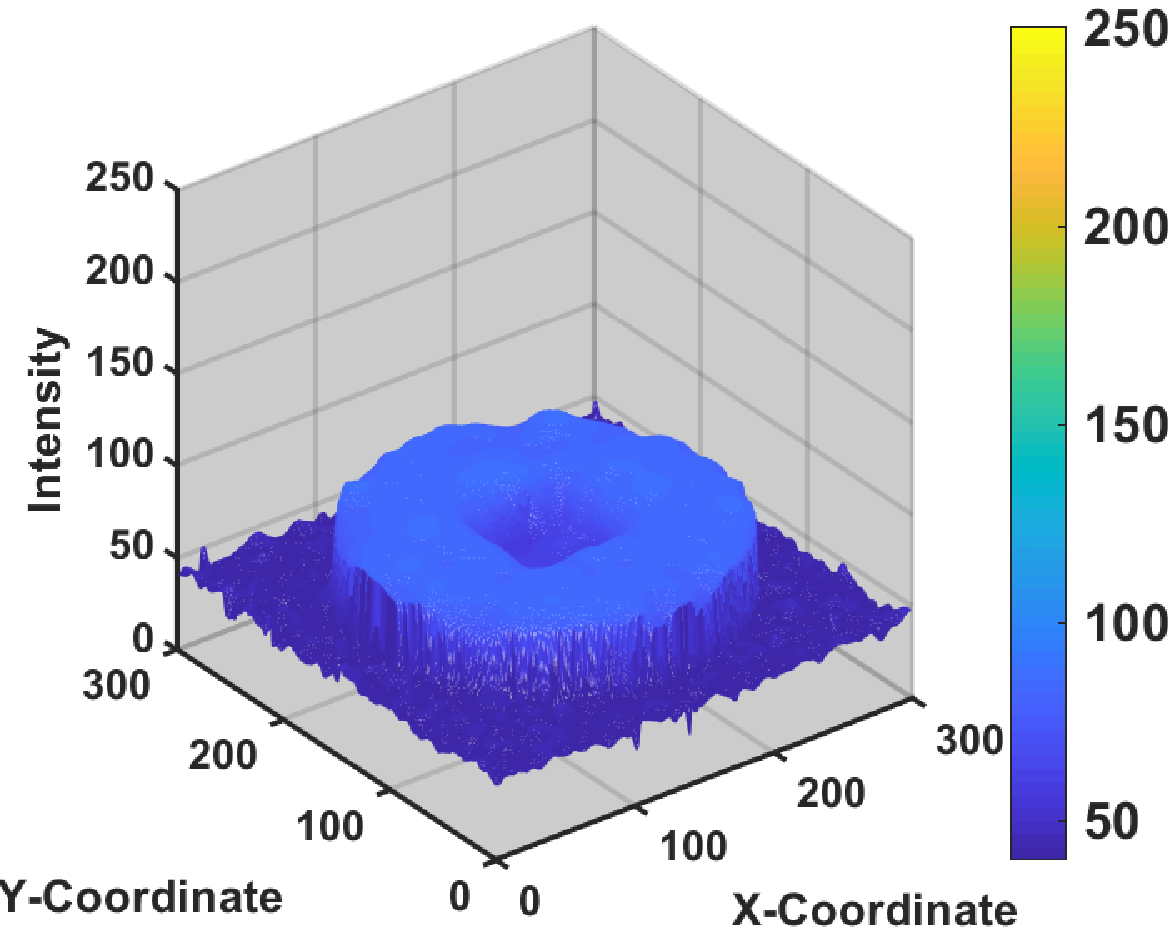}               
                \caption{Shan}
                \label{fig:2c_3d}
        \end{subfigure}% 
        \begin{subfigure}[b]{0.35\textwidth}           
                \includegraphics[scale=0.4]{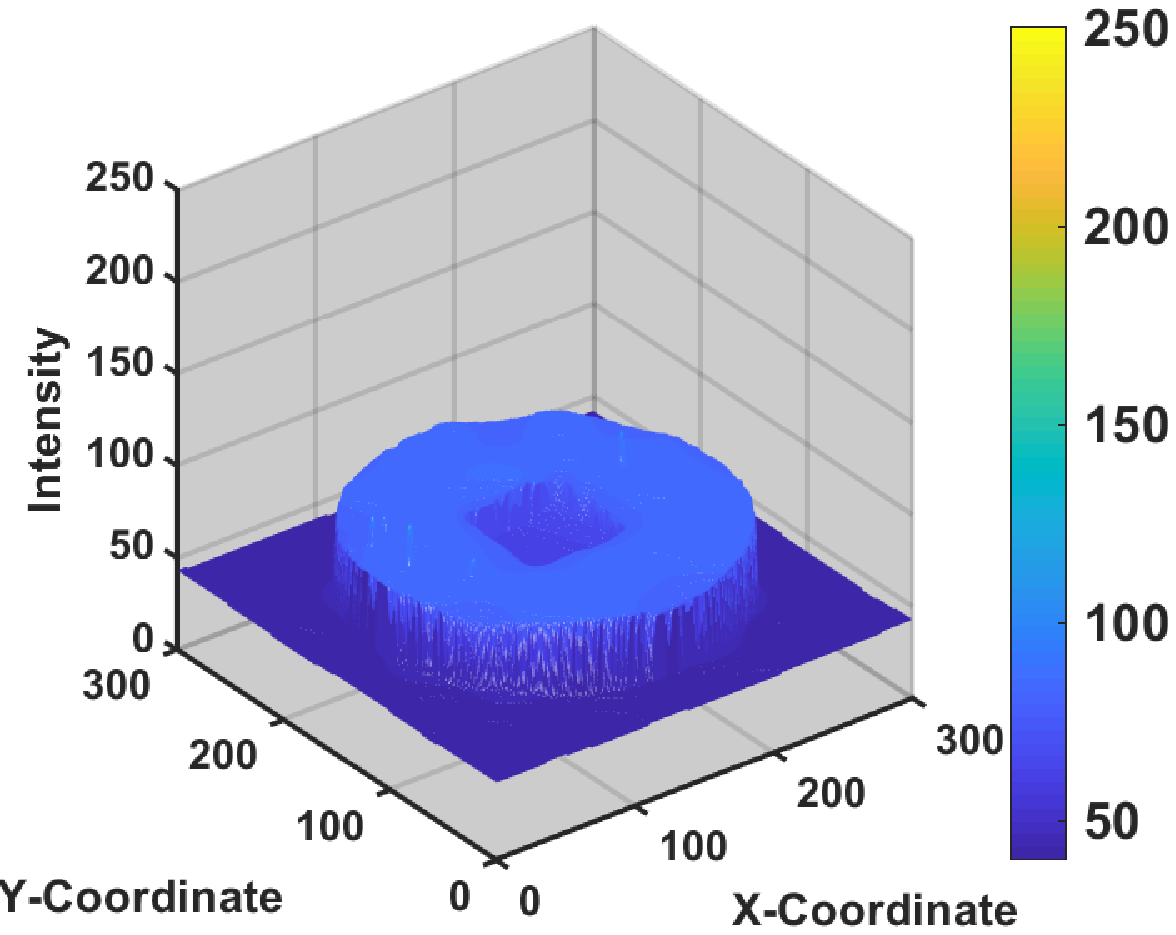}               
                \caption{TDM}
                \label{fig:2d_3d}
         \end{subfigure}% 
        \begin{subfigure}[b]{0.35\textwidth}           
                \includegraphics[scale=0.4]{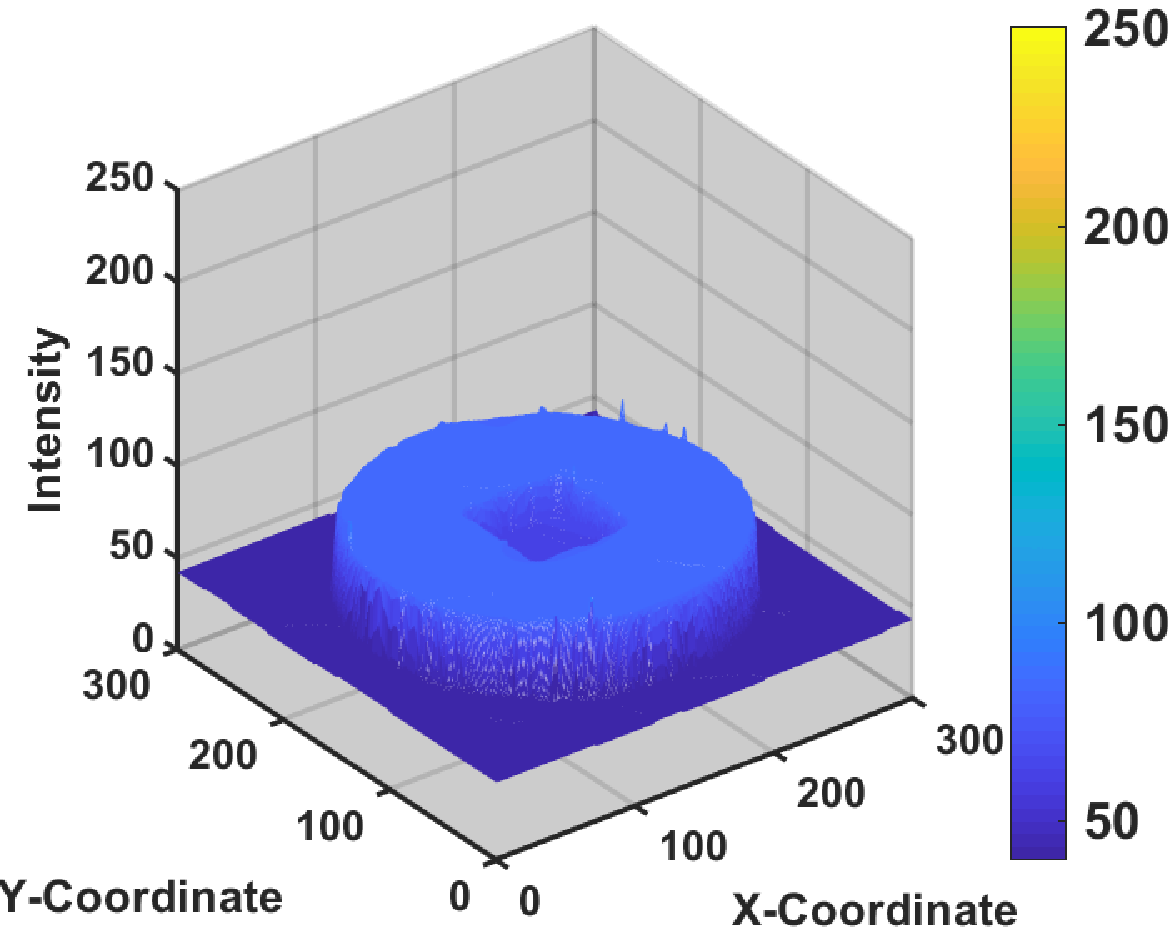}               
                \caption{Proposed}
                \label{fig:2e_3d}
        \end{subfigure}% 
       
\caption{ 3D surface plots of the restored images in figure \ref{circle_5}.}\label{circle_5_3d}
\end{figure}

To check the more reconstruction ability of the present model in figures \ref{circle_1}- \ref{circle_5_3d} illustrate the qualitative results of a Circle image (Synthetic Image) which is corrupted by speckle noise with  $L = \lbrace 1,3,5 \rbrace$.

In the figures \ref{circle_1}- \ref{circle_5} we demonstrate the despeckled images, and in the figures \ref{circle_5_cont}-\ref{circle_5_3d} we illustrate the contour maps and 3D surface plots when the image is corrupted by $L=5$. One can observe that from the contour maps, and 3D surface plots, Shan and TDM models left some speckles in the homogeneous regions, but the present model produces fewer artifacts with better edge preservation.

Computational Values of PSNR and MSSIM are presented in the Table \ref{tab:psnr_ssim_parameter}.  The highest values of PSNR and MSSIM for each noise level clearly shows that the suggested model is better than the other two models.

Conclusively, summarize the quantitative and qualitative results, we can confirm that the present model performance better than the other discussed models.
%===================================================================================================================

%===================================================================================================================
%==============================================================================

%==============================================================================

%============================================================================================================

%======================================================================================================

%======================================================================================================
%======================================================================================================
%======================================================================================================
%======================================================================================================
%======================================================================================================

%=========================================
%============================================================================================================

%===================================================================================================================
\begin{center}
\begin{table}
\caption{Left table: Comparison of MSSIM and PSNR values of despeckled images. Right table: Parameter values for the numerical experiments.}
\label{tab:psnr_ssim_parameter}
\begin{tabular}{ll}
\scalebox{0.6}{
\begin{tabular}[t]{llcccccccccccccc}
\toprule
    \multirow{2}[8]{*}{Image} & \multirow{2}[8]{*}{$L$}  & \multicolumn{2}{c}{Shan Model\cite{shan2019multiplicative}}  & \multicolumn{2}{c}{ TDM \cite{majeetdm2019speckle} }  & \multicolumn{2}{c}{Proposed Model} \\
\cmidrule(r){3-4}
		\cmidrule(r){5-6}
		\cmidrule(r){7-8}
		\cmidrule(r){9-10}  
		\cmidrule(r){11-12} 
		\cmidrule(r){13-14} 
		
&        & \multicolumn{1}{c}{MSSIM} & \multicolumn{1}{c}{PSNR} & \multicolumn{1}{c}{MSSIM} & \multicolumn{1}{c}{PSNR} & \multicolumn{1}{c}{MSSIM}  & \multicolumn{1}{c}{PSNR}  \\
    \midrule
Circle     & 1             & 0.9581 & 34.29  & 0.9643 & 34.69 & \textbf{0.9651} & \textbf{35.43}\\
           & 3             & 0.9734 & 38.09  & 0.9771 & 39.52 & \textbf{0.9782} & \textbf{40.00}\\
           & 5             & 0.9764 & 39.35  & 0.9805 & 40.72 & \textbf{0.9810} & \textbf{41.23}\\
           &               &        &        &        &       &                 &                \\
Texture    & 1             & 0.8125 & 27.61  & 0.8355 & 27.71 & \textbf{0.8360} & \textbf{27.83}\\
           & 3             & 0.8782 & 30.83  & 0.8925 & 31.05 & \textbf{0.8967} & \textbf{31.31}\\
           & 5             & 0.8979 & 31.87  & 0.9109 & 32.15 & \textbf{0.9164} & \textbf{32.47}\\
           &               &        &        &         &      &                 &                \\
Peppers    & 1             & 0.5827 & 17.56  & 0.5895 & 17.64 & \textbf{0.5905} & \textbf{17.83}\\
           & 3             & 0.7018 & 22.46  & 0.7019 & 22.55 & \textbf{0.7061} & \textbf{22.85}\\
           & 5             & 0.7155 & 23.73  & 0.7334 & 24.22 & \textbf{0.7395} & \textbf{24.67}\\
           &               &        &        &         &      &                 &                \\       
\midrule        
%\bottomrule
\end{tabular}
}
&
\scalebox{0.7}{
\begin{tabular}[t]{llccccccccccc}
\toprule
    \multirow{2}[4]{*}{Image} & \multirow{2}[4]{*}{ $L$ } &\multicolumn{2}{c}{Shan\cite{shan2019multiplicative}}& \multicolumn{3}{c}{TDM \cite{majeetdm2019speckle}} & \multicolumn{5}{c}{Proposed}      \\
    
        \cmidrule(r){3-4}
		\cmidrule(r){5-7}
		\cmidrule(r){8-12}
		    &   & \multicolumn{1}{c}{$\alpha$} & \multicolumn{1}{c}{$\beta$} & \multicolumn{1}{c}{$\gamma$} & \multicolumn{1}{c}{$\nu$} & \multicolumn{1}{c}{$K$} &  \multicolumn{1}{c}{$\gamma$} & \multicolumn{1}{c}{$\alpha$} & \multicolumn{1}{c}{$\beta$} & \multicolumn{1}{c}{$\iota$} & \multicolumn{2}{c}{$\nu$}\\
  \midrule
Circle      & 1          & 1.5    & 2       & 10      & 1       & 1      & 1   &  1.5  & 1.8  & 2.5  & 0.1  \\
            & 3          & 1.5    & 2       & 10      & 1       & 1      & 2   &  1.7  & 2    & 2.5  & 0.1  \\
            & 5          & 2      & 2.25    & 5       & 1       & 1      & 2   &  1.7  & 2.2  & 2    & 0.1  \\
\midrule
Texture     & 1          & 1.5    & 1.8     & 2       & 1.5     & 2      & 1   &  2    & 1    & 3    & 0.1  \\
            & 3          & 1.8    & 2       & 2       & 1.5     & 2      & 2   &  2    & 1    & 3    & 0.1  \\
            & 5          & 1.8    & 2       & 2       & 1.5     & 2      & 5   &  2.5  & 1    & 3    & 0.1  \\
\midrule
Peppers     & 1          & 1.5    & 2       & 1       & 1.5     & 2      & 1   &  2    & 1    & 3    & 0.1  \\
            & 3          & 1.5    & 2       & 1       & 1.5     & 2      & 2   &  2    & 1    & 3    & 0.1  \\
            & 5          & 2      & 2.4     & 2       & 1.5     & 2      & 2   &  2.5  & 1    & 3    & 0.1  \\
\bottomrule
\end{tabular}
}
\end{tabular}
\end{table}
\end{center}

%====================================================================================================================
\section{Conclusion}
\label{sec:Conclusion}
In this work, we present a non-linear hyperbolic-parabolic coupled system applied to image despeckling. Such a improve method preserves the image characteristics in the noise removal process.
To the best of our knowledge, coupled hyperbolic-parabolic PDE based model has not been used before for image speckle reduction. Moreover, we establish the well-posedness of the present model, show
 the boundedness of the weak solution. We compare the experimental results with two recently developed models and arrive at the conclusion that the proposed model
well recovered the corrupted images without introducing undesired artifacts than that of existing models.  

\thispagestyle{empty}

\bibliographystyle{unsrt}

\end{document}